# GENERALIZED LÜROTH PROBLEMS, HIERARCHIZED I: SBNR - STABLY BIRATIONALIZED UNRAMIFIED SHEAVES AND LOWER RETRACT RATIONALITY

NORIHIKO MINAMI

Abstract. This is the first of a series of papers, where we investigate hierarchies of generalized Lüroth problems on the hierarchy of rationality, starting with the obvious hierarchy between the rationality and the ruledness.

Our primary goal here was to construct very general necessary conditions for a smooth, not necessary proper, scheme of finite type over a perfect base field $k$ to be "retract $(-i)$-rational,."

We achieve this goal by, for any Morel's unramified sheaf $S$, constructing **SBNR** - stably birationalized (Nisnevich) subsheaf $S_{sb}$ of the unramified sheaf $S$, in such a way that $S_{sb}$ coincides with $S$ for any proper smooth $k$-scheme of finite type.

Such a stably birationalized Nisnevich subsheaf $S_{sb}$ sheds a new light on the familiar irrational examples of Artin-Mumford, Saltman, Colliot-Thélène-Ojanguren, Bogomolov, Peyre, Colliot-Thélène-Voisin, and many other retract irrational classifying spaces of finite group, presented as counterexamples to the Noether problem of finite groups with the base complex number field $\mathbb{C}$. In fact, for all of these examples, the game is not over from our hierarchical perspective!

An immediate consequence of our construction of $S_{sb}$ is the stably birational invariance of an arbitrary unramified sheaf $S$ on proper smooth $k$-schemes of finite type. This in particular implies that, for any generalized motivic cohomology theory, its naively defined unramified (resp. stably birationalized) generalized motivic cohomology theory is stably birational invariant on smooth proper $k$-schemes of finite type (resp. arbitrary smooth $k$-schemes of finite type).

In the course of constructing $S_{sb}$, we have also shown a general local uniformization theorem of the first kind for arbitrary geometric valuations.

## 1. Introduction

In 1875, Lüroth [L875] asked the converse of the obvious impliation

$$\text{rational} \implies \text{unirational}$$

This is what is now called the **Lüroth problem**, and has drawn lots of attention, Nowadays, as is surveyed e.g. in [K96][B16][AB17] [P18][V19] , the Lüroth problem has been generalized to problems which ask the converse of various implications in the following refined implications:

(1)
$$\text{rational} \implies \text{stable rational} \implies \text{retract rational}$$
$$\implies \text{separably unirational} \implies \text{separably rationally connected} \implies \text{rationally connected}$$

These problems are what we mean by **generalized Lüroth problems**.

On the other hand, there is an obvious, though rarely featured, hierarchy which interpolates the rationality and the ruledness in the spirit of:

(2) $$\boxed{\text{lower rationality} = \text{higher ruledness}}$$

(For a precise definition, see Definition 1.1(i) below.)

In this series of papers, we embark on an organized study of the hierarchy of hierarchies of rationality, obtained by introducing an appropriate hierarchy in the spirit of (2) to each of (1). As a kick-off, we shall concentrate on the first three of (1). Fix the base field $k$, which we shall soon (in fact, after the following Definition 1.1) assume to be perfect. Then we have the following hierachies in the spirit of (2) for the first three of (1):

**Definition 1.1.** *For a $n$-dimensional $k$-variety* [1] *$X$, let us say:*

(i) $X$ *is* ***$(-i)$-rational*** *or* ***$(n-i)$-ruled*** $(0 \leq i \leq n)$

*if there exist an $i$-dimensional smooth $k$-variety $Z^i$, where $Z^0$ is taken to be* Spec $k$ *for $i = 0$, and a birational map*
$$\mathbb{A}^{n-i} \times Z^i - - > X.$$

(ii) $X$ *is* ***stable $(-i)$-rational*** *or* ***stable $(n-i)$-ruled*** $(0 \leq i \leq n)$

*if there exist $N \in \mathbb{Z}_{\geq n}$, an $i$-dimensional smooth $k$-variety $Z^i$, where $Z^0$ is taken to be* Spec $k$ *for $i = 0$, and a birational map*
$$\mathbb{A}^{N-n} \times \mathbb{A}^{n-i} \times Z^i - - > \mathbb{A}^{N-n} \times X.$$

(iii) $X$ *is* ***retract $(-i)$-rational*** *or* ***retract $(n-i)$-ruled*** $(0 \leq i \leq n)$

*if there exist $N \in \mathbb{Z}_{\geq n}$, an $i$-dimensional smooth $k$-variety $Z^i$, where $Z^0$ is taken to be* Spec $k$ *for $i = 0$, and rational maps*
$$f : X - - > \mathbb{A}^{N-i} \times Z^i, \quad g : \mathbb{A}^{N-i} \times Z^i - - > X$$
*such that the composition*
$$g \circ f : X - - > X$$
*is defined to be an identity on a dense open subset of $X$.*

**Remark 1.2.** We may also consider a slightly restrictive variant of the above definitions by demanding $Z^i$ $(i \in \mathbb{Z}_{\geq 1})$ to be further projective (of course, this is the same as simply demanding properness by Chow's lemma [sp18, 0200]), When *char* $k = 0$, definitions are unchanged by this extra requirement, thanks to Hironaka [H64]. However, it is not clear for the case *char* $k > 0$, and we shall mostly adapt this more restrictive definitions for the hierarchies in the sequel, where we are more concerned with smooth projective varieties.

[1] By a variety, we mean an integral $k$-scheme, which is separatred and of finite type as in [sp18, Tag 020D].

From now on, we assume that **the base field $k$ is perfect**. Then, the concepts presented in the above Definition 1.1(ii)(iii) are invariant with respect to the following standard equivalence relation (see e.g. [CTS07, §1]), as we shall show in Proposition 1.4(iv) below:

**Definition 1.3.** *Two varieties of possibly different dimensions $X$ and $Y$ are said to be* ***stable birational equivalent*** *if for some natural numbers $r, s$, $X \times \mathbb{A}^r$ and $Y \times \mathbb{A}^s$ are birationally equivalent.*

In fact, these hierarchies for the first three of (1) enojy satisfactory properties:

**Proposition 1.4.** (i) *When $i = 0$, those concepts presented in Definition 1.1 reduce to the usual classical concepts (mentioned in the first line of* (1)*).*

*0-rational = rational; stable 0-rational = stable rational;*
*retract 0-rational = retract rational.*

(ii) *Each concept in Definition 1.1 is a hierarchy; i.e. for any $0 \le i \le j \le n$,*

$(-i)$*-rational* $\Longrightarrow$ $(-j)$*-rational;*
*stable* $(-i)$*-rational* $\Longrightarrow$ *stable* $(-j)$*-rational;*
*retract* $(-i)$*-rational* $\Longrightarrow$ *retract* $(-j)$*-rational.*

(iii) *Concepts in Definition 1.1 define a hierarchy of hierarchies stated in above (ii); i.e. for any $0 \le i \le n$,*

$(-i)$*-rational* $\Longrightarrow$ *stable* $(-i)$*-rational* $\Longrightarrow$ *retract* $(-i)$*-rational*

(iv) *For any $0 \le i \le n$, stable $(-i)$-rationality and retract $(-i)$-rationality are stable birational invariants in the sense of Definition 1.3.*

**Remark 1.5.** (i) By definition, $X$ is stable rational if and only if $X$ is stable birationally equivalent to the point Spec $k$.
(ii) By definition, the concept of the stable birational equivalence is canonically extended from the usual varieties to ind-$k$-varieties $\{W_m\}_m$ with each structure morphism $i_m : W_m \to W_{m+1}$ equipped with a retraction $r_{m+1} : W_{m+1} \to W_m$ so that $r_{m+1} \circ i_m = id_{W_m}$, making $k(W_{m+1})$ purely transcendental over $k(W_m)$. For instance, for a linear algebraic group $G$ and a scheme $X$ with $G$-action, enjoying one of the conditions in [EG98, Proposition 23], the geometric classifying spaces $BG$ of [T99, §1] [MV99, §4] and the Borel construction $X_G = EG \times_G X$ (see e.g. [K18][CJ19, 2.1]) are such ind-$k$-varieties which possess respective canonical stable birational types, independent of particular ind-$k$-variety representations, thanks to Bogomolov's double fibration argument. (This stable birational independence can be argued and stated purely in terms of the corresponding pro-$k$-algebras. Bogomolov's double fibration argument is a culmination of many works [S919][EM73][L74] [BK85][CGR06, Lemma 4.4] [CTS07, 3.2], and christened "no-name lemma" by [D87].) This fact, amongst of all, enables us to consider the classical **Noether's problem** [2] of a finite group $G$ in terms of its classifying space $BG$, either by its geometric ind-$k$-variety model or its algebraic pro-$k$-algebra model. (e.g. [L05, Prop.9.4.4][M17, 4.2]).

**Corollary 1.6.** *For $BG$, and the Borel construction $X_G = EG \times_G X$, the concepts of stable $(-i)$-rationality and retract $(-i)$-rationality are well-defined stable birational invariants.* □

[2]In fact, the retract rationality, which is the case $i = 0$ of Definition 1.1 (iii), was introduced by Saltman [S84a, S84b] in his study of Noether's problem (for Noether's problem, consult surveys of [S83][CTS07][H14][H20], or books of [JLY02][GMS03][L05]).

For our later proof of Proposition 1.4 and especially to understand retract $(-i)$-rationality better, we now state some remarks:

**Remark 1.7.** (i) Generalizing the setting of the retract $(-i)$-rationality in Definition 1.1 (iii), let us consider the setting of the **rational retraction**. In other words, let us consider the situation when $k$-varieties $X, Y$, with rational maps

$$f : X - - > Y, \quad g : Y - - > X,$$

are given so that the composition

$$g \circ f : X - - > X$$

is defined to be an identity on a dense open subset of $X$. This is the same as saying that there are some dense opens $U \subseteq U' \subseteq X$, $V \subseteq Y$ such that $f, g$ restrict to morphisms

$$f : U \to V, \quad g : V \to U',$$

whose composition

$$g \circ f : U \to U'$$

is the canonical inclusion in $X$.

However, we sometimes wish to restrict our attention to a particular dense open $\widetilde{U} \subseteq U (\subseteq U' \subseteq X)$. Fortunately, we immediately get the following retraction given by morphisms, satisfying our desire:

$$\widetilde{U} \xrightarrow[f]{} g^{-1}(\widetilde{U}) \xrightarrow[g]{} \widetilde{U}, \qquad g \circ f = id_{\widetilde{U}}$$

(ii) For any $k$-variety $X$ and any $a \in \mathbb{A}^d$, we have an obvious retraction:

$$\begin{array}{ccccc} X & \overset{i_a}{\hookrightarrow} & \mathbb{A}^d \times X & \overset{p_X}{\twoheadrightarrow} & X \\ x & \mapsto & (a, x) & \mapsto & x \end{array}$$

However, we sometimes wish to restrict our attention to a particular dense open $V \subseteq \mathbb{A}^d \times X$. Fortunately, we can still choose appropriate $a \in \mathbb{A}^d$ and a non empty open (consequently dense, because $X$ is a $k$-variety) $U \subseteq X$, $V' \subseteq V \subseteq \mathbb{A}^d \times X$ with a restricted retraction given by honest morphisms of the form:

$$\begin{array}{ccccc} U & \overset{i}{\hookrightarrow} & V' & \overset{r}{\twoheadrightarrow} & U \\ \downarrow & & \downarrow & & \downarrow \\ X & \underset{x \mapsto (a,x)}{\overset{i_a}{\hookrightarrow}} & \mathbb{A}^d \times X & \overset{p_X}{\twoheadrightarrow} & X, \end{array}$$

which we think as a particular instance of rational retraction.

To see this, just choose any $a \in \mathbb{A}^r$ with $i_a^{-1}(V) \cong V \cap (\{a\} \times X) \neq \emptyset$. Then, just proceed as in (i) by setting:

$$V' := V \cap p_X^{-1}\left(i_a^{-1}(V)\right) = V \cap \left(\mathbb{A}^r \times i_a^{-1}(V)\right)$$
$$U := i_a^{-1}(V)$$

Those dense open subschemes we wish to restrict our attention in the above Remark 1.7 are actually some open subschemes of the smooth locus $U$ of the given variety $X$. For this purpose, the smooth locus $U$ should be non empty and dense open, which is guaranteed if $X$ is geometrically reduced. But this geometrically reduced condition is automatically satisfied under our perfect base field $k$ assumption. For this, see [BLR90, §2.2, Proposition 16] [G65, p.68, Proposition 4.6.1] [P17, Proposition 3.5.64, Warning 3.5.18, Prop.2.2.20]. We record the situation below following [sp18, 0B8X] for our later purpose:

**Proposition 1.8.** *Let $k$ be a perfect field, i.e.* [3] *every field extension of $k$ is separable over $k$. Let $X$ be a locally algebraic, i.e.* [4] *locally of finite type $k$-scheme as in* (12)*, which is furthermore reduced, for example a variety over $k$. Then we have*

$$\{x \in X \mid X \to \operatorname{Spec}(k) \textit{ is smooth at } x\} = \{x \in X \mid \mathcal{O}_{X,x} \textit{ is regular}\},$$

*which we call the* <u>*smooth locus*</u> *of $X$, is a dense open subscheme of $X$.* □

Motivated by Propositin 1.4 and Noether's problem, where $BG$ is approximated by non proper smooth varieties and the fact that the retract rationality implies the higher vanishings of the unramified cohomologies was so successfully exploited to produce many counter examples of Noether's problem for the case $k = \mathbb{C}$ [S84b] [B89] [P93] [CTS07], to cite just a few (for more, consult Hoshi's surveys [H14][H20] and [HKY20]) we are naturally led to the following problem:

The original motivation of this paper

Extend the following implications of the hierarchies in Propositin 1.4 for not necessarily proper varieties to some hierarchy ? of algebraic, i.e. "cohomological", stably birational invariants:

$$\{(-i)\text{-rational}\}_{i\in\mathbb{Z}_{\geq 0}} \implies \{\text{stable } (-i)\text{-rational}\}_{i\in\mathbb{Z}_{\geq 0}} \implies \{\text{retract } (-i)\text{-rational}\}_{i\in\mathbb{Z}_{\geq 0}} \implies \boxed{?} \tag{3}$$

To search after ?, let us suppose a variety $X$ over a perfect field $k$ is retract $(-i)$-rational. Then, for some $Z^i \in Sm_k^{ft}$, we have the corresponding rational retraction as in Remark 1.7(i) with $Y = \mathbb{A}^{N-i} \times Z^i$. and $\widetilde{U} \subseteq U (\subseteq U' \subseteq X)$ the dense smooth locus of $U$, as is guaranteed by Proposition 1.8.

Thus, denoting $\widetilde{U}$ simply by $U$ for ease of notation, we have the following diagram of $k$-varieties:

(4)
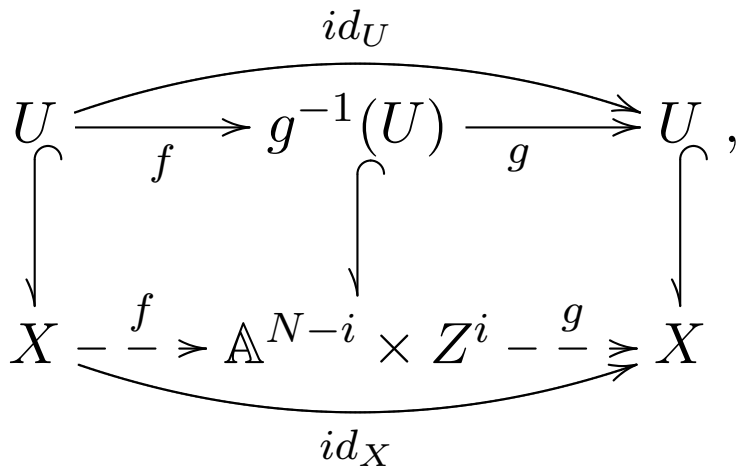

where vertical arrows are smooth dense open inclusions.

Then the following concept of Asok-Morel [AM11] is clearly of fundamental importance for our purpose:

[3] See [sp18, 030Y,030Z].
[4] See [sp18, 06LG].

**Definition 1.9.** ([AM11, Definition 6.1.1]) *Let $Sm_k^{ft}$ be the category of smooth* $\boxed{\textit{finite type}}$ *$k$-schemes. Then a presheaf*

$$M : \left(Sm_k^{ft}\right)^{op} \to \mathcal{C}(= Sets, Groups, ...),$$

*is called* birational *if the following properties are satisfied:*

**(B0):** *For any $X \in Sm_k$, denote by $X^{(0)}$ its set of irreducible components so that $X = \cup_{\alpha \in X^{(0)}} X_\alpha$. then the obvious map*

$$S(X) \to \prod_{\alpha \in X^{(0)}} S(X_\alpha)$$

*is a bijection;*

**(B1):** *For any $X \in Sm_k$ and its arbitrary everywhere dense open subscheme $U \subset X$, the restriction map*

$$S(X) \to S(U)$$

*is a bijection.*

**Remark 1.10.** (i) ([AM11, Lemm.6.1.2] [MV99, §3, Proposition 1.4]) Any birational presheaf is automatically a Nisnevich sheaf.
(ii) ([KS15, p.330, Appendix A by Jean-Louis Colliot-Thélène]) Any birational presheaf $S$ is automatically stably birational invariant, thanks to a beautiful observation of Colliot-Thélène [KS15, p.330, Appendix A by Jean-Louis Colliot-Thélène], whose proof can be succinctly summarized as follows:

Consider the birational correspondence from $\mathbb{P}^1 \times \mathbb{P}^1$ to $\mathbb{P}^2$ given by the surface expressed in terms of the multihomogeneous coordinates:

$$W := V(uT - vX, wT - zY) \subset \{(u, v; w, z; X, Y, T)\} = \mathbb{P}^1 \times \mathbb{P}^1 \times \mathbb{P}^2.$$

This realizes the composition of the blowup of $\mathbb{P}^1 \times \mathbb{P}^1$ at the $k$-point

$$M := (1, 0; 1; 0) \in \mathbb{P}^1 \times \mathbb{P}^1$$

and the contraction of the proper transforms

$$E_1 := \overline{\{(u, v; 1, 0; X, Y, 0) \mid vX = 0\} \setminus \{(1, 0; 1, 0; X, Y, 0)\}} = \{(u, v; 1, 0; 0, Y, 0)\} \simeq \mathbb{P}^1;$$

$$E_2 := \overline{\{(1, 0; w, z; X, Y, 0) \mid zY = 0\} \setminus \{(1, 0; 1, 0; X, Y, 0)\}} = \{(1, 0; w, z; X, 0, 0)\} \simeq \mathbb{P}^1,$$

which are non intersecting exceptional curves of the first kind, of the following rational curves in $\mathbb{P}^1 \times \mathbb{P}^1$ :

$$L_1 := \{(u, v; 1, 0)\} = \mathbb{P}^1 \times \{(1, 0)\} \subset \mathbb{P}^1 \times \mathbb{P}^1,$$
$$L_2 := \{(1, 0; w, z)\} = \{(1, 0)\} \times \mathbb{P}^1 \subset \mathbb{P}^1 \times \mathbb{P}^1.$$

to the following $k$-points of $\mathbb{P}^2$, respectively:

$$M_1 := \{(0, 1, 0)\} \subset \mathbb{P}^2,$$
$$M_2 := \{(1, 0, 0)\} \subset \mathbb{P}^2.$$

Of course, this conversely implies $W$ is the blowup of $\mathbb{P}^2$ at $M_1$ and $M_2$.

Now the stably birational invariance of a birational presheaf $S$ immediately follows from the following commutative diagram (The point is to observe the surjectivity of $S(W) \to S(E_1)$ from the left hand side of the commutative diagram, and apply it in the right hand side of the commutative diagram to deduce the surjectivity of $S(\mathrm{Spec}(k)) \cong S(M_1) \to S(E_1) \cong S(\mathbb{P}^1)$, exploiting the birational invariance of $S$: $S\left(\mathbb{P}^2\right) \xrightarrow{\simeq} S(W)$. ):

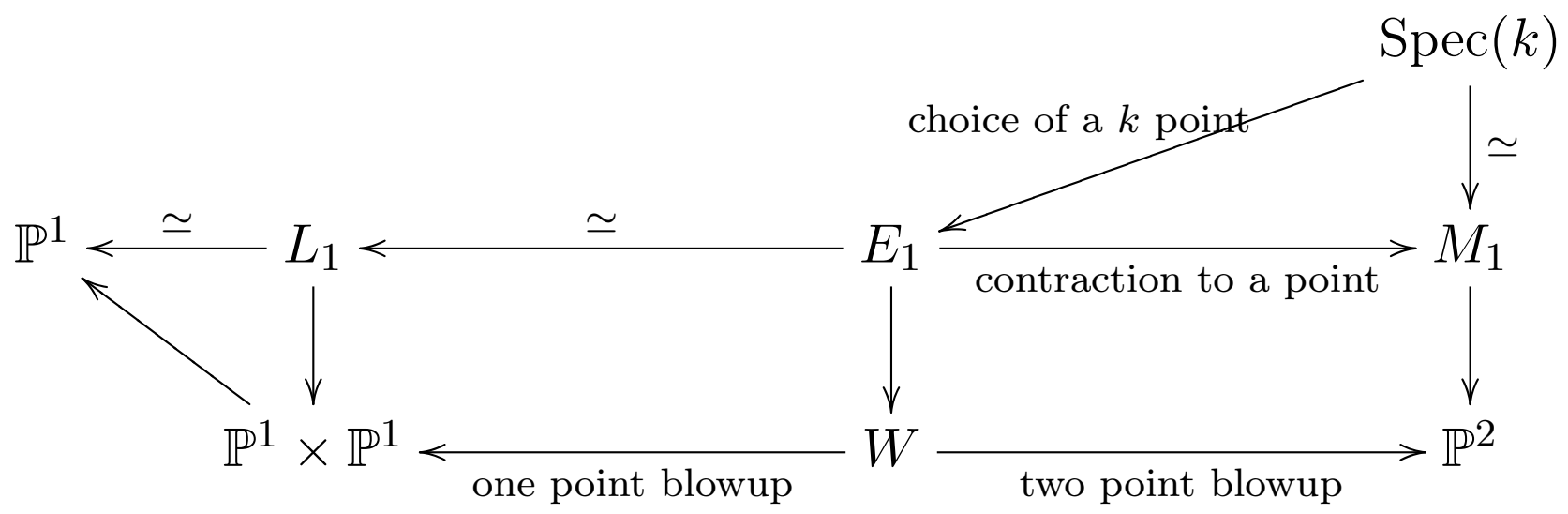

Applying a birational Nisnevich sheaf $S$ to the commutative diagram (4), the above Remark 1.10(ii) tells us the following important conclusion:

**Theorem 1.11.** *If $X \in Sm_k^{ft}$ is retract $(-i)$-rational, then $S(X)$ is a direct summand of $S\left(Z^i\right)$ for some $Z^i \in Sm_k^{ft}$ of dimension $i$.*

Then, in search of useful algebraic, i.e. "cohomological", stable biratoinal invariants [?] for the purpose of (3), we turn our attention to unramified sheaves in the sense of Morel [M12, Chapter 2]. While we shall delay its precise definition until Definition 4.1, we emphasize that there are considerably many examples of unramified sheaves, including the unramified cohomology, as we recall in Example 4.4.

**SBNR** - Stably birationalized Nisnevich subsheaf $S_{sb}$ of an unramified Zariski-sheaf $S$

**Theorem and Definition 1.12.** (Theorem and Definition 4.5)
(i) Given an unramified Zariski sheaf $S$ on $Sm_k$, let us set

$$S_{sb}(K/k) := \cap_{\substack{v,\ \text{divisorial} \\ \text{valuation of } K/k}} S(V_v) \tag{5}$$

for any finitely generated field extension $K/k$. Then the correspondence

$$U \mapsto S_{sb}(U) := S_{sb}(k(U)/k) = \cap_{\substack{v,\ \text{divisorial} \\ \text{valuation of } k(U)/k}} S(V_v)$$

$$\left( \subseteq S(U) := \cap_{\substack{v,\ \text{divisorial valuation of } k(U)/k \\ \text{s.t. } V_v = \mathcal{O}_{U,x} \text{ for some } x \in U^{(1)}}} S(V_v) = \cap_{x \in U^{(1)}} S(\mathcal{O}_{U,x}) \subset S(k(U)) \right)$$

defines a birational subsheaf $S_{sb}$ of an unramified Zariski sheaf $S$ on $Sm_k^{ft}$ which we call the *stably birationalized subsheaf* $S_{sb}$ of an unramfied sheaf $S$, or simply **SBNR**, because $S_{ab}$ is, once its presheaf $\left(Sm_k^{ft}\right)^{op} \to \mathcal{C}$ property is established, immediately seen to be a birational presheaf, which is actually a Nisnevich sheaf by [AM11, Lemm.6.1.2] [MV99, §3, Proposition 1.4] as was recalled in Remark 1.10 (i), and is staby birational invariant by [KS15, p.330, Appendix A by Jean-Louis Colliot-Thélène] as was recalled in Remark 1.10 (ii).
(ii) For any smooth proper $k$-scheme $X$,

$$S_{sb}(X) = S(X).$$

Consequently, any unramified sheaf $S$ is stably birational invariant among smooth proper $k$-schemes.

While we shall defer the precise definition of "divisorial valuation" and underlying the basic valuation theory until the next section, we here give a very short summary of necessary notations and terminologies, including "$V_v$ for a divisorial valuation", which showed up in the above Theorem and Definition 1.12:

- For a valuation $v$ for a finitely generated field extension $K/k$, denote by $V_v, \mathfrak{m}_v, \Gamma_v$ its valuation ring, the maximal ideal, and the value group, respectively, so that

$$\begin{cases} v : K^* \twoheadrightarrow \Gamma_v \\ V_v := v^{-1}\big( (\Gamma_v)_{\geq 0} \big) \\ \mathfrak{m}_v = v^{-1}\big( (\Gamma_v)_{>0} \big) \\ k \subseteq v^{-1}(0) \end{cases}$$

- Set $rank(v) := \dim\ V_v$, then we have the following inequality:

$$rank(v) \ + \ tr.deg_k\, V_v/\mathfrak{m}_v \ \leq \ tr.deg_k\, K. \tag{6}$$

- If an equality holds in the above (6), then the valuation $v$ is called *geometric*. A rank 1 geometric valuation $v$ is also called a *divisorial valuation*.
- The *local uniformization problem* of an algebraic variety $X$ over $k$ asks, for each valuation on $k(X)/k$, existence of a finitely generated $k$-subalgebra $R$ $(\subseteq V_v \subset k(X))$ with $Frac(R) = k(X)$ such that the following is regular:

$$R_{\mathfrak{m}_v \cap R} \left( \subseteq V_v \subset k(X) = Frac(R) \right)$$

Now, the essence of our proof of Theorem and Definition 1.12 is a very special kind of, i.e. *first kind*, local uniformization theorem for general geometric valuation, whose summary is recalled below:

A condensed version of Theorem 3.6

**Theorem 1.13.** *For any rank $r$ geometric valuation $v = v_1 \circ \cdots \circ v_r$ on a finitely generated field extensin $K/k$ and a fnitely generated $k$ algebra $S$ such that*

$$S \subset V_{v_1 \circ \cdots \circ v_r} \subset K = Frac(S)$$

*there exists a* *local blowup of $S$ with respect to $v_1 \circ \cdots \circ v_r$, i.e. an obviously defined local ring homomorphisms*

$$\begin{aligned}&(S_{\mathfrak{m}_{v_1 \circ \cdots \circ v_r} \cap S}, \mathfrak{m}_{v_1 \circ \cdots \circ v_r} \cap S_{\mathfrak{m}_{v_1 \circ \cdots \circ v_r} \cap S}) \\ \to &\left( S\left[\frac{a_1}{b_1}, \ldots, \frac{a_r}{b_r}\right]_{\mathfrak{m}_v \cap S\left[\frac{a_1}{b_1}, \ldots, \frac{a_r}{b_r}\right]}, \mathfrak{m}_{v_1 \circ \cdots \circ v_r} \cap S\left[\frac{a_1}{b_1}, \ldots, \frac{a_r}{b_r}\right]_{\mathfrak{m}_v \cap S\left[\frac{a_1}{b_1}, \ldots, \frac{a_r}{b_r}\right]} \right) \\ \to &(V_{v_1 \circ \cdots \circ v_r}, \mathfrak{m}_{v_1 \circ \cdots \circ v_r})\end{aligned}$$

*for some $a_i, b_i \in S$ $(i = 1, \ldots, r)$ with*

$$(v_1 \circ \cdots \circ v_r)(a_i) \;\geq\; (v_1 \circ \cdots \circ v_r)(b_i) \quad (i = 1, \ldots, r),$$

*such that the following conditions are satisfied:*

(a) *There exists a finitely generated $k$ subalgebra $R_f$*

$$R_f \;\subset\; S\left[\frac{a_1}{b_1}, \ldots, \frac{a_r}{b_r}\right]_{\mathfrak{m}_v \cap S\left[\frac{a_1}{b_1}, \ldots, \frac{a_r}{b_r}\right]} \;\subset\; V_{v_1 \circ \cdots \circ v_r} \tag{7}$$

*such that $\mathfrak{m}_{v_1 \circ \cdots \circ v_r} \cap R_f \subset R_f$ is a height $r$ prime ideal.*

(b) *For $R = R_f$, $S\left[\frac{a_1}{b_1}, \ldots, \frac{a_r}{b_r}\right]_{\mathfrak{m}_v \cap S\left[\frac{a_1}{b_1}, \ldots, \frac{a_r}{b_r}\right]}$,*

- *For each $1 \leq i \leq j \leq r$, the local homomorphism induced by* (7)

$$\left(\frac{R}{\mathfrak{m}_{v_1 \circ \cdots \circ v_{i-1}} \cap R}\right)_{\overline{\mathfrak{m}_{v_1 \circ \cdots \circ v_j} \cap R}} \to V_{v_i \circ \cdots \circ v_j}$$

*has its source $\left(\frac{R}{\mathfrak{m}_{v_1 \circ \cdots \circ v_{i-1}} \cap R}\right)_{\overline{\mathfrak{m}_{v_1 \circ \cdots \circ v_j} \cap R}}$ regular local (actually essentially smooth).*

- *For each $1 \leq i \leq r$, the following first kind condition is satisfied:*

$$\begin{array}{ccccc} Frac\left(\frac{R}{\mathfrak{m}_{v_1 \circ \cdots \circ v_{i-1}} \cap R}\right) & \hookleftarrow & \left(\frac{R}{\mathfrak{m}_{v_1 \circ \cdots \circ v_{i-1}} \cap R}\right)_{\overline{\mathfrak{m}_{v_1 \circ \cdots \circ v_i} \cap R}} & \twoheadrightarrow & Frac\left(\frac{R}{\mathfrak{m}_{v_1 \circ \cdots \circ v_i} \cap R}\right) \\ \| & & \text{"first kind" condition} \; \| & & \| \\ \kappa(v_{i-1}) & \hookleftarrow & V_{v_i} & \twoheadrightarrow & \kappa(v_i) \end{array}$$

This theorem generalizes a classical theorem of Zariski [Z39], which corresponds to the case $rank(v) = 1$, to arbitrary rank geometric valuations, and should be of independent interest. For the precise statement, including notations and backgrounds, please consult §3, especially Theorem 3.6.

Going back to our main Theorem and Definition 1.12, we have:

**Remark 1.14.** (i) I believe the above Theorem and Definition 1.12, which is applicable to arbitrary smooth, not necessary proper, $k$-scheme, is the most conceptually transparent result along the line of [CT95] (whose philosophy might go back to [G68] which deals with Brauer Groups).
(ii) When the unramified sheaf $S$ is provided by a Rost cycle module $M$, i.e. when

$S(X) = A^0(X, M)$, $S_{sb}$ was already defined by Merkurjev [M08, p.56, 2.2] and an intimately related claim [M08, p.61, Proposition 2.15] was already stated [M08, p.56, 2.2; p.61, Proposition 2.15]. [5]

This was also considered by Kahn-Ngan [KN14, Definition 6.1], but be aware that different notations are used for this case:

$$\underbrace{A^0_{sb}(-, M)}_{\text{our notation}} = \underbrace{M(-)_{nr}}_{\text{Merkurjev's notation}} = \underbrace{A^0_{nr}(-, M)}_{\text{Kahn-Ngan's notation}}$$

We opted not to adopt the notation "$S_{nr}$" for our $S_{sb}$ because $S$ is already an unramified sheaf. "$S_{NR}$" is not appropriate either, because of a potential confusion with [KN16, p.274, Definition 3.1].

(iii) Stably birational invariance of an unramified sheaf among smooth proper $k$-schemes, claimed in Theorem and Definition 1.11 (ii), was proved for Rost's cycle modules mentioned in Example 4.4 (ii) by Rost [R96, Corollary (12.10)], and for Feld's Milnor-Witt cycle modules mentioned in Example 4.4 (iii) by Feld [F21a, Theorem 5.3.1]. However, from these approaches which stick only to smooth and proper $k$-schemes, we can never deduce Corollary 1.15 which is valid for arbitrary smooth $k$-schemes.

Now, in view of Theorem 1.11, we immediately obtain the following;

**Corollary 1.15.** (Corollary 4.6) *Let $S$ be an unramified Zariski sheaf, and $S_{sb}$ be its associated stably birational Nisnevich sheaf. Then, for any retract $(-i)$-rational $X$ in $Sm^{ft}_k$, or more generally in the category of* ind $-k$ *varieties which includes $BG$ and $Y_G = EG \times_G Y$ as in Remark 1.5(ii), the following hold:*

**The case $i = 0$:** *The following canonical map is an isomorphism:*

$$S(\text{Spec } k) = S_{sb}(\text{Spec } k) \xrightarrow{\cong} S_{sb}(X).$$

**The case $i > 0$:** *For some smooth $Z^i$ of dimension $i$. $S_{sb}(X)$ is a direct summand of $S_{sb}(Z^i)$.*

□

In this paper, we are mostly concerned with the case $i > 0$, and we wish to have some kind of vanishing theorems of

$$S_{sb}(Z^i) \quad \left( \subseteq S_{nr}(Z^i) \subseteq S\left(k(Z^i)\right) \right)$$

to exploit Corollary 1.15. Fortunately, there are some:

- By choosing $r = p^m, D = \mu_{p^m}^{l-j}$ with $(char\ k, p) = 1$ in Example 4.4(ii), we find the classical unramified cohomology

$$Sm^{ft}_k \ni X \quad \mapsto \quad H^j_{nr}\left(X/k, \mu_{p^m}^{\otimes l}\right)$$

[5]After this paper was written up, the author realized Kahn-Sujatha [KS17, §7] also proved very similar results when the unramified sheaf $S$ is provided by a homotopy invariant Nisnevich sheaf with transfers.

is unramified. So, we may consider its associated stably birationalised Nisnevich sheaf with $X = Z^i$, which is clearly contained in the Galois group:[6]

$$\left(H^j_{nr}\left(k(Z^i)/k, \mu_{p^m}^{\otimes l}\right)\right)_{sb} \subset H^j_{nr}\left(Z^i/k, \mu_{p^m}^{\otimes l}\right) \subset H^j_{Gal}\left(k(Z^i)/k, \mu_{p^m}^{\otimes l}\right)$$

Here, by the classical Tate's thoerem (formally a conjecture of Grothendieck) [S72, p.119, §4, Th. 28] (see also [AGV73, Exposé xiv 3]) about the cohomological dimension, the rightmost galois group vanishes when $j > i + cd_p(k)$. In particular, we obtain the following vanishing theorem:

When $(char\ k, p) = 1$, the following vanishing theorem holds for any $j > i + cd_p(k)$ and $1 \leq m \leq +\infty, l \in \mathbb{Z}$:

$$\left(H^j_{nr}\left(k(Z^i)/k, \mu_{p^m}^{\otimes l}\right)\right)_{sb} = H^j_{nr}\left(Z^i/k, \mu_{p^m}^{\otimes l}\right) = 0.$$

In particular, when $k = \mathbb{C}$, this vanishing theorem holds for any $j > i$ and $1 \leq m \leq +\infty, l \in \mathbb{Z}$.

- As is recalled in Example 4.4(iv), the motivic (stable) homotopy theory offers an unlimited supply of unramified sheaves. In fact, it is shown in Example 4.4(iv) that, for any generalized motivic cohomology theory $E^s(-)_t$ represented by a $S^1$-spectrum and for any generalized motivic cohomology theory $\mathbb{M}^{p,q}(-)$ represented by a $\mathbb{P}^1$-spectrum, the naively defined their unramified analogues

(8)
$$\begin{aligned} Sm^{ft}_k &\to AbelianGroups \\ X &\mapsto E^s(X)_t,\ \mathbb{M}^{p,q}(X) \end{aligned}$$

are both unramified, and so define their stably birationalized Nisnevich subsheaves, together with the canonical maps from the original presheaves (8):

$$E^s_{sb}\left(k(X)\right)_t \hookrightarrow E^s_{nr}(X)_t \xleftarrow{\exists 1} E^s(X)_t$$
$$\mathbb{M}^{p,q}_{sb}\left(k(X)\right) \hookrightarrow \mathbb{M}^{p,q}_{nr}(X) \xleftarrow{\exists 1} \mathbb{M}^{p,q}(X)$$

Now, Voevodsky [V98, Theorem 4.14] and Morel [M05b, Lemma 4.3.1] proved the vanishing theorem which implies the following:

Suppose a generalized motivic cohomology theory $E^*(-)_t$ is represented by a connective $S^1$-spectrum, and a generalized motivic cohomology theory $\mathbb{M}^{*,*}(-)$ is represented by a connective $\mathbb{P}^1$-spectrum. Then the following vanishing theorems hold for any $s > i$ and $p - q > i$ :

$$E^s_{sb}\left(k(Z_i)/k\right)_t = E^s_{nr}\left(Z_i/k\right)_t = 0$$
$$\mathbb{M}^{p,q}_{sb}\left(k(Z_i)/k\right) = \mathbb{M}^{p,q}_{nr}\left(Z_i/k\right) = 0$$

- As is recalled in Example 4.4(v), the sheaf of Kähler differentials $\Omega^j$

$$Sm^{ft}_k \ni X \quad \mapsto \quad H^0(X, \Omega^j)$$

[6] We warn the readers that our $\left(H^j_{nr}\left(k(Z^i)/k, \mu_{p^m}^{\otimes l}\right)\right)_{sb}$ is simply denoted by $H^j_{nr}\left(k(Z^i)/k, \mu_{p^m}^{\otimes l}\right)$ and called unramified cohomology in some literature like [S21b], and the original literature [CTO89].

is unramified. Obviously, we have the following vanishing theorem for this:

$$H^0\left(Z^i, \Omega^j\right) = 0, \quad \forall j > i.$$

This naive vanishing theorem turns out to play a very important role when we shall upgrade Totaro's hypersurface irrationality theorem [T16] in our hierarchical setting, as is briefly previewed at the end of this introduction.

Then our Theorem and Definition 1.12 and the above vanishing theorems lead to the following very useful conclusion:

**Theorem 1.16.** *Let $p$ be a prime number such that $(char\ k, p) = 1$, $E^*(-)_t$ be a generalized motivic cohomology theory represented by a connective $S^1$-spectrum, and $\mathbb{M}^{*,*}(-)$ be a generalized motivic cohomology theory represented by a connective $\mathbb{P}^1$-spectrum. Suppose either one of the following non-vanishing conditions holds for $X \in Sm_k^{ft}$ :*

$$\begin{cases} \left(H^j_{nr}\left(k(X)/k, \mu_{p^m}^{\otimes l}\right)\right)_{sb} & \neq 0 \qquad (\exists j > i + cd_p(k)) \\ E^s_{sb}\left(k(X)/k\right)_t & \neq 0 \qquad (\exists s > i) \\ \mathbb{M}^{p,q}_{sb}\left(k(X)/k\right) & \neq 0 \qquad (\exists p > i + q) \end{cases}$$

*Then $X$ is not retract $(-i)$-rational.*
*Here, when $X$ is furhter $\boxed{\text{proper}}$, recall from Theorem and Definition 1.12 (ii) that either one of the above left hand side stably birationalized Nisnevich sheaf is equal to its unramified analogue.*

The most successful instances when the above Theorem 1.16 is applicable are found among the unramified cohomology theories with $k = \mathbb{C}$, when $\mu_{p^m}$ becomes the trivial Galois module:

**Corollary 1.17.** *When $k = \mathbb{C}$, if*

$$\left(H^j_{nr}\left(k(X)/k, \mu_{p^m}^{\otimes j}\right)\right)_{sb} \neq 0 \qquad (1 \leq \exists m \leq +\infty)$$

*for some*

$$j > i,$$

*then $X$ is not retract $(-i)$-rational.*
*Here, when $X$ is further $\boxed{\text{proper}}$, recall from Theorem and Definition 1.12 (ii) that*

$$\left(H^j_{nr}\left(k(X)/k, \mu_{p^m}^{\otimes j}\right)\right)_{sb} = H^j_{nr}\left(X_k, \mu_{p^m}^{\otimes j}\right).$$

Now, let us state three messages delivered by the philosophy of this paper.

1.1. **The game is not over!:** Starting with the pioneering work of Artin-Mumford [AM72], unramified cohomologies have been so successfully applied to produce many non retract rational examples (e.g. [S84b][CTO89][P93][CTS07], to quote just a few). Here, the following very important message is delivered by our Theorem 1.16 (and Corollary 1.17):

The game is not over!

Up until now, whenever a smooth variety $X$ is found to enjoy

(9) $$\left(H^j_{nr}\left(k(X)/k, \mu_{p^m}^{\otimes l}\right)\right)_{sb} \neq 0 \qquad (1 \leq \exists m \leq +\infty)$$

for some

$$j > i + cd_p(k),$$

it was immediately declared that the game is over simply because $X$ is now known to be non retract rational.

However, from our perspective, this is nothing but a prelude to the new game. In fact, under the assumption, we know much more that $X$ is not even retract $(-i)$-rational. And so, we are prompted to answer the next hierarchical level question:

*Is $X$ retract $(-(i+1))$-rational?*

*To start with, we should answer the following question:*

*Is* $\left(H^{j+1}_{nr}\left(k(X)/k, \mu_{p^m}^{\otimes l}\right)\right)_{sb} \overset{?}{\neq} 0 \ (1 \leq \exists m \leq +\infty)$

Here, when $X$ is further proper, recall from Theorem and Definition 1.12 (ii) that

$$\left(H^j_{nr}\left(k(X)/k, \mu_{p^m}^{\otimes l}\right)\right)_{sb} = H^j_{nr}\left(X_k, \mu_{p^m}^{\otimes l}\right).$$

Of course, there are so many such an example of $X$ such that $\left(H^j_{nr}\left(k(X)/k, \mu_{p^m}^{\otimes l}\right)\right)_{sb} \neq 0$ for some $j > 0$. For instance, when $k = \mathbb{C}$, Asok [A13, Theorem 4.2], building upon [P93], constructed for each $n \in \mathbb{N}$ a unirational smooth proper $X$ such that $H^n_{nr}\left(X_{\mathbb{C}}, \mu_n^{\otimes n}\right) \neq 0$. However, unlike our point of view, the next concern of Asok [A13, Theorem 4.2] was the vanishings of $H^j_{nr}\left(X_{\mathbb{C}}, \mu_n^{\otimes j}\right)$ for $j < n$. Of course, our next concern here is the vanishings of $H^j_{nr}\left(X_{\mathbb{C}}, \mu_n^{\otimes j}\right)$ for $j > n$. This is a testimony how our very naive looking hierarchical consideration is actually revolutionary; completely opposite to the traditional thinking.

Without the proper assumption, many finite group $G$ have been shown to yield a counter-examples to Noehter's problem, by proving $\left(H^j_{nr}\left(\mathbb{C}(BG)/\mathbb{C}, \mu_{p^m}^{\otimes j}\right)\right)_{sb} \neq 0$ for some small $j$, as are recalled in Hoshi's excellent surveys [H14][H20]. For all of these groups, the most natural problem from our hierarchical point of view, i.e. to find the largest $j$ such that $\left(H^j_{nr}\left(\mathbb{C}(BG)/\mathbb{C}, \mu_{p^m}^{\otimes j}\right)\right)_{sb} \neq 0$, is very much wide open!

1.2. **Counter examples to Integral Hodge conjecture:** For a complex projective manifold $X$, set $Z^{2k}(X)$ to be the cokernel of the cycle map *cycl*:

$$Z^{2k}(X) := \mathrm{Coker}\left(cycl : \mathrm{CH}^k(X) \to H^{2k}(X, \mathbb{Z}) \cap H^{k,k}(X)\right)$$

The the *Hodge conjecture* $HC^k(X)_{\mathbb{Q}}$ and the *integral Hodge conjecture* $HC^k(X)_{\mathbb{Z}}$ respectively predict the triviality of $Z^{2k}(X)_{\mathbb{Q}}$ and $Z^{2k}(X)$:

$$\begin{cases} HC^k(X)_{\mathbb{Q}} & \overset{\text{Definition}}{\Longleftrightarrow} Z^{2k}(X)_{\mathbb{Q}} = 0 \\ HC^k(X)_{\mathbb{Z}} & \overset{\text{Definition}}{\Longleftrightarrow} Z^{2k}(X) = 0 \end{cases}$$

Then, Atiyah-Hirzebruch [AH62] found a counter-examples to the integral Hodge conjecture $HC^2(X)_{\mathbb{Z}}$, and more examples were constructed later (e.g. [K92][T99] [CTV12], to quote just a few). For some time, there were some hope that such counter-examples to the integral Hodge conjecture $HC^k(X)_{\mathbb{Z}}$ might give some insight toward the original Hodge conjecture $HC^k(X)_{\mathbb{Q}}$.

However, such a hope was gone when Rosenschon-Srinivas [RS16] reformulated the integral Hodge conjecture étale theoretically $HC^k_L(X)_{\mathbb{Z}}$:

$$\begin{cases} HC^k_L(X)_{\mathbb{Z}} \overset{\text{Definition}}{\iff} Z^{2k}_L(X) = 0 \\ \text{where } Z^{2k}_L(X) \text{ is defined to be:} \\ \text{Coker}\left(\mathrm{CH}^k_L(X) := H^{2k}_L(X, \mathbb{Z}(n)) := \mathbb{H}^m_{\text{ét}}\left(X, z_n(-, \bullet)[-2k]\right) \to H^{2k}(X, \mathbb{Z}) \cap H^{k,k}(X)\right) \end{cases}$$

enjoying

$$HC^k(X)_{\mathbb{Q}} \iff HC^k_L(X)_{\mathbb{Z}}.$$

In other words, counter-examples to the integral Hodge conjecture $HC^k(X)_{\mathbb{Z}}$ do not give any insight toward the original Hodge conjecture $HC^k(X)_{\mathbb{Q}}$. (For the precise meaning of these notations, please consult the original paper [RS16], where a similar result for the Tate conjecture is also obtained.)

Then we naturally wonder if there is any implication of counter-examples to the integral Hodge conjecture or not. Now, Colliot-Thélène-Voisin [CTV12] proved the following theorem, building upon the solution of the Bloch-Kato conjecture by Voevodsky [V03][V11] Rost [R98] [R02] (see also [SJ06][W09][HW19]):

**Theorem 1.18.** [CTV12, Theorem 3.7] *For a smooth projective connected $\mathbb{C}$-manifold, we have the following short exact sequence:*

$$0 \to H^3_{nr}\left(X_{\mathbb{C}}, \mathbb{Z}(2)\right) \otimes \mathbb{Q}/\mathbb{Z} \to H^3_{nr}\left(X_{\mathbb{C}}, \mathbb{Q}/\mathbb{Z}(2)\right) \to Z^4(X)\{tors\} \to 0.$$

From Corollary 1.17, the above Theorem 1.18 immediately furnishes the following implication for many integral Hodge conjecture counter-examples:

**Corollary 1.19.** *Let $X$ be a smooth projective connected $\mathbb{C}$-manifold $X$, which yields an integral Hodge conjecture counter-example such that $Z^4(X)\{tors\} \neq 0$. Then $X$ is not retract $(-2)$-rational.*

We shall say more on this topic from the viewpoint of the rational connectedness and the rational $(-2)$-connectedness discussed in the sequel, which is briefly previewed at the end of this introduction.

1.3. **Homotopy group sheaf:** Recently, Isaksen-Wang-Xu [IWX20] made a spectacular advance in the computation of the two primary stable homotopy groups of the spheres from the stem 61 to the stem 90. In their calculation, the *motivic stable homotopy groups of the spheres for $\mathbb{C}$* plays some important role. Now, the motivic stable homotopy groups of the spheres $\pi^{\mathbb{A}^1}_{p,q}(\mathbb{S})$, where $\mathbb{S}$ is the motivic sphere spectrum, is in fact an unramified sheaf, obtained as the Nisnevich sheafication of the presheaf:

$$\begin{aligned} \pi_{p,q}(\mathbb{S}) : Sm^{ft}_k &\to \text{Abelian Groups} \\ U &\mapsto Hom_{\mathcal{SH}(k)}\left(\Sigma^{p-q}_{S^1}\Sigma^q_{\mathbb{G}_m}(U_+), \mathbb{S}\right), \end{aligned}$$

as we shall see in Example 4.4(iv). And what Isaksen-Wang-Xu [IWX20] exploited are just their values at the base field for the case $k = \mathbb{C} : \pi_{p,q}^{\mathbb{A}^1}(\mathbb{S})(\mathrm{Spec}\ \mathbb{C})$. Then we are naturally led to a suspicion that what homotopy theorists are enjoying from algebraic geometry is just information from the point Spec $\mathbb{C}$.

Now, our Theorem and Definition 1.12(ii) and Example 4.4(iv) impliy $\pi_{p,q}^{\mathbb{A}^1}(\mathbb{S})(\mathrm{Spec}\ \mathbb{C})$ are actually reflecting all the smooth proper $\mathbb{C}$-schemes:

$$\pi_{p,q}^{\mathbb{A}^1}(\mathbb{S})(\mathrm{Spec}\ \mathbb{C}) \ = \ \pi_{p,q}^{\mathbb{A}^1}(\mathbb{S})(X), \quad (X \text{ is any smooth proper } \mathbb{C}\text{-scheme.})$$

In fact, motivic (stable) homotopy theory offers an unlimited supply of stable birational invariants of smooth proper $k$-schemes, as we see in Example 4.4(iv).

The rest of this paper is organized as follows:

Most of the material of §2 and a good portion of §3 should be standard for experts of algebraic geometry. However, the author opted to leave them as they are presented in this paper, because many of the expected readers of this paper would be algebraic number theorists and algebraic topologists, who are not so familiar with algebraic geometry.

Now, since this is the first of a series of papers, let us close this introduction by giving a short previews of the sequel:

1.4. **<u>A preview of the sequel:</u>** In this paper, we have focused upon the hierarchies (Definition 1.1) for the first three of (1) for general smooth schemes without properness assumption. In the sequel, we shall introduce the hierarchies, extending Definition 1.1, for the latter three of (1). and study their properties and interrelationship for smooth $\boxed{\text{proper}}$ $k$-varieties. Thus, we shall establish the concepts of

$$\begin{aligned} \cdots \implies \text{retract } (-i)\text{-rational} \implies \text{separably } (-i)\text{-unirational} \\ \implies \text{separably } (-i)\text{-rationally connected} \implies \text{rationally } (-i)\text{-connected} \end{aligned}$$

When $k = \mathbb{C}$, we shall give a criterion of rationally $(-i)$-connectedness, by interpolating the uniruledness criterion of [BDPP13] and the rationally connectedness criterion of [CDP15]. We shall also incorporate the hierarchical generalization of the classical decomposition of the diagonal [BS83] studied in [CTV12] [V14][CL17], together with [KS16][KOY21], to establish the following implication diagram for a smooth proper variety $X$ (where we have restricted to the case *char* $k = 0$ for simplicity):

$$\begin{array}{ccccc} \text{retract } (-i)\text{-rational} & \Longrightarrow & \exists\text{level } i \text{ integral decomposition of diagonal} & \Longrightarrow & \begin{array}{c}\text{analogue of Corollary 1.15} \\ \text{for } \mathbb{P}^1\text{-invariant Nisnevich} \\ \text{sheaves with transfer}\end{array} \\ \Downarrow & & \Downarrow & & \Downarrow \text{if } k = \mathbb{C} \\ \text{rationally } (-i)\text{-connected} & \Longrightarrow & \exists\text{level } i \text{ rational decomposition of diagonal} & \overset{\text{if } k = \mathbb{C}}{\Longrightarrow} & H^j_{nr}(X, \mathbb{Z}(q)) = 0, \ \forall j > i \end{array}$$

The reason we are interested in the level $i$ decomposition of diagonal is because, in the specialization setting:

$$\begin{array}{ccccccc}
Z & \xrightarrow[\substack{\text{birational} \\ \text{total.}CH_0\text{-trivial}}]{q} & Y \left(\substack{\text{pssibly} \\ \text{singular}}\right) & \hookrightarrow & \mathcal{X} & \hookleftarrow & X \left(\substack{\text{smooth} \\ \text{proper}}\right) \\
 & & \downarrow & & \text{flat proper} \downarrow \pi & & \downarrow \\
 & & \operatorname{Spec} k \left(\substack{k\text{: residue} \\ \text{field}}\right) & \longrightarrow & \operatorname{Spec} R \left(\substack{R\text{: regular} \\ \text{local}}\right) & \longleftarrow & \operatorname{Spec} K \left(\substack{K\text{: quotient} \\ \text{field}}\right)
\end{array},$$

[CL17] proved the following very useful implication:

$$\begin{pmatrix} X \text{ admits a level } i \text{ integral (resp. rational)} \\ \text{decomposition of diagonal} \end{pmatrix} \implies \begin{pmatrix} Z \text{ also admits a level } i \text{ integral (resp. rational)} \\ \text{decomposition of diagonal} \end{pmatrix}$$

A very interesting consequence of these considerations is the following hierarchical upgrades of the (complex) hypersurface irrationality theorems of Totaro [T16] and Schreieder [S19]:

Totaro v.s. Schreieder, hierarchically upgraded

**Theorem 1.20.** *For very general smooth hypersurfaces*

$$X_d \subset \mathbb{P}^{n+1}$$

*defined over* $\mathbb{C}$, *the followings hold:*

**(Totaro-type):** *If*

$$d \geq 2\lceil \frac{n+2}{3} \rceil,$$

*then* $X_d$ *is not* $\big($*retract* $2$-*ruled* $=$ *retract* $-(n-2)$-*rationall*$\big)$,

**(Schreieder-type):** *For* $n \in \mathbb{N}$, *define* $\ell_2 n \ \big( \leq \lceil \log_2 n \rceil \big) \ \in \mathbb{N}$ *as follows:*

$$\begin{cases} n = l + r \\ 2^{l-1} - 2 \leq r \leq 2^l - 2 \end{cases}$$

$$\implies \ell_2 n := l = \min \Big\{ l \in \mathbb{N} \mid l + 2^l - 2 \geq n \Big\} \ \ \big( \leq \lceil \log_2 n \rceil \big)$$

*Then, if*

$$d \geq 2 + \ell_2 n \quad \textit{(which is satisfied if } d \geq 2 + \log_2 n\textit{)},$$

*then* $X_d$ *is not* $\big($*retract* $(n+1-L_2 n)$-*ruled* $=$ *retract* $-(L_2 n - 1)$-*rational*$\big)$,

If we do not care the hierarchy and concentrate in the irrationality problem, then Schreieder's theorem is clearly a magnificent improvement over Totaro's as Schreieder eloquently claimed in the introduction of [S19]. However, if we take into account the hierarchy, we immediately find Totaro and Schreieder are playing with different ball games: For an extremely difficult course of "rationality", whereas Professor Totaro's final exam was fairy easy and so just a small number of hypersurfaces failed, Professor Schreieder's final exam was much tougher and so much more hypersurfaces failed. Who is a better professor then? Of course, there is no clear answer for such a question. We can only say: both are good professors! To sum up, when the hypersurface irrationality theorems of Totaro and Schreieder are hierarchized, neither one of them supersedes the other: both are great theorems!

In another sequel, we shall present a sufficient criterion for a smooth projective $n$-dimensional variety $X$ over $k = \mathbb{C}$ to admit a dominant *morphism*

$$\mathcal{T}^k \times Z^{n-k} \to X, \tag{10}$$

where $\mathcal{T}^k$ is a $k$-dimensional generalized Bott tower in the sense of [GK94] [CMS10] which is a special kind of a toric manifold. Although this condition (10) is a priori not birational for it is descried in terms of *morphism*, it clearly implies the (separably) $k$-uniruledness = (separably) $\big(-(n-k)\big)$-unirationality of $X$, which is birationally invariant. Just like Suzuki's paper [S21], which deals with a problem of a similar problem, our proof is a repeated application of [AC12] which builds upon earlier work of [M79][MM86][dJS07]. However, [S21] could unfortunately offer just somewhat miserable general result. For instance, for $\mathbb{P}^{98}$, the general result of [S21] could only claim there passes through (just) $\mathbb{P}^{11}$ (not what everyone expects $\mathbb{P}^{98}$) at general points of $\mathbb{P}^{98}$. (Thus, [S21] could show $\mathbb{P}^{98}$ is just 11-uniruled!) Still, we can manage to offer a general sufficient condition, which when applied to $\mathbb{P}^n$, guarantees (10) with $k = n$, as everyone would require. Our sufficient criterion is too complicated to be stated here, but it gives us a very faithful estimate of times we can iterate the VMRT (=variety of minimal rational tangent) construction of [HM04]. For more precise account of our sufficient criterion, please consult our survey paper on this topic [M19].

The author would like to express his gratitude to Mark Spivakovsky for this encouragement and reference information concerning the local uniformization.

## 2. Glosarries

To state unramified sheaves and our main theorem in the next section §3, we here collect some conventions, notations and terminologies, together with some related background materials:

### 2.1. Pro-objects, Noetherian approximation,and Essentially smooth schemes of finite presentation.

- For a scheme $S$, let $Sch_S$ be the category of schemes over $S$ [GD71, p.215, I Def.2.1.2; p.226, I Def.2.6.1] which contains the following full subcategories:

(11)

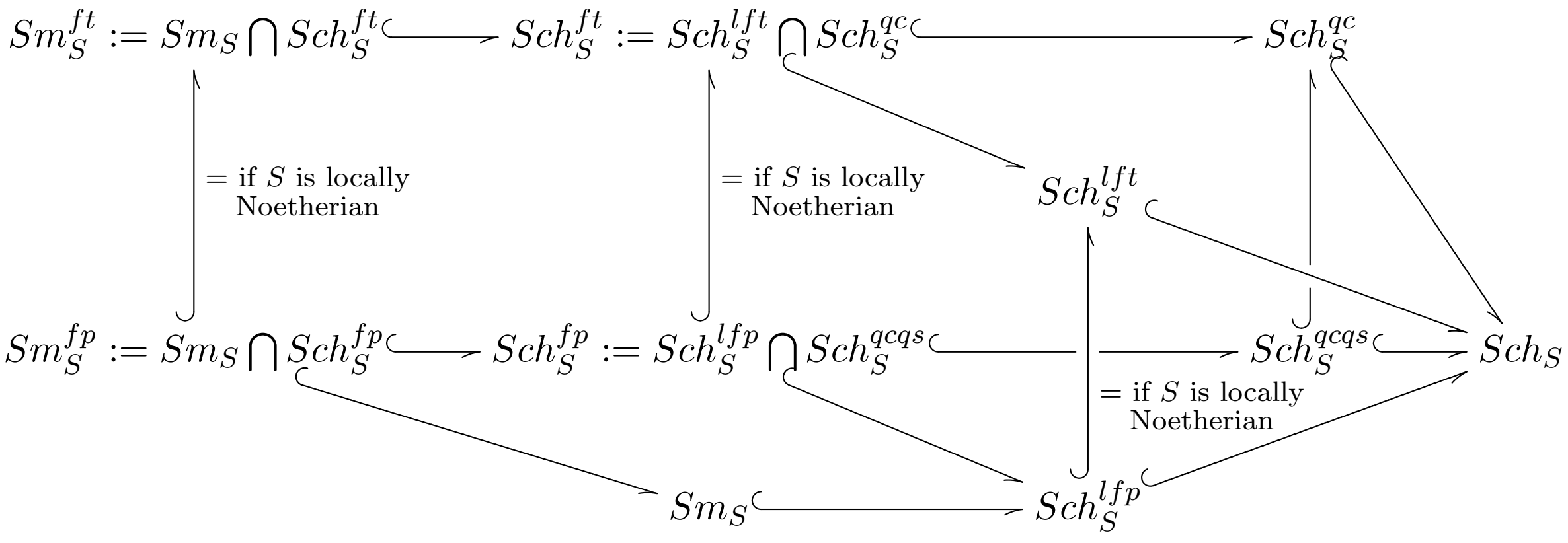

  - $Sch_S^{qcqs}$ is the full subcategory of $Sch_S$, consisting of *quasi-compact* [H77, p.80, II Exe 2.13; p.91, II Exe 3.2] and *quasi-separated* (i.e. the canonical morphism to Spec($\mathbb{Z}$) is quasi-compact) [G64, p.226,Def.1.2.1] [GD71, p.291,I Def.6.1.3] [sp18, 01KK] $S$-schemes.

- $Sch_S^{lfp}$ is the full subcategory of $Sch_S$, consisting of *locally finitely presented* [G64, p.230,Def.1.4.2] [GD71, p.297,I Def.6.2.1] [sp18, 01TP] $S$-schemes.
- $Sch_S^{fp}$ is the full subcategory of $Sch_S$, consisting of *finitely presented* (i.e. locally finitely presented and quasi-compact quasi-separated) [G64, p.234,Def.1.6.1] [GD71, p.305,I Def.6.3.7] [sp18, 01TP] $S$-schemes.
- $Sch_S^{lft}$ is the full subcategory of $Sch_S$, consisting of *locally of finite type* [G64, p.228,Def.1.3.2] [GD71, p.297,I Def.6.2.1] [sp18, 01T1] $S$-schemes.
- $Sch_S^{ft}$ is the full subcategory of $Sch_S$, consisting of $S$-schemes *of finite type* (i.e. locally of finite type and quasi-compact) [G64, p.233,Def.1.5.2] [GD71, p.304,I Def.6.3.2] [sp18, 01T1].
-
  * A locally finitely presented (resp. finitely presented) $S$-scheme is easily seen to be of locally finite type (resp. of finite type) [sp18, 01TW].
  $$Sch_S^{lft} \subset Sch_S^{lfp}, \qquad Sch_S^{ft} \subset Sch_S^{fp}.$$
  * Conversely, if the base scheme $S$ is locally Noetherian [sp18, 01OV], then any $S$-scheme of locally finite type (resp. of finite type) is locally finitely presented (resp. finitely presented) [sp18, 01TX]:
  $$Sch_S^{lft} = Sch_S^{lfp}, \qquad Sch_S^{ft} = Sch_S^{fp}.$$
  * In particular, when $S = \text{Spec } k$ for a field $k$, members of
  $$Sch_k^{lft} = Sch_k^{lfp} \tag{12}$$
  and
  $$Sch_k^{ft} = Sch_k^{fp} \tag{13}$$
  are also called [sp18, 06LG] locally algebraic $k$-scheme and algebraic $k$-scheme, respectively.
- $Sm_S$ is the full subcategory of $Sch_S$, consisting of *smooth* (i.e. locally finitely presented and formally smooth [G67, p.56,Def.17.1.1] [sp18, 02H0] ) [G67, p.61,Def.17.3.1] [sp18, 01V5,02H6] $S$-schemes.
- $Sm_S^{fp}$ is the full subcategory of $Sch_S$, consisting of finitely presented smooth $S$-schemes.
- $Sm_S^{ft}$ is the full subcategory of $Sch_S$, consisting of smooth $S$-schemes of finite type.

• For a category $\mathcal{C}$, let $Pro(\mathcal{C})$ be the category of *pro-objects* [sp18, Tag 05PX] of $\mathcal{C}$, whose object is a cofiltered diagram [sp18, Tag 04AZ] $X : \mathcal{I} \to \mathcal{C}$ (which we sometimes denote by $\left(X_\gamma, u_{\gamma,\gamma'} : X_{\gamma'} \to X_\gamma\right)_{(\gamma' \geq \gamma) \subset \mathcal{I}}$), and whose

morphism is given by

$$\begin{aligned}&\mathrm{Hom}_{Pro(\mathcal{C})}\left(\left(X_\gamma, u_{\gamma,\gamma'}: X_{\gamma'} \to X_\gamma\right)_{(\gamma'\geq\gamma)\subset\mathcal{I}}, \left(Y_\delta, v_{\delta,\delta'}: Y_{\delta'} \to Y_\delta\right)_{(\delta'\geq\delta)\subset\mathcal{J}}\right)\\ &\quad := \varprojlim_{\delta\in\mathcal{J}} \varinjlim_{\gamma\in\mathcal{I}} \mathrm{Hom}_{\mathcal{C}}\left(X_\gamma, Y_\delta\right).\end{aligned}$$

For a class $\mu$ of morphisms closed with respect to compositions in $\mathcal{C}$, let us denote by

$$Pro^{\mu}(\mathcal{C}) \subset Pro(\mathcal{C})$$

the full subcategory consisting of those $\left(X_\gamma, u_{\gamma,\gamma'}: X_{\gamma'} \to X_\gamma\right)_{(\gamma'\geq\gamma)\subset\mathcal{I}}$ with some $\alpha \in \mathcal{I}$ such that the transition morphism $u_{\gamma,\gamma'}: X_{\gamma'} \to X_\gamma$ belongs to $\mu$ for any $(\gamma' \geq \gamma \geq \alpha) \subset \mathcal{I}$.

- For the class $aff$ of affine morphisms [GD71, p.354,Def.9.1.1] [sp18, 01S6,01SC] in $Sch_S$, the projective limit exists [G66, p.8, Prop.8.2.3]:

(14)
$$\begin{aligned}\varprojlim : Ob\left(Pro^{aff}(Sch_S)\right) &\to Ob\left(Sch_S\right)\\ \left(X_\gamma, u_{\gamma,\gamma'}: X_{\gamma'} \to X_\gamma\right)_{(\gamma'\geq\gamma)\subset\mathcal{I}} &\mapsto \varprojlim_{\gamma\in\mathcal{I}} X_\gamma\end{aligned}$$

- When $S$ is quasi-compact and quasi-separated, (14) leads to the following very pleasing categorical equivalence:

**Theorem 2.1.** (Categorified "Noetherian approximation") [a]
When $S$ is quasi-compact and quasi-separated, we have the category equivalence:

(15)
$$\begin{aligned}\varprojlim : Pro^{aff}(Sch_S^{fp}) &\overset{\cong}{\to} Sch_S^{qcqs}\\ \left(X_\gamma, u_{\gamma,\gamma'}: X_{\gamma'} \to X_\gamma\right)_{(\gamma'\geq\gamma)\subset\mathcal{I}} &\mapsto \varprojlim_{\gamma\in\mathcal{I}} X_\gamma\end{aligned}$$

In particular, we have the following isomorphisms of morphism sets:

$$\begin{aligned}&\mathrm{Hom}_{Sch_S^{qcqs}}\left(\varprojlim_{\gamma\in\mathcal{I}} X_\gamma,\ \varprojlim_{\delta\in\mathcal{J}} Y_\delta\right) = \varprojlim_{\delta\in\mathcal{J}} \mathrm{Hom}_{Sch_S^{qcqs}}\left(\varprojlim_{\gamma\in\mathcal{I}} X_\gamma,\ Y_\delta\right)\\ &\quad = \varprojlim_{\delta\in\mathcal{J}} \varinjlim_{\gamma\in\mathcal{I}} \mathrm{Hom}_{Sch_S^{qcqs}}\left(X_\gamma, Y_\delta\right) =:\\ &\mathrm{Hom}_{Pro^{aff}(Sch_S^{fp})}\left(\left(X_\gamma, u_{\gamma,\gamma'}: X_{\gamma'} \to X_\gamma\right)_{(\gamma'\geq\gamma)\subset\mathcal{I}}, \left(Y_\delta, v_{\delta,\delta'}: Y_{\delta'} \to Y_\delta\right)_{(\delta'\geq\delta)\subset\mathcal{J}}\right)\end{aligned}$$

[a] Usually, "Noetherian approximation" is stated simply at the level of objects, especially when $S$ is a Noetherian [G60, p.140,Def.6.1.1, p.141,Prop.6.4.2] [sp18, 01OV] (i.e. locally-Noetherian and quasi-compact) scheme (then $S$ automatically becomes quasi-separated too [G64, p.228,Cor.1.2.8][sp18, 01OY]).
This is because the primary interest of these authors (e.g. [TT90, p.421,Th.C.9] [T11, p.5,Th.1.1.2] [CLO12, Th.1.2.2] [R15, p.108,Th.D]). is to generalize various results of finitely presented $S$-schemes (with $S$ Noetherian), which are still Noetherian by [sp18, 01T6] as being of finitely presented [sp18, 01TP] implies [sp18, 01TW] being of finite type [sp18, 01T1], to general quasi-compact and quasi-separated $S$-schemes.
However, all the necessary ingredient to derive our categorified version are already available in the literature as we see below.

*Proof.* First, when $S$ is quasi-compact and quasi-separated, Grothendieck [G66, p.50, (8.13.4), P.51, Prop.(8.13.5) ] proved the correspondence $\varprojlim$ in (14) induces a categorical fully faithful embedding

$$\varprojlim : Pro^{aff}\left(Sch_S^{fp}\right) \hookrightarrow \mathcal{C} \overset{\bigstar}{\subseteq} Sch_S^{qcqs} \ (\subseteq Sch_S)\,, \tag{16}$$

where $\mathcal{C} \subseteq Sch_S$ is the full subcategory consisting of $S$-schemes $X$, whose $S$-scheme structure morphism admits a decomposition of the following form:

$$X \xrightarrow{\text{affine}} X_0 \xrightarrow{\text{finitely presented}} S, \tag{17}$$

and the inclusive relation $\mathcal{C} \overset{\bigstar}{\subseteq} Sch_S^{qcqs}$ follows from the following facts:

- An affine morphism is quasi-compact and separated [sp18, 01S7], and so quasi-separated by their definitions [sp18, 01KK].
- A finitely presented morphism is defined to be a locally finitely presented morphism which is also quasi-compact and quasi-separated [G64, p.234,Def.1.6.1] [GD71, p.305,I Def.6.3.7] [sp18, 01TP].
- A composition of quasi-compact morphisms is quasi-compact [sp18, 01K6].
- A composition of quasi-separated morphisms is quasi-separated [sp18, 01KU].

Thus, to complete the proof of the category equivalence (15), it suffices to show every object of $Sch_S^{qcqs}$ is expressed as $\varprojlim_{\gamma\in\mathcal{I}} X_\gamma$ for some object $\left(X_\gamma, u_{\gamma,\gamma'} : X_{\gamma'} \to X_\gamma\right)_{(\gamma'\geq\gamma)\subset\mathcal{I}}$ in $Pro^{aff}\left(Sch_S^{fp}\right)$. However, up to some standard facts in the scheme theory (i.e. [sp18, 01TR,01S7,01SG]), this is taken care of in Temkin's work [T11, p.5,Th.1.1.2, p.4,Lem.1.1.1] [sp18, 01ZA,09MV], which generalizes the earlier work of Thomason-Trobaugh [TT90, p.421,Th.C.9] for the particular case of $S$ being the affine scheme corresponding to a Noetherian ring. □

- Suppose $S$ is quasi-compact and quasi-separated. Then, for a class $\mu$ of scheme morphisms closed with respect to compositions and a full subcategory $\mathcal{C} \subseteq Sch_S^{fp}$, let us call those objects of $Sch_S^{qcqs}$, which are contained in the essential image of

$$Pro^{aff\cap\mu}(\mathcal{C}) \subseteq Pro^{aff}(Sch_S^{fp}) \xrightarrow[(15)]{\cong} Sch_S^{qcqs},$$

<u>*$\mu$-essentially $\mathcal{C}$*</u>. When $aff \cap \mu = aff$, let us call them simply <u>*essentially $\mathcal{C}$*</u>. For instance, let us consider inclusions of subcategories

$$\begin{aligned}&\mathcal{C} \hookrightarrow Pro^{aff\cap oi}(\mathcal{C}) \subset Pro^{aff\cap ét}(\mathcal{C}) \subset Pro^{aff}(\mathcal{C})\\ &\subseteq Pro^{aff}(Sch_S^{fp}) \xrightarrow[(15)]{\cong} Sch_S^{qcqs},\end{aligned} \tag{18}$$

where *oi* is the collection of open immersions [sp18, 01IO] and *ét* is the collection of étale (i.e. locally finitely presented and formally étale [G67, p.56,Def.17.1.1]) morphisms [G67, p.61,Def.17.3.1,p.62,Cor.17.3.5]. (<u>*Warning:*</u> *Our terminologies here are not standard ones, though there is no such standard ones!* ) [7]

[7] For instance, "essentially affine" in the sense of Grothendieck [G66, p.50,8.13.4] is our essentially finitely presented [GD71, p.322,Prop.6.9.16(iii)] [T11, p.4,Lem.1.1.1] (which turned out to exhaust entire $Sch_S^{qcqs}$ by the aforecited Temkin's theorem [T11, p.5,Th.1.1.2]), "essentially of finite type"

- From the above discussion and the categorified "Noetherian approximation" Theorem 2.1, we immediately obtain the following very important corollary:

**Corollary 2.2.** *When the base scheme $S$ is quasi-compact and quasi-separated, let us suppose a subcategory $\mathcal{C} \subseteq Sch_S^{fp}$ and its presheaf $F$ valued in a cocomplete category $\mathcal{D}$ is given:*

$$F : \mathcal{C}^{op} \rightarrow \mathcal{D} \tag{19}$$

*Then, $F$ admits a canonical extension to the full subcategory of essentially $\mathcal{C}$ schemes $Pro^{aff}(\mathcal{C}) \hookrightarrow Sch_S^{qcqs}$ by the left Kan extension :*

$$\begin{array}{rccc} F: & Pro^{aff}(\mathcal{C})^{op} & \rightarrow & \mathcal{D} \\ & \varprojlim_{\gamma \in \mathcal{I}} X_\gamma & \mapsto & \varinjlim_{\gamma \in \mathcal{I}} F(X_\gamma) \end{array} \tag{20}$$

*Note that* (20) *gives the canonical extension of* (19) *to the essential image of $Pro^{aff}(\mathcal{C}) \subseteq Sch_S^{qcqs}$ in* (18)*, thanks to Theorem 2.1, the categorified "Noetherian approximation."*

**From now on, we shall assume our base scheme $S$ to be the affine scheme $Spec\ k$ corresponding to a perfect base field $k$.**

Amongst of all, we shall later apply Corollary 2.2 with $S$ a perfect field $k$ and $\mathcal{C} = Sm_k^{ft}$ when we study Morel's unramified sheaf.

2.2. **Regular local ring, Essentially smooth algebra of finite presentation, Popescu's theorem.** In order to treat a given regular local ring as an essentially smooth affine scheme, we may turn our attention to *Swan's question* [P89, p.121, Question (1.2)], which asks us to express a regular local ring as a filtered inductive limit of regular local rings essentially of finite type. While Swan's question had been known to be a tough problem, Popescu managed to prove the so-called the *General Néron desingularization theorem* (a.k.a. *Popescu's theorem*):

**Theorem 2.3.** (Popescu [P86, Theorem 1.8] [P90]) *For a ring homomorphism $\sigma : A \rightarrow B$ of Noetherian rings, the following statements are equivalent:*

(lim): *$B$ is a filtered colimit of smooth $A$-algebras of finite type;*
(reg): *$\sigma$ is regular* [M89, p.256][MR10, p.83, Definition 5.3.1][sp18, 07BZ].

□

Actually, the implication: (lim) $\Longrightarrow$ (reg) can be poven easily by elementary commutative ring theory [sp18, 07EP]. However, the most concise and

---

in the sense of [HKK17, Def.2.5], where the base scheme is assumed to be (separated) Noetherian which allowed us to apply [sp18, 01TX], is our open-immersive-essentially of finite presentation, and more generally, "essentially $\mathcal{C}$" in the sense of [DFJK21, p.567,Def.B.2] is our étale-essentially $\mathcal{C}$ [DFJK21, p.567,Prop.B.5]. Now, let us suppose the base scheme $S$ is locally Noetherian, which allows us to apply [sp18, 01TX,01TW] to conclude being of finitely presented is equivalent to being of finite type, some literature imposed finite type assumption to their definition of smooth schemes. Then, "essentially smooth" in the sense of [M12, p.vi] is our étale-essentially smooth of finite presentation, but "essentially smooth" in the sense of [B21, p.4,1.4] [EKW21, p.5,Rem.2.2] is our essentially smooth of finite presentation.

instantaneous, but much deeper proof of this implication is given by Andrés characterization of regular homomorphisms via the André-Quillen homology [A74, p.331, Supplement, Theorem 30] (where readers might wish to see more commutative ring theory results like [sp18, 07PM,07C0], but such a result is explicitly given in [MR10, p.108, Proposition 5.6.2]):

$$\sigma : A \to B \text{ is regular} \iff H_1(A, B, -) = 0,$$

and the fact that the André-Quillen homology commutes with colimits

$$H_1(A, \varinjlim_{\lambda} B_\lambda, -) = \varinjlim_{\lambda} H_1(A, B_\lambda, -)$$

[A74, p.45, III Proposition 35] (c.f. [T95, p.274]). Anyway, Popescu's proof [P86, Theorem 1.8] [P90] and subsequent improvements of the proof and the statements by others [O94][S98][S99] (see also [MR10, p.83, Theorem 5.3.2] [sp18, 07GC]), concentrating upon the other implication: (reg) $\Longrightarrow$ (lim).

In fact, the implication: (reg) $\Longrightarrow$ (lim) is exactly what concerns us in our case: $A = k, B = R$, a Noetherian local ring with maximal ideal $\mathfrak{m}_R$ and residue field $k_R = R/\mathfrak{m}_R$. Here, we have the following very useful characterizations of the regularity of $\sigma : k \to R$, due to Grothendieck:

**Proposition 2.4.** *Let $k$ be a field with characteristic exponent $p$, and let $R$ be a Noetherian local $k$-algebra. Then the following are equivalent:*
- (r) *The structure homomorphism $\sigma : k \to R$ of the $k$ algebra $R$ is a regular homomorphism.*
- (g) *$R$ is geometrically regular over $k$.*
- (p) *for all finite extension $k'$ of $k$ such that $k'^p \subset k$, $R \otimes_k k'$ is a regular ring.*
- (f) *$R$ is formally smooth over $k$ in the $\mathfrak{m}_R$ topology (i.e. $R$ is $\mathfrak{m}_R$-smooth over $k$).*

*Proof.* Actually, the equivalence $(r) \iff (g)$ follows immediately from the definition of regularity of a ring homomorphism [M89, p.256] and of geometric regularity of a $k$-algebra [M89, p.255], and the rest of the equivalences $(g) \iff (p) \iff (f)$ are shown by Grothendieck [G64, p.114, 0-Definition 19.9.1, p.204, 0-Theorem 22.5.8] (see also [M89, p.213, p.219, Theorem 28.7] for Faltings' short proof).

□

Since the above condition (p) in Proposition 2.4 is automatically satisfied when the base field $k$ is perfect and $R$ is a regular local $k$-algebra, together with the Popescu's theorem Theorem 2.3, we immediately obtain the following very useful corollary:

**Corollary 2.5.** *Any regular local ring over a perfect base field $k$ is a filtered colimit of smooth $k$-algebras of finite type.*
*In particular, the affine scheme of an arbitrary regular local ring over a perfect base field $k$ belongs to $Pro^{Aff \cap oi}\left(Sm_k^{ft}\right) \subset Pro^{Aff}\left(Sm_k^{ft}\right)$.*

We note that the condition (p) in Proposition 2.4 is another reason, in addition to Proposition 1.8, why we require our base field $k$ to be perfect.

While this Corollary 2.5 of Popescu's deep Theorem 2.3 is theoretically very pleasing, we can get away with applying it for the purpose of this paper, thanks to our Theorem 3.6.

## 3. Local uniformization, Abhyankaer and Geometric valuation

- For any $x \in X$, let $\mathfrak{m}_x$ be the maximal ideal of the local ring $\mathcal{O}_{X,x}$, and denote its residue field by $\kappa(x) := \mathcal{O}_{X,x}/\mathfrak{m}_x$.
- [sp18, 09FZ] We say $K/k$ is a *finitely generated field extension* if there exists a finite subset $S \subset K$ such that the smallest subfield $k(S)$ of $K$ containing $k$ and $S$ is $K$.
- Given a valuation $v : K \to \Gamma \cup \{+\infty\}$ on a field $K$ to a totally ordered commutative group $\Gamma$ (and the extra greatest element $+\infty$), we obtain the following familiar algebraic objects (see e.g. [ZS60b, p.34, VI, §8] [B72, VI, §3, Definition 1, Proposition 2] [M89, §10] ):
  - The *valuation ring* $V_v := \{x \in K \mid v(x) \geq 0\}$, which is a local ring with the maximal ideal $\mathfrak{m}_v := \{x \in K \mid v(x) > 0\}$ and the *residue field* $\kappa(v) := V_v/\mathfrak{m}_v$;
  - The *value group* $\Gamma_v := v(K^*) \cong K^*/(V_v)^*$, which is a totally ordered commutative group.

  From these, we have the following numerical invariants of a valuation $\nu$ :
  - ([ZS60b, p.50, VI, §10] [V06, p.485, Definition]) *The rational rank of $v$* is defined by $rat.rank(v) := \dim_{\mathbb{Q}} (\Gamma_v \otimes_{\mathbb{Z}} \mathbb{Q})$ ;
  - ([ZS60b, p.39, VI, §10; p.9, VI, §3, Definition 1; p.40, VI, §10, Theorem 15] [V06, p.482, Definition, Corollary; p.483, Proposition 1.8; p.484, lines 3-5]) *The rank of $v$*, denoted by $rank(v)$, is defined by the maximal length $r \in \mathbb{Z}_{\geq 0}$ of either one of the following three chains:
    - $\Gamma_v = \Delta_0 \supsetneq \Delta_1 \supsetneq \cdots \supsetneq \Delta_{r-1} \supsetneq \Delta_r = 0$, a chain of *isolated subgroups* [8] of the totally ordered commutative group $\Gamma_v$.
    - $V_v = V_r \subsetneq V_{r-1} \subsetneq \cdots \subsetneq V_1 \subsetneq V_0 = K$, a chain of valuation rings of $K$, containing the valuation ring $V_v$.
    - $\mathfrak{m}_v = \mathfrak{p}_r \supsetneq \mathfrak{p}_{r-1} \supsetneq \cdots \supsetneq \mathfrak{p}_1 \supsetneq \mathfrak{p}_0 = 0$, a chain of prime ideals of $V_v$.

[8]([ZS60b, p.40,VI, §10] [V06, p.481-482, Definition]) A subgroup $H$ of a totally ordered commutative group $A$ is called *isolated* if any $a \in A$ with $0 \leq a \leq h \in H$ is contained in $H$, i.e. $a \in H$.

Actually, these chains are interrelated as follows:

(21)

$$\{\Delta_i\}_{0\le i\le r}$$

$$\mathfrak{p}_i \mapsto \operatorname{Ker}\Big(\Gamma_v \xleftarrow[\cong]{v^{-1}} K^*/V_v^* \twoheadrightarrow K^*/((V_v)_{\mathfrak{p}_i})^*\Big) =: \Delta_i$$

$$\Delta_i \mapsto \Big(K^* \twoheadrightarrow K^*/V_v^* \xrightarrow[\cong]{v} \Gamma_v \twoheadrightarrow \Gamma_v/\Delta_i\Big)^{-1}(\ge 0) \cup \{0\} =: V_i$$

$$\Delta_i := \operatorname{Ker}\Big(\Gamma_v \xleftarrow{v^{-1}} [\cong] K^*/V_v^* \twoheadrightarrow K^*/V_i^*\Big) \hookleftarrow V_i$$

$$\mathfrak{p}_i := \Big(K^* \twoheadrightarrow K^*/V_v^* \xrightarrow[\cong]{v} \Gamma_v \twoheadrightarrow \Gamma_v/\Delta_i\Big)^{-1}(>0) \cup \{0\} \hookleftarrow \Delta_i$$

$$\{\mathfrak{p}_i\}_{r\ge i\ge 0} \rightleftarrows \{V_i\}_{r\ge i\ge 0}$$

$$\mathfrak{p}_i \mapsto (V_r)_{\mathfrak{p}_i} =: V_i$$

$$\mathfrak{p}_i := (\text{maximal ideal of } V_i) \hookleftarrow V_i$$

We also note that $rank(v)$ is the Krull dimension of $V_v$ :

$$rank(v) = \dim V_v.$$

– ( [ZS60b, p.41, §10, Theorem 16, p.42, §10, just before Corollary 1] [M89, p.78, Theorem 11.1]) The valuation ring $V_v$ is Noetherian, if and only if $\Gamma_v \cong \mathbb{Z}$ as ordered abelian groups, i.e. $v$ is a *discrete valuation of rank* 1 [9].

– Between the above numerical invariants, the following inequality holds ([B72, p.438, VI, §10.2, Corollary, Proposition 4]):

(22) $$rank(v) \le rat.rank(v).$$

When $\Gamma_v$ is finitely generated, if the equality $rank(v) = rat.rank(v)$ holds in (22), then $\Gamma_v$ is *discrete* in the sense of [ZS60b, p.48, -3rd line] [V06, p.484, Definition], i.e. $\Gamma_v \cong \mathbb{Z}^{rat.rank(v)}$ as totally ordered commutative groups, where $\mathbb{Z}^{rat.rank(v)}$ is given the lexicographical order.

– Suppose $K$ is a field extension of its subfield $k$, then a valuation $v : K \to \Gamma \cup \{+\infty\}$ is called a *valuation of $K/k$*, if $v(k^*) = 0$ [V06, p.481, 1.1].

– Given a valuation $v : K \to \Gamma \cup \{+\infty\}$ of $K/k$, we have a couple of field extensions:

$$K/k, \qquad \kappa(v)/k,$$

for which let us consider respective transcendental degrees:

(23) $$tr.deg_k K, \qquad tr.deg_k \kappa(v).$$

The latter is called *the dimension of $v$* ([ZS60b, p.34, VI, §8] [V06, p.494]), and denoted by

(24) $$\dim_k(v) := tr.deg_k \kappa(v).$$

– ([A56a][ZS60b, p.331, Appendix 2, Proposition 2; p.335, Appendix 2, Proposition 3] [B72, p.439, VI, §10.3, Corollary 1]) For any valuation $v$ of $K/k$, the following *Abhyankar's inequality* holds:

(25) $$rat.rank(v) + \dim_k(v) \ \le\ tr.deg_k K$$

[9] This is what is usually called DVR in most literature, especially in the area of algebraic number theory.

When $K/k$ is a finitely generated field extension, if the equality $rat.rank(v) + \dim_k(v) = tr.deg_k K$ holds in Abhyankar's inequality (25), then $v$ is called an <u>*Abhyankar valution*</u> on $K/k$.

- For a valuation $v$ on $K/k$, we have the following inequality from (22) and Abhyankar's inequality (25):

$$rank(v) + \dim_k(v) \ \leq\ tr.deg_k K \tag{26}$$

When $K/k$ is a finitely generated field extension, if the equality $rank(v) + \dim_k(v) = tr.deg_k K$ holds in (26), then $v$ is called a <u>*geometric valuation*</u> on $K/k$ in some motivic literature (e.g. [R96, p.328, 3rd and 4th paragraphs] for the rank 1 case and [M08, p.53, 2nd paragraph] for arbitrary rank cases). If $v$ is a geometric valuation of $K/k$, then $v$ is of course Abhyankar. Furthermore, $\Gamma_v$ is not only finitely generated, but also discrete: $\Gamma_v \cong \mathbb{Z}^{rat.rank(v)}$ as totally ordered commutative groups, where $\mathbb{Z}^{rat.rank(v)}$ is given the lexicographical order (for a direct proof of this fact, see [ZS60b, p.90, VI §14, Corollary]).

- A rank 1 geometric valuation $v$ is nothing but a rational rank 1 Abhyankar valuation, and was called a <u>*prime divisor*</u> [ZS60b, p.88 VI, §14; p.89, VI, §14, Theorem 31], but is also called a <u>*divisorial valuation*</u> in many recent literature (see e.g. the proceedings [v14]).
- [ZS60b, p.43–44, VI, §10] (see also [M89, p.72, Theorem 10.1] [V00, §4] [V06, p.485–487])
  - For a valuation $v_1$ on $K/k$ and a valuation $\overline{v}$ on $\kappa(v_1)/k$, their <u>*composite valuation*</u> $v_1 \circ \overline{v}$ on $K/k$ is defined so that its valuation ring $V_{v_1\circ\overline{v}}$ is given by $\pi_{v_1}^{-1}(V_{\overline{v}})$ [M89, p.72, 4, Theorem 10.1]:

$$\begin{array}{ccccc} K & \hookleftarrow & V_{v_1} & \hookleftarrow & V_{v_1\circ\overline{v}} = \pi_{v_1}^{-1}(V_{\overline{v}}) \\ & & \downarrow{\scriptstyle \pi_{v_1}} & & \downarrow{\scriptstyle \pi_{v_1}} \\ & & \kappa(v_1) = V_{v_1}/\mathfrak{m}_{v_1} & \hookleftarrow & V_{\overline{v}} \end{array} \tag{27}$$

and its value group $\Gamma_{v_1\circ\overline{v}}$ to enjoy the following short exact sequence [ZS60b, p.43, VI, §10 Theorem 17]:

$$0 \to \Gamma_{\overline{v}} \to \Gamma_{v_1\circ\overline{v}} \to \Gamma_{v_1} \to 0$$

If this sequence splits (simply as abelian groups), then

$$\Gamma_{v_1\circ\overline{v}} \ \cong\ \Gamma_{v_1} \times \Gamma_{\overline{v}}$$

as totally ordered commutative groups, where $\Gamma_{v_1} \times \Gamma_{\overline{v}}$ is endowed with the lexicographic order [V00, p.553, Proposition 4.3, p.554, Remarque.]. We also note from (27) that

$$\begin{cases} V_{\overline{v}} & \overset{\cong}{\longleftarrow}\ V_{v_1\circ\overline{v}}/\mathfrak{m}_{v_1} \cap V_{v_1\circ\overline{v}} = V_{v_1\circ\overline{v}}/\left(\underline{\mathfrak{m}_{v_1}}\right) \\ \mathfrak{m}_{\overline{v}} & \overset{\cong}{\longleftarrow}\ \mathfrak{m}_{v_1\circ\overline{v}}/\mathfrak{m}_{v_1} \cap \mathfrak{m}_{v_1\circ\overline{v}} = \mathfrak{m}_{v_1\circ\overline{v}}/\left(\underline{\mathfrak{m}_{v_1}}\right) = \overline{\mathfrak{m}_{v_1\circ\overline{v}}}, \end{cases} \tag{28}$$

where, here and after whenever there is no danger of confusion, we shall adapt the following simplified notations of ideals:

  - ∗ If a commutative ring $A$ contains an ideal $\mathfrak{a}$ and a subring $B$, denote by $\underline{\mathfrak{a}} := \mathfrak{a} \cap B$, the restriction of the ideal $\mathfrak{a} \subset A$ to $B$.

* If a commutative ring $B$ admits a quotient surjective ring homomorphism $\pi : B \twoheadrightarrow C$ and an ideal $\mathfrak{b}$ containing Ker$\pi$, denote by $\overline{\mathfrak{b}} := \pi(\mathfrak{b})$, the image of the ideal $\mathfrak{b} \subset B$ in $C$.

(29)
$$\begin{array}{ccccccc} & & & & \underline{\mathfrak{a}} := B \cap \mathfrak{a} & \hookrightarrow & \mathfrak{a} \\ & & & & \downarrow & & \downarrow \\ \text{Ker } \pi & \hookrightarrow & \mathfrak{b} & \hookrightarrow & B & \hookrightarrow & A \\ & & \downarrow & & \downarrow \pi & & \\ & & \overline{\mathfrak{b}} := \pi(\mathfrak{b}) & \hookrightarrow & C & & \end{array}$$

– Conversely, let $v$ be a valuation on $K$ whose rank is strictly larger than 1, or whose value group $\Gamma_v$ contains a proper nontrivial isolated subgroup $\Delta$. Then:
  * Define a valuation $v_1$ on $K$ by the composition:
$$v_1 : K^* \xrightarrow{v} \Gamma_v \twoheadrightarrow \Gamma_v/\Delta =: \Gamma_{v_1}.$$
  * Define a group homomorphism $\overline{v} : \kappa(v_1)^* = (V_{v_1}/\mathfrak{m}_{v_1})^* \to K^*/(V_v)^*$ so that the following diagram commutes:
$$\begin{array}{ccc} (V_{v_1})^* & \xrightarrow{\text{canonical inclusion}} & K^* \\ \pi_{v_1} \downarrow & & \downarrow v \\ \kappa(v_1)^* = (V_{v_1}/\mathfrak{m}_{v_1})^* & \xrightarrow{\overline{v}} & K^*/(V_v)^* = \Gamma_v \end{array}$$
  To justify this definition of $\overline{v}$, it suffices to show:
  *For any* $1 + x \in 1 + \mathfrak{m}_{v_1} = \text{Ker } \pi_{v_1}$, $v(1 + x) = 0$.
  In fact, since $\mathfrak{m}_{v_1} \subset \mathfrak{m}_v$, we see $1+x \in 1+\mathfrak{m}_v \subset V_v$, i.e. $v(1+x) \geq 0$, in addition to $v(-x) > 0$. Thus, if we suppose $v(1 + x) > 0$ contrary to the above claim, then
$$0 = v(1) = v\left((1 + x) + (-x)\right) \geq \min\left(v(1 + x), v(-x)\right) > 0,$$
  an obvious contradiction.
  * $\overline{v}$ induces a valuation on $\kappa(v_1)$, whose valuation ring shall be denoted by $V_{\overline{v}}$ $(\subset \kappa(v_1))$, and so Ker $\overline{v} = (V_{\overline{v}})^*$. Consequently, we have the following short exact sequence:

(30)
$$\begin{array}{ccccccccc} 0 & \longrightarrow & \Gamma_{\overline{v}} \cong \Delta & \longrightarrow & \Gamma_v & \longrightarrow & \Gamma_v/\Delta =: \Gamma_{v_1} & \longrightarrow & 0 \\ & & \cong \uparrow & & \cong \uparrow & & \cong \uparrow & & \\ 0 & \longrightarrow & \kappa(v_1)^*/(V_{\overline{v}})^* \cong (V_{v_1})^*/(V_v)^* & \longrightarrow & K^*/(V_v)^* & \longrightarrow & K^*/(V_{v_1})^* & \longrightarrow & 0 \end{array}$$

  * We then see $v = v_1 \circ \overline{v}$.

– ([ZS60b, p.43] [V00, p.552, Corollaire] [V06, p.486, Corollary]) For the composite valuation $v = v_1 \circ \overline{v}$, we have the equalities:

(31)
$$\begin{aligned} rank(v) &= rank(v_1) + rank(\overline{v}) \\ rat.rank(v) &= rat.rank(v_1) + rat.rank(\overline{v}) \end{aligned}$$

– In general, any rank $r$ valuation $v$ on $K$, where we shall also denote $K$ by $\kappa(v_0)$, can be expressed (inductively with respect to $r$) uniquely as

(32) $$v = v_1 \circ \cdots \circ v_r,$$

where $v_i$ is a valuation on $\kappa(v_{i-1})$ for each $1 \leq i \leq r$.

– In terms of the valuations $v_i$ on $\kappa(v_{i-1})$ $(1 \leq i \leq r)$ and their composition, the "trinity correspondences"(21) can be restated via (30) (28) as follows:

(33) $$\begin{cases} \Gamma_{v_1\circ\cdots\circ v_r} = \Delta_0 \supsetneq \Delta_1 \supsetneq \cdots \supsetneq \Delta_{r-1} \supsetneq \Delta_r = 0 \\ \Delta_i = \Gamma_{v_{i+1}\circ\cdots\circ v_r} = \mathrm{Ker}\left(\Gamma_{v_1\circ\cdots\circ v_r} \twoheadrightarrow \Gamma_{v_1\circ\cdots\circ v_i}\right) \end{cases}$$
$$\begin{cases} V_{v_1\circ\cdots\circ v_r} = V_r \subsetneq V_{r-1} \subsetneq \cdots \subsetneq V_1 \subsetneq V_0 = K \\ V_i = V_{v_1\circ\cdots\circ v_i} = (V_{v_1\circ\cdots\circ v_r})_{\mathfrak{m}_{v_1\circ\cdots\circ v_i}\cap V_{v_1\circ\cdots\circ v_r}} \end{cases}$$
$$\begin{cases} \mathfrak{m}_{v_1\circ\cdots\circ v_r} = \mathfrak{p}_r \supsetneq \mathfrak{p}_{r-1} \supsetneq \cdots \supsetneq \mathfrak{p}_1 \supsetneq \mathfrak{p}_0 = 0 \\ \mathfrak{p}_i = (v_1 \circ \cdots \circ v_i)^{-1}\left((\Gamma_v/\Delta_i)_{>0}\right) = \mathfrak{m}_{v_1\circ\cdots\circ v_i} \cap V_{v_1\circ\cdots\circ v_r} \end{cases}$$

We now review some basic results of local uniformization following [NS14]:

- (cf. [NS14, p.410, Definition 2.8]) When $R$ is a subring of $K$ and $\mathfrak{m}_v \cap R = \mathfrak{p} \subset R$, we say $v$ has a *center on $R$ at $\mathfrak{p}$*.

  Let $v$ be a valuation on a field $K$ with $\mathfrak{m}_v \subset V_v \subset K = Frac(V_v)$ its maximal ideal $\mathfrak{m}_v$ and valuation ring $V_v$, as usual. Then, let us say the valuation ring $V_v$ is centerd in an integral domain $S$, if

(34) $$S \subseteq V_v \subset K = Frac(V_v) = Frac(S),$$

where we shall also call $\mathfrak{m}_v \cap S \subset S$ the center of the valuation $v$.

Given such an integral domain $S$ where a valuation ring $V_v$ is centered in, we define a local blowing up of $S$ with respect to $v$ to be a local ring $\left(S^{(1)}, \mathfrak{m}^{(1)}\right)$ of the form

(35) $$\begin{cases} S \subseteq S^{(1)} \subseteq V_v \subset K = Frac(V_v) = Frac(S^{(1)}) = Frac(S), \\ (S, \mathfrak{m}_v \cap S) \ \to \ \left(S^{(1)}, \mathfrak{m}^{(1)} = \mathfrak{m}_v \cap S^{(1)}\right) \ \to \ (V_v, \mathfrak{m}_v), \end{cases}$$

which is constructed in the following two steps:

**Step 1:** Pick elements $a_i, b_i \in S, i = 1, \ldots, s$ such that $v(b_i) \leq v(a_i)$ for all $i = 1, \ldots, s$ and define a pair of a subdomain and its prime ideal $(S', \mathfrak{m}')$ of $(V_v, \mathfrak{m}_v)$ by

(36) $$S' = S\left[\frac{a_1}{b_1}, \cdots, \frac{a_s}{b_s}\right] \quad \text{and} \quad \mathfrak{m}' = \mathfrak{m}_v \cap S'.$$

**Step 2:** The the desired local domain $\left(S^{(1)}, \mathfrak{m}^{(1)}\right)$ is the localization of $S'$ with respect to the prime ideal $\mathfrak{m}'$, that is,

(37) $$\begin{cases} S^{(1)} = S'_{\mathfrak{m}'} = \left\{ \frac{x}{y} \in F \;\middle|\; x \in S', y \in S' \setminus \mathfrak{m}' \right\} \\ \text{and} \quad \mathfrak{m}^{(1)} = \mathfrak{m}' S'_{\mathfrak{m}'} = \mathfrak{m}_v \cap S^{(1)}. \end{cases}$$

In other words, by the following composition:

$$
\begin{aligned}
&(S, \mathfrak{m}_v \cap S) \ \to \ \left( S\left[\frac{a_1}{b_1}, \cdots, \frac{a_s}{b_s}\right], \mathfrak{m}_v \cap S\left[\frac{a_1}{b_1}, \cdots, \frac{a_s}{b_s}\right] \right) \\
&\to \ \left( S\left[\frac{a_1}{b_1}, \cdots, \frac{a_s}{b_s}\right]_{\mathfrak{m}_v \cap S\left[\frac{a_1}{b_1}, \cdots, \frac{a_s}{b_s}\right]}, \mathfrak{m}_v \cap \left( S\left[\frac{a_1}{b_1}, \cdots, \frac{a_s}{b_s}\right]_{\mathfrak{m}_v \cap S\left[\frac{a_1}{b_1}, \cdots, \frac{a_s}{b_s}\right]} \right) \right) \qquad (38) \\
&=: \left( S^{(1)}, \mathfrak{m}^{(1)} = \mathfrak{m}_v \cap S^{(1)} \right) \ \to \ (V_v, \mathfrak{m}_v)
\end{aligned}
$$

While a general local blowing up might look complicated, it is always a composition of <u>simple</u> local blowing ups, i.e. those with $s = 1$ and so of the form

$$
(S, \mathfrak{m}_v \cap S) \ \to \ \left( S\left[\frac{a}{b}\right]_{\mathfrak{m}_v \cap S\left[\frac{a}{b}\right]}, \mathfrak{m}_v \cap \left( S\left[\frac{a}{b}\right]_{\mathfrak{m}_v \cap S\left[\frac{a}{b}\right]} \right) \right), \qquad (39)
$$

as the following proposition shows:

**Proposition 3.1.** (cf. [NS14, p.410, Lemma 2.9.]) *Every local blowing up can be decomposed as a finite sequence of simple local blowing ups.* □

On the other hand, given an input (34), a local blowing up with respect to a valuation $v$ returns an output (35). Of course, we may iterate this procedure, but, it turns out that any composition of such iterations can be realized in a single step:

**Proposition 3.2.** (i) (cf. [NS14, p.411, Lemma 2.10.]) *Suppose a valuation $v$ centers on a domain $S$ as in* (34). *Then the composition of every sequence of local blowing ups with respect to $v$*

$$
(S, \mathfrak{m}_v \cap S) \to \left( S^{(1)}, \mathfrak{m}^{(1)} = \mathfrak{m}_v \cap S^{(1)} \right) \to \ldots \to \left( S^{(t)}, \mathfrak{m}^{(t)} = \mathfrak{m}_v \cap S^{(t)} \right)
$$

*can be expressed as a single local blowing up*

$$
(S, \mathfrak{m}_v \cap S) \to \left( \widetilde{S}^{(1)}, \widetilde{\mathfrak{m}}^{(1)} = \mathfrak{m}_v \cap \widetilde{S}^{(1)} \right) \qquad (40)
$$

*with respect to $v$ , with $\widetilde{S}^{(1)} = S^{(t)}$.*
(ii) *In (i), suppose $v$ is a valuation on a finitely generated field extension $K/k$ and $S$ is a localization of a finitely generated $k$-algebra, then the outcome $S^{(t)}$ of a composition of blowing ups starting with $S$ is a localization of a finitely generated $k$-algebra.*

*Proof of (i).* Actually, the proof [NS14, p.411, Proof of Lemma 2.10.] of a slightly weaker [10] statement works without modification.

□

*Proof of (ii).* By (i), the outcome $S^{(t)}$ is realized by a single local blowing up (40), which by (38) can be written as

$$
S\left[\frac{a_1}{b_1}, \cdots, \frac{a_s}{b_s}\right]_{\mathfrak{m}_v \cap S\left[\frac{a_1}{b_1}, \cdots, \frac{a_s}{b_s}\right]},
$$

[10] In fact, the only difference between our Proposition 3.2 and [NS14, p.411, Lemma 2.10.] is that $S$ is allowed to be non Noetherian local in Proposition 3.2 unlike [NS14, p.411, Lemma 2.10.]. The same comment is applied to the above Proposition 3.1 and [NS14, p.410, Lemma 2.9.].

which is clearly a localization of a finitely generated $S$-algebra. Since $S$ is a localization of a finitely generated $k$-algebra, we see by induction that $S^{(t)}$ is also a localization of a finitely generated $k$-algebra, as desired.

□

- (cf. [NS14, p.412, Lemma 2.13, p.413, Definition 2.14]) When a valuation $v$ centers on a Noetherian domain $S$ as in (34), there exists a canonical local blowing up of $S$ with respect to $v$ for each ideal $I \subset S$ :

$$(S, \mathfrak{m}_v \cap S) \to \left(S_I^{(1)}, \mathfrak{m}_I^{(1)} = \mathfrak{m}_v \cap S_I^{(1)}\right) \to (V_v, \mathfrak{m}_v),$$

said to be the local blowing up of $S$ with respect to $v$ along $I$ which is obtained by carefully picking elements of $S$ in Step 1 of the construction (36) as follows:

**Step 1:** Pick $u_i \in I$ $(0 \le i \le s)$ so that
– $v(u_0) = \min\{v(u) \mid u \in I\}$;
– $u_i \in I$ $(0 \le i \le s)$ generate $I$ as an $S$-module. [11]

then define a pair of subdomain and its prime ideal $(S'_I, \mathfrak{m}'_I)$ of $(V_v, \mathfrak{m}_v)$ by

(41) $$S'_I = S\left[\frac{u_1}{u_0}, \cdots, \frac{u_s}{u_0}\right] \quad \text{and} \quad \mathfrak{m}'_I = \mathfrak{m}_v \cap S'_I.$$

**Step 2:** The the desired local domain $\left(S_I^{(1)}, \mathfrak{m}_I^{(1)}\right)$ is given by:

(42) $$\begin{cases} S_I^{(1)} = (S'_I)_{\mathfrak{m}'_I} = \left\{ \frac{x}{y} \in F \;\middle|\; x \in S'_I, y \in S'_I \setminus \mathfrak{m}'_I \right\} \\ \text{and} \quad \mathfrak{m}_I^{(1)} = \mathfrak{m}_v \cap S_I^{(1)}. \end{cases}$$

**Proposition 3.3.** (cf. [NS14, p.412, Lemma 2.13, p.413, Definition 2.14]) *The above construction of local blowing up with respect to $v$ is independent of particular choices of $u_i \in S$ $(0 \le i \le s)$ in (41). This justifies the terminology: the local blowing up of $S$ with respect to $v$ along $I$.*

*Proof.* Once again, the proof [NS14, p.412, Proof of Lemma 2.13.] of a slightly weaker [12] statement works without modification. □

- When $(S, \mathfrak{m}_S)$ is a local ring whose quotient field is $K$ and $\mathfrak{m}_S$ the center of $v$, the *local uniformization problem* [Z40] [V06] [NS16] asks us to find a local blowing up of $(S, \mathfrak{m}_S)$ with respect to $v$ [NS14, p.410, Definition 2.8]

(43) $$(S, \mathfrak{m}_S) \to (R, \mathfrak{m}_R) \ (\to (V_v, \mathfrak{m}_v))$$

such that $(R, \mathfrak{m}_R)$ is a regular local ring. Local uniformization may be regarded as a local version of the resolution of singularities, and was proved affirmatively when the base field $k$ is characteristic 0 by Zariski [Z40].

To state some positive results of local uniformization, let us review Abhyankar valuations in more detail.

- ([A56a] [ZS60b, p.335, Appendix 2, Proposition 3], see also [C21, 1.Introduction] and our (21) (33).) For an Abhyankar valuation $v$ of $K/k$,

[11] We have used the Noetherian assumption of $S$ here.

[12] Once again, the only difference between our situation and [NS14, p.412, Lemma 2.13.] is that $S$ is allowed to be non local in our situation unlike [NS14, p.412, Lemma 2.13.].

– the unique composition (32)

$$v = v_1 \circ v_2 \circ \cdots \circ v_{rank(v)},$$

is given so that each $v_i$ $(1 \leq i \leq rank(v))$ is a rank 1 Abhyankar valuation of a finitely generated field extension $\kappa(v_{i-1})/k$;

– its longest chain of isolated subgroups of the value group $\Gamma_v$ (21) (33):

$$0 = \Delta_{rank(v)} \subsetneq \Delta_{rank(v)-1} \subsetneq \cdots \subsetneq \Delta_1 \subsetneq \Delta_0 = \Gamma_v,$$

is given uniquely so that $\Gamma_v/\Delta_i$ becomes the value group of $v_1 \circ \cdots \circ v_i$ $\quad (1 \leq i \leq rank(v))$, enjoying the following properties:

  * each $\Delta_{i-1}/\Delta_i$ $(1 \leq i \leq rank(v))$ is the value group of the rank 1 Abhyankar valuation $v_i$ of $\kappa(v_{i-1})/k$ and can be expressed as $\Delta_{i-1}/\Delta_i \cong \mathbb{Z}^{rr_i}$ for some $rr_i \in \mathbb{N}$ as abelian groups, and can be embedded $\Delta_{i-1}/\Delta_i \subset \mathbb{R}$ as totally ordered commutative groups [ZS60b, p.45, VI §10].
  * $\Gamma_v \cong \mathbb{Z}^{rat.rank(v)}$ $(rat.rank(v) = \sum_{1\leq i\leq r} rr_i)$ as abelian groups, and can be embedded $\Gamma_v \subset \mathbb{R}^{rank(v)}$ as totally ordered commutative groups, where $\mathbb{R}^{rank(v)}$ is given the lexicographical order.

– any intermediate composite valuation $v_i \circ \cdots \circ v_j$ of $\kappa(i-1)/k$ $\quad (1 \leq i \leq j \leq rank(v))$ is still Abhyankar, and its value group $\Gamma_{v_i\circ\cdots\circ v_j}$ and valuation ring $V_{v_i\circ\cdots\circ v_j}$ are given as follows (see also (30) (28), in addition to (21) (33)):

$$\Gamma_{v_i\circ\cdots\circ v_j} = \Gamma_{v_i\circ\cdots\circ v_{rank(v)}}/\Delta_j = \Delta_{i-1}/\Delta_j;$$

$$V_{v_i\circ\cdots\circ v_j} = \begin{cases} \dfrac{V_{v_1\circ\cdots\circ v_j}}{\mathfrak{m}_{v_1\circ\cdots\circ v_{i-1}} \cap V_{v_1\circ\cdots\circ v_j}} = \dfrac{\left(V_{v_1\circ\cdots\circ v_{rank(v)}}\right)_{\mathfrak{m}_{v_1\circ\cdots\circ v_j}}}{\mathfrak{m}_{v_1\circ\cdots\circ v_{i-1}} \cap \left(V_{v_1\circ\cdots\circ v_{rank(v)}}\right)_{\mathfrak{m}_{v_1\circ\cdots\circ v_j}}} \\ \left(V_{v_i\circ\cdots\circ v_{rank(v)}}\right)_{\mathfrak{m}_{v_i\circ\cdots\circ v_j} \cap V_{v_i\circ\cdots\circ v_{rank(v)}}} = \left(\dfrac{V_{v_1\circ\cdots\circ v_{rank(v)}}}{\mathfrak{m}_{v_1\circ\cdots\circ v_{i-1}} \cap V_{v_1\circ\cdots\circ v_{rank(v)}}}\right)_{\overline{\left(\underline{\mathfrak{m}_{v_1\circ\cdots\circ v_j}}\right)}} \end{cases}$$

$$= \begin{cases} \dfrac{(V_v)_{\mathfrak{m}_j}}{\left(\underline{\mathfrak{m}_{i-1}}\right)} \\ \left(\dfrac{V_v}{\left(\underline{\mathfrak{m}_{i-1}}\right)}\right)_{\overline{\left(\underline{\mathfrak{m}_j}\right)}} \end{cases}$$

Here, we have abbreviated so that $v := v_1 \circ \cdots \circ v_{rank(v)}$, $\mathfrak{m}_i = \mathfrak{m}_{v_1\circ\cdots\circ v_i}$ and adapted the convention of (29) as in (28).

They are interrelated so as to enjoy the following equivalences, where we set $r := rank(v)$ for ease of notation:

**Systems of value groups:**

(44)
$$\left\{ \Gamma_{v_i\circ\cdots\circ v_{j-1}} \twoheadleftarrow \Gamma_{v_i\circ\cdots\circ v_j} \hookleftarrow \Gamma_{v_{i+1}\circ\cdots\circ v_j} \right\}_{1\leq i\leq j\leq r}$$
$$= \left\{ \Delta_{i-1}/\Delta_{j-1} \twoheadleftarrow \Delta_{i-1}/\Delta_j \hookleftarrow \Delta_i/\Delta_j \right\}_{1\leq i\leq j\leq r}$$

Here at the entry $(i,j) = (i, i-1)$ $(1 \leq i \leq r+1$, we assign $\Gamma_{triv} = \Delta_{i-1}/\Delta_{i-1}$. More diagrammatically,

$$
\begin{array}{ccccccccccc}
\Gamma_{triv} & \twoheadleftarrow \cdots \twoheadleftarrow & \Gamma_{v_1\circ\cdots v_{i-1}} & \twoheadleftarrow & \Gamma_{v_1\circ\cdots v_i} & \twoheadleftarrow & \Gamma_{v_1\circ\cdots v_{i+1}} & \twoheadleftarrow \cdots \twoheadleftarrow & \Gamma_{v_1\circ\cdots v_{rank(v)}} = \Gamma_v \\
 & & \cup & & \cup & & \cup & & \cup \\
\ddots & \cdots & \vdots & \cdots & \vdots & \cdots & \vdots & \cdots & \vdots \\
 & & \uparrow & & \uparrow & & \uparrow & & \uparrow \\
 & & \Gamma_{triv} & \twoheadleftarrow & \Gamma_{v_i} & \twoheadleftarrow & \Gamma_{v_i\circ v_{i+1}} & \cdots \twoheadleftarrow & \Gamma_{v_i\circ\cdots v_{rank(v)}} \\
 & & & & \uparrow & & \uparrow & & \uparrow \\
 & & & & \Gamma_{triv} & \twoheadleftarrow & \Gamma_{v_{i+1}} & \cdots \twoheadleftarrow & \Gamma_{v_{i+1}\circ\cdots v_{rank(v)}} \\
 & & & & & & \uparrow & & \uparrow \\
 & & & & & & \Gamma_{triv} & \cdots \twoheadleftarrow & \Gamma_{v_{i+2}\circ\cdots v_{rank(v)}} \\
 & & & & & & & & \cup \\
 & & & & & & & \ddots & \vdots \\
 & & & & & & & & \uparrow \\
 & & & & & & & & \Gamma_{triv}
\end{array}
\tag{45}
$$

is identified with

$$
\begin{array}{ccccccccc}
\Delta_0/\Delta_0 & \twoheadleftarrow \cdots \twoheadleftarrow & \Delta_0/\Delta_{i-1} & \twoheadleftarrow & \Delta_0/\Delta_i & \twoheadleftarrow & \Delta_0/\Delta_{i+1} & \twoheadleftarrow \cdots \twoheadleftarrow & \Delta_0/\Delta_{rank(v)} \\
 & & \cup & & \cup & & \cup & & \cup \\
\ddots & \cdots & \vdots & \cdots & \vdots & \cdots & \vdots & \cdots & \vdots \\
 & & \uparrow & & \uparrow & & \uparrow & & \uparrow \\
 & & \Delta_{i-1}/\Delta_{i-1} & \twoheadleftarrow & \Delta_{i-1}/\Delta_i & \twoheadleftarrow & \Delta_{i-1}/\Delta_{i+1} & \cdots \twoheadleftarrow & \Delta_{i-1}/\Delta_{rank(v)} \\
 & & & & \uparrow & & \uparrow & & \uparrow \\
 & & & & \Delta_i/\Delta_i & \twoheadleftarrow & \Delta_i/\Delta_{i+1} & \cdots \twoheadleftarrow & \Delta_i/\Delta_{rank(v)} \\
 & & & & & & \uparrow & & \uparrow \\
 & & & & & & \Delta_{i+1}/\Delta_{i+1} & \cdots \twoheadleftarrow & \Delta_{i+1}/\Delta_{rank(v)} \\
 & & & & & & & & \cup \\
 & & & & & & & \ddots & \vdots \\
 & & & & & & & & \uparrow \\
 & & & & & & & & \Delta_{rank(v)}/\Delta_{rank(v)}
\end{array}
\tag{46}
$$

**Systems of valuation rings:**

$$
\begin{aligned}
&\left\{ V_{v_i\circ\cdots\circ v_{j-1}} \hookleftarrow V_{v_i\circ\cdots\circ v_j} \twoheadrightarrow V_{v_{i+1}\circ\cdots\circ v_j} \right\}_{1\leq i\leq j\leq r} \\
&= \left\{ \left(\frac{V_v}{(\underline{\mathfrak{m}_{i-1}})}\right)_{\overline{(\underline{\mathfrak{m}_{j-1}})}} \hookleftarrow \left(\frac{V_v}{(\underline{\mathfrak{m}_{i-1}})}\right)_{\overline{(\underline{\mathfrak{m}_j})}} \twoheadrightarrow \left(\frac{V_v}{(\underline{\mathfrak{m}_i})}\right)_{\overline{(\underline{\mathfrak{m}_j})}} \right\}_{1\leq i\leq j\leq r}
\end{aligned}
\tag{47}
$$

Here at the entry $(i,j) = (i, i-1)$ $(1 \leq i \leq r+1$, we assign

$$
\begin{cases} K & \text{if } i = 1 \\ \kappa(v_{i-1}) := \left(\frac{V_v}{(\underline{\mathfrak{m}_{i-1}})}\right)_{\overline{(\underline{\mathfrak{m}_{i-1}})}} & \text{if } 2 \leq i \leq r+1 \end{cases}
$$

More diagrammatically,

(48)
$$
\begin{array}{ccccccccccc}
K & \hookleftarrow \cdots \hookleftarrow & V_{v_1\circ\cdots v_{i-1}} & \hookleftarrow & V_{v_1\circ\cdots v_i} & \hookleftarrow & V_{v_1\circ\cdots v_{i+1}} & \hookleftarrow \cdots \hookleftarrow & V_{v_1\circ\cdots v_{rank(v)}} \\
 & & \downarrow\!\!\!\downarrow & & \downarrow\!\!\!\downarrow & & \downarrow\!\!\!\downarrow & & \downarrow\!\!\!\downarrow \\
 & \ddots & \vdots & \cdots & \vdots & \cdots & \vdots & \cdots & \vdots \\
 & & \downarrow\!\!\!\downarrow & & \downarrow\!\!\!\downarrow & & \downarrow\!\!\!\downarrow & & \downarrow\!\!\!\downarrow \\
 & & \kappa(v_{i-1}) & \hookleftarrow & V_{v_i} & \hookleftarrow & V_{v_i\circ v_{i+1}} & \cdots \hookleftarrow & V_{v_i\circ\cdots v_{rank(v)}} \\
 & & & & \downarrow\!\!\!\downarrow & & \downarrow\!\!\!\downarrow & & \downarrow\!\!\!\downarrow \\
 & & & & \kappa(v_i) & \hookleftarrow & V_{v_{i+1}} & \cdots \hookleftarrow & V_{v_{i+1}\circ\cdots v_{rank(v)}} \\
 & & & & & & \downarrow\!\!\!\downarrow & & \downarrow\!\!\!\downarrow \\
 & & & & & & \kappa(v_{i+1}) & \cdots \hookleftarrow & V_{v_{i+2}\circ\cdots v_{rank(v)}} \\
 & & & & & & & & \downarrow\!\!\!\downarrow \\
 & & & & & & & \ddots & \vdots \\
 & & & & & & & & \downarrow\!\!\!\downarrow \\
 & & & & & & & & \kappa(v_{rank(v)})
\end{array}
$$

is identified with

(49)
$$
\begin{array}{ccccccccc}
Frac\,(V_v) \hookleftarrow \cdots \hookleftarrow & (V_v)_{\mathfrak{m}_{i-1}} & \hookleftarrow & (V_v)_{\mathfrak{m}_i} & \hookleftarrow & (V_v)_{\mathfrak{m}_{i+1}} & \hookleftarrow \cdots \hookleftarrow & (V_v)_{\mathfrak{m}_{rank(v)}} = V_v \\
 & \downarrow\!\!\!\downarrow & & \downarrow\!\!\!\downarrow & & \downarrow\!\!\!\downarrow & & \downarrow\!\!\!\downarrow \\
\ddots & \vdots & \cdots & \vdots & \cdots & \vdots & \cdots & \vdots \\
 & \downarrow\!\!\!\downarrow & & \downarrow\!\!\!\downarrow & & \downarrow\!\!\!\downarrow & & \downarrow\!\!\!\downarrow \\
 & Frac\left(\frac{V_v}{(\mathfrak{m}_{i-1})}\right) & \hookleftarrow & \left(\frac{V_v}{(\mathfrak{m}_{i-1})}\right)_{\overline{(\mathfrak{m}_i)}} & \hookleftarrow & \left(\frac{V_v}{(\mathfrak{m}_{i-1})}\right)_{\overline{(\mathfrak{m}_{i+1})}} & \cdots \hookleftarrow & \left(\frac{V_v}{(\mathfrak{m}_{i-1})}\right)_{\overline{(\mathfrak{m}_{rank(v)})}} \\
 & & & \downarrow\!\!\!\downarrow & & \downarrow\!\!\!\downarrow & & \downarrow\!\!\!\downarrow \\
 & & & Frac\left(\frac{V_v}{(\mathfrak{m}_i)}\right) & \hookleftarrow & \left(\frac{V_v}{(\mathfrak{m}_i)}\right)_{\overline{(\mathfrak{m}_{i+1})}} & \cdots \hookleftarrow & \left(\frac{V_v}{(\mathfrak{m}_i)}\right)_{\overline{(\mathfrak{m}_{rank(v)})}} \\
 & & & & & \downarrow\!\!\!\downarrow & & \downarrow\!\!\!\downarrow \\
 & & & & & Frac\left(\frac{V_v}{(\mathfrak{m}_{i+1})}\right) & \cdots \hookleftarrow & \left(\frac{V_v}{(\mathfrak{m}_{i+1})}\right)_{\overline{(\mathfrak{m}_{rank(v)})}} \\
 & & & & & & & \downarrow\!\!\!\downarrow \\
 & & & & & & \ddots & \vdots \\
 & & & & & & & \downarrow\!\!\!\downarrow \\
 & & & & & & & Frac\left(\frac{V_v}{(\mathfrak{m}_{rank(v)})}\right)
\end{array}
$$

– generalizing and building upon earlier works of Zariski [Z40] and Abhyankar [A56b], Knaf and-Kuhlmann [KK05, p.835, Theorem 1.1] (see also Temkin [T13, p.114, Theorem 5.5.2] and Cutkosky [C21, Theorem 1.3]) proved the following strong form of local uniformization (43) theorem for Abhyankar valuations: [13]

[13] While Knaf and Kuhlmann [KK05, p.835,Theorem 1.1] (and also Temkin [T13, p.114, Theorem 5.5.2]) realized their model by $\mathcal{O}_{X,x}$ at a smooth point $x \in X$ of some variety when $\kappa(v)$ is separable over $k$, Cutkosky [C21, Theorem 1.1, Theorem 1.3] realized his model without any separable assumption by a regular local ring that is open-immersed-essentially of finite presentation over arbitrary field $k$ [C21, 2.Notation]. However, under our base field perfect [sp18, 030Y,030Z], such a separable assumption is irrelevant.

**Theorem 3.4.** [KK05, p.835, Theorem 1.1] ([T13, p.114, Theorem 5.5.2] [C21, Theorem 1.3], see also [NS14, p.405, Theorem 1.1]) *Suppose that $K$ is an algebraic function field over a field $k$ and $v = v_1 \circ \cdots \circ v_r$ ($r =$ rank $v$) is an Abhyankar valuation of $K/k$;*
*Then there exists a regular local ring $(R, \mathfrak{m}_R)$ such that:*

* $*$ *$R$ is a filtered colimit of smooth $k$-algebras of finite type;*
* $*$ *there is a local homomorphism $i_R : (R, \mathfrak{m}_R) \to (V_v, \mathfrak{m}_v)$ such that*

(50) $$Frac(R) \xrightarrow[\cong]{i_R} Frac(V_v) = K.$$

* $*$ *$(R, \mathfrak{m}_R)$ has a regular system of generators $\{x_{i,j}\}_{\substack{1\le i\le rank(v)\\ 1\le j\le rr_i}}$ such that for each $1 \le i \le j \le r$,*
  * $\cdot$ *$\{x_{a,b}\}_{\substack{i\le a\le j\\ 1\le b\le rr_a}} \subset i_R^{-1}\left(V_{v_1\circ\cdots\circ v_j}\right)$ and their images by the composite $\{x_{a,b}\}_{\substack{i\le a\le j\\ 1\le b\le rr_a}} \xrightarrow{i_R} V_{v_1\circ\cdots\circ v_j} \twoheadrightarrow V_{v_i\circ\cdots\circ v_j} \xrightarrow{v_i\circ\cdots\circ v_j} \Gamma_{v_i\circ\cdots\circ v_j}$ $= \Gamma_{i-1}/\Gamma_j \cong \mathbb{Z}^{\sum_{i\le a\le j} rr_a}$, where the last isomorphism is just as abelian groups, not necessarily as ordered commutative groups, form a $\mathbb{Z}$-basis;*
  * $\cdot$ *the set $\left\{\pi_{v_1\circ\cdots\circ v_{i-1}}\big|_R (x_{a,b})\right\}_{\substack{i\le a\le j\\ 1\le b\le rr_a}}$ forms a regular system of generators of $\left(\frac{R}{R\cap\mathfrak{m}_{v_1\circ\cdots\circ v_{i-1}}}\right)_{\overline{R\cap\mathfrak{m}_{v_1\circ\cdots\circ v_j}}} \left(\to V_{v_i\circ\cdots\circ v_j}\right)$.*
* $*$ *$i_R$ induces (in addition to (50) isomorphisms:*

(51) $$Frac\left(\frac{R}{R\cap\mathfrak{m}_{v_1\circ\cdots\circ v_i}}\right) \xrightarrow[\cong]{i_R} V_{v_i}/\mathfrak{m}_{v_i} = \kappa(v_i).$$

□

While this wonderful Theorem 3.4 is very deep (and so not surprisingly difficult to digest its proof [14]), we can relatively easily prove a sufficient result for our purpose in the framework of the geometric valuations, as we shall soon see in Theorem 3.6..

So, let us turn our attention to geometric valuations:

* The terminologies of <u>*divisorial divisor*</u> or <u>*prime divisor*</u> for a rank 1 geometric valuation $v$ are justified by the following classical result of Zariski [Z39] (slightly modified along the line of [NS14]):

**Theorem 3.5.** (Local uniformization of first kind with respect to a divisorial valuation)
(i) *Any divisorial valuation $v$ of a finitely generated field extension $K/k$ with the valuation ring $V_v$ and the maximal ideal $\mathfrak{m}_v$ can be "geometrically realized." More precisely, there exist a normal integral subdomain $R \subseteq V_v$ which is finitely generated over $k$ posessing $\mathfrak{m}_v \cap R \subset R$ as its height 1 prime ideal, with the corresponding integral affine variety and its prime divisor* Spec $R \supset V(\mathfrak{m}_v \cap R)$, *yielding the following commutative diagram with vertical equalities holds:*

[14] see e.g. [S17, The last paragraph in p.6]. However, consult [T14].

$$
\begin{array}{ccccc}
k(\mathrm{Spec} R) & \hookleftarrow & \mathcal{O}_{\mathrm{Spec} R, V(\mathfrak{m}_v \cap R)} & \twoheadrightarrow & k\left(V(\mathfrak{m}_v \cap R)\right) \\
\| & & \| & & \| \\
Frac(R) & \hookleftarrow & R_{\mathfrak{m}_v \cap R} & \twoheadrightarrow & Frac\left(R/\mathfrak{m}_v \cap R\right) \\
\| & & \| & & \| \\
K & \hookleftarrow & V_v & \twoheadrightarrow & \kappa(v)
\end{array} \tag{52}
$$

(ii) *Let $v$ be a divisorial valuation of a finitely generated field extension $K/k$. Then, for any finitely generated $k$-algebra $S$ such that*

$$S \subset V_v \subset K = Frac(S),$$

*there is a* local blowup of $S$ with respect to $v$, *i.e. an obviously defined local ring homomorphism*

$$
\begin{aligned}
&(S_{\mathfrak{m}_v \cap S}, \mathfrak{m}_v \cap S_{\mathfrak{m}_v \cap S}) \\
&\to \left( S\left[\frac{a_1}{b_1}, \ldots, \frac{a_s}{b_s}\right]_{\mathfrak{m}_v \cap S\left[\frac{a_1}{b_1}, \ldots, \frac{a_s}{b_s}\right]}, \mathfrak{m}_v \cap S\left[\frac{a_1}{b_1}, \ldots, \frac{a_s}{b_s}\right]_{\mathfrak{m}_v \cap S\left[\frac{a_1}{b_1}, \ldots, \frac{a_s}{b_s}\right]} \right)
\end{aligned}
$$

*for some $a_i, b_i \in S$ $(i = 1, \ldots, s)$ with $v(a_i) \geq v(b_i)$ $(i = 1, \ldots, s)$, such that the following so-called the* first kind *equality holds:*

$$
\begin{aligned}
&\left( S\left[\frac{a_1}{b_1}, \ldots, \frac{a_s}{b_s}\right]_{\mathfrak{m}_v \cap S\left[\frac{a_1}{b_1}, \ldots, \frac{a_s}{b_s}\right]}, \mathfrak{m}_v \cap S\left[\frac{a_1}{b_1}, \ldots, \frac{a_s}{b_s}\right]_{\mathfrak{m}_v \cap S\left[\frac{a_1}{b_1}, \ldots, \frac{a_s}{b_s}\right]} \right) \\
&= (V_v, \mathfrak{m}_v).
\end{aligned}
$$

*Furthermore, there exists a normal integral subdomain*

$$R \subset S\left[\frac{a_1}{b_1}, \ldots, \frac{a_s}{b_s}\right]_{\mathfrak{m}_v \cap S\left[\frac{a_1}{b_1}, \ldots, \frac{a_s}{b_s}\right]} = V_v,$$

*which is finitely generated over $k$ possessing $\mathfrak{m}_v \cap R \subset R$ as its height 1 prime ideal and realizes the commutative diagram* (52) *in* (i).

Although the proofs of Zariski's local uniformization of the first kind with respect to a divisorial valuation are given in literature like [ZS60b, p.89, VI, §14, Theorem 31] [A86b, p.230, Theorem 5.2] [KM98, p.61, Lemma 2.45] (see also [V06, p.495, Example 5]), we shall still present a proof of Theorem 3.5 below, because Theorem 3.5 is sligtly modified along the line of [NS14] so as to be generalized to the case of geometric valuation of arbitrary rank. For this purpose, we now recall some basic results of commutative ring theory :

- As in [M80, p.71–72, §12] [M89, p.30, §5], we immediately obtain the following equalities and and an inequality for any commutative ring $A$ and its prime ideal $\mathfrak{p} \subset A$ from the definitions:

$$ht\ \mathfrak{p} = \dim A_{\mathfrak{p}}, \quad coht\ \mathfrak{p} = \dim A/\mathfrak{p} \quad \text{and} \quad ht\ \mathfrak{p} + coht\ \mathfrak{p} \leq \dim A. \tag{53}$$

If $A$ is an integral domain, which is finitely generated over a field $k$, then the above naive inequality (53) turns out to be an equality as (54) below, for any prime ideal $\mathfrak{p} \subset A$ [S00, p.45, Propsition 15] [M80, p.92, (14.H) Corollary 3] [M89, p.137, Theorem 17.7, Theorem 17.9]:

(54) $$ht\ \mathfrak{p} + \dim A/\mathfrak{p} = \dim A.$$

– [S00, p.44, Propsition 14] [M80, p.91, (14.G) Corollary 1] [M89, p.34, Theorem 5.6] For an integral domain $A$ which is finitely generated over a field $k$,

(55) $$\dim A = tr.deg_k Frac(A).$$

– [AM69, p.61, Proposition 5.6 ii)][M89, p.65, Example 3] [sp18, 0307] "Integral closure commutes with localization": Denote by $\widetilde{(-)}$ the integral closure of a ring $(-)$ in its given ring extension. Then, for $S \subset A \hookrightarrow B$ a multiplicative closed subset $S$ of an integral domain $A$ and an injective ring homomorphism to $B$,

(56) $$\widetilde{S^{-1}A} = S^{-1}\widetilde{A}.$$

– [AM69, p.66, Corollary 5.22] [M89, p.73, Theorem 10.4] Let $A$ be a subring of a field $K$. Then the integral closure $\widetilde{A}$ of $A$ in $K$ is the intersection of all the valuation rings of $K$ which contain $A$.
– [ZS60a, p.267, Chap V, Theorem 9] [S00, p.46, Proposition 16] [M89, p.262–263] For an integral domain $A$ which is a finitely generated algebra over a field $k$, its integral closure $\widetilde{A}$ in a finite algebraic extension of $Frac(A)$ is a finitely generated $A$-module and is also a finitely generated $k$-algebra.

*Proof of Theorem 3.5.*

**Reduction of (i) to (ii):** Since $K/k$ is finitely generated, $K$ is generated by some $x_i \in K \setminus \{0\}$ $(1 \leq i \leq m)$, and since $Frac(V_v) = K$, there are some $n_i, d_i \in V_v \setminus \{0\}$ $(1 \leq i \leq m)$ such that $x_i = \frac{n_i}{d_i}$ $(1 \leq i \leq m)$. Now, set $S$ be the $k$-subalgebra of $K$, generated by $\{n_i, d_i \in V_v \setminus \{0\} \mid 1 \leq i \leq m\}$, then clearly,

$$S \subset V_v \subset K = Frac(S),$$

which is the situation of (ii), as desired.

So, we now concentrate in proving (ii):

**Construction of a weakly increasing sequence: $S = S^{(0)} \subseteq S^{(1)} \subseteq \cdots \subseteq V_v$ :** We first set $S^{(0)} = S \subseteq V_v$, and we inductively define $S = S^{(0)} \subseteq \left(S^{(1)}, \mathfrak{m}^{(1)}\right) \to \left(S^{(2)}, \mathfrak{m}^{(2)}\right) \to \cdots \to \left(S^{(n)}, \mathfrak{m}^{(n)}\right) \to \left(S^{(n+1)}, \mathfrak{m}^{(n+1)}\right) \to$

$\cdots \subseteq V_v$ so that

(57)
$$\begin{cases} S^{(n+1)} := \left\{ \frac{x}{y} \in F \;\middle|\; x \in S^{(n)}\left[\frac{u_1}{u_0}, \cdots, \frac{u_s}{u_0}\right], \right. \\ \qquad\qquad\qquad \left. y \in S^{(n)}\left[\frac{u_1}{u_0}, \cdots, \frac{u_s}{u_0}\right] \setminus \left(\mathfrak{m}_v \cap S^{(n)}\left[\frac{u_1}{u_0}, \cdots, \frac{u_s}{u_0}\right]\right) \right\} \\ \quad \text{with } u_i \in \mathfrak{m}^{(n)} = \mathfrak{m}_v \cap S^{(n)} \ (0 \le i \le s) \text{ s.t.} \\ \qquad v(u_0) = \min\{v(u) \mid u \in \mathfrak{m}_v \cap S^{(n)}\} \text{ and } u_i \in \mathfrak{m}^{(n)} = \mathfrak{m}_v \cap S^{(n)} \\ \qquad (0 \le i \le s) \text{ generate } \mathfrak{m}^{(n)} = \mathfrak{m}_v \cap S^{(n)} \text{ as an } S\text{-module.} \\ \mathfrak{m}^{(n+1)} := \mathfrak{m}_v \cap S^{(n+1)} = \{x \in S^{(n+1)} \mid v(x) > 0\}, \\ \qquad \text{and so } S^{(n+1)}_{\mathfrak{m}^{(n+1)}} = S^{(n+1)}. \end{cases}$$

i.e. the local blowing up of $S^{(n)}$ with respect to $v$ along $\mathfrak{m}^{(n)} = \mathfrak{m}_v \cap S^{(n)}$ as in (42). This gives us the desired increasing sequence of Noetherian [15] (and also local, except possibly $S = S^{(0)}$) rings.

$$S = S^{(0)} \subseteq S^{(1)} \subseteq \cdots \subseteq S^{(n)} \subseteq S^{(n+1)} \subseteq \cdots \subseteq V_v.$$

**Proof of $\cup_n S^{(n)} = V_v$:** For $n \in \mathbb{N}$, suppose there exists $\alpha \in V_v \setminus S^{(n)}$.

- Then, there exist $a, b \in S^{(n)}$ such that $\alpha = \frac{a}{b}$, because $\alpha \in V_v \setminus S^{(n)} \subset F = Frac\left(S^{(n)}\right)$.
- Here, $b \in \mathfrak{m}^{(n)}$, because, otherwise, $b$ would belong to $S^{(n)} \setminus \mathfrak{m}^{(n)}$, and so, $\alpha = \frac{a}{b} \in \left(S^{(n)}\right)_{\mathfrak{m}^{(n)}} = S^{(n)}$ : a contradiction to our choice of $\frac{a}{b} = \alpha \in V_v \setminus S^{(n)}$.
- $a \in \mathfrak{m}^{(n)} = \{x \in S^{(n)} \mid v(x) > 0\}$ too, because $0 \overset{\alpha \in V_v}{\le} v(\alpha) = v(a) - v(b)$ and $b \in \mathfrak{m}^{(n)} = \{x \in S^{(n)} \mid v(x) > 0\}$.
- Then, there exist $a_1, b_1 \in S^{(n+1)}$ such that $\alpha = \frac{a}{b} = \frac{a_1}{b_1}$, such that $0 \le v(b_1) < v(b)$, because:
  * By (57), we may express $b, a \in \mathfrak{m}^{(n)} = m_v \cap S^{(n)}$ as:
  $$b = \sum_{0 \le i \le s} s_i'' u_i, \ a = \sum_{0 \le i \le s} s_i' u_i \text{ for some } s_i', s_i'' \in S^{(n)} \ (0 \le i \le s).$$
  * If we set
  $$b_1 = \frac{b}{u_0} \sum_{0 \le i \le s} s_i'' \frac{u_i}{u_0}, \ a_1 = \sum_{0 \le i \le s} s_i' \frac{u_i}{u_0} \ \in \ S^{(n+1)},$$
  we clearly have $\alpha = \frac{a}{b} = \frac{a_1}{b_1}$ and $0 \le v(b_1) = v(b) - v(u_0) < v(b)$, as desired.
- We may iterate the above construction inductively from the $i$-th step to $(i+1)$-th step, as far as $\alpha = \frac{a}{b} = \frac{a_1}{b_2} = \cdots = \frac{a_i}{b_i} \notin S^{(n+i)}$; then we may find some $a_{i+1}, b_{i+1} \in S^{(n+i+1)}$ such that
  $$\frac{a_{i+1}}{b_{i+1}} = \frac{a_i}{b_i} = \cdots \frac{a_1}{b_2} = \frac{a}{b} = \alpha$$
  with $0 \le v(b_{i+1}) < v(b_i) < \cdots < v(b_2) < v(b_1) < v(b)$.

[15]Since $S$ is Noetherian by assumption, this follows immediately by applying [AM69, p.80, Proposition 7.3, p.81, Corollary 7.7] [sp18, 00FN] inductively.

– Since the above iteration can not last forever (i.e. not possible for $i > v(b)$), there exists some $j \in \mathbb{N}$ such that $v(b_j) = 0$, and so $\alpha = \frac{a}{b} = \frac{a_1}{b_2} = \cdots = \frac{a_j}{b_j} \in S^{(n+j)}_{\mathfrak{m}^{(n+j)}} = S^{(n+j)} \subseteq \cup_n S^{(n)}$, as desired.

$\exists n_0 \in \mathbb{N}$ **s.t.** $S^{(n_0)} \subseteq \forall$**ring** $B \subseteq V_v,\ tr.deg_k\,(B/(\mathfrak{m}_v \cap B)) = tr.deg_k\,(V_v/\mathfrak{m}_v)$ **:**

Since $tr.deg_k\,(V_v/\mathfrak{m}_v) = tr.deg_k \kappa(v) \overset{(24)}{=} \dim(v) \overset{v:\text{ divisorial valuation}}{=} tr.deg_k K - 1 < \infty$, there exist finitely many transcendence basis of $V_v/\mathfrak{m}_v$ over $k$. Then, as $V_v = \cup_n S^{(n)}$, we only have to choose $n_0$ large enough so that $S^{(n_0)}/(S^{(n_0)} \cap \mathfrak{m}_v)$ contains these finitely many transcendence basis.

**If further** $B = A_{\mathfrak{m}_v \cap A}$ **for a fin.gen.** $k$**-subalgebra** $A \subset K$, **then** $\dim B = 1$**:**

This simply follows from:

$$\begin{aligned}
\dim B = \dim A_{\mathfrak{m}_v \cap A} &\overset{(53)(54)}{=} \dim A - \dim A/(\mathfrak{m}_v \cap A) \\
&\overset{(55)}{=} tr.deg_k Frac\,(A) - tr.deg_k Frac\,(A/(\mathfrak{m}_v \cap A)) \\
&= tr.deg_k K - tr.deg_k\,(B/(\mathfrak{m}_v \cap B)) = tr.deg_k K - tr.deg_k\,(V_v/\mathfrak{m}_v) \\
&= tr.deg_k K - tr.deg_k \kappa(v) \overset{v:\text{ divisorial}}{=} 1.
\end{aligned}$$

**For the integral closure** $\widetilde{S^{(n_0)}}$ **of** $S^{(n_0)}$ **in** $K$, $\widetilde{S^{(n_0)}} \subseteq V_v, \dim\left(\widetilde{S^{(n_0)}}\right)_{\mathfrak{m}_v \cap \widetilde{S^{(n_0)}}} = 1$ **:**

– By [AM69, p.66, Corollary 5.22][M89, p.73, Theorem 10.4], $\widetilde{S^{(n_0)}} \subseteq \bigcap_{\substack{V:\text{valuation ring} \\ \text{of } K \text{ s.t. } S^{(n_0)} \subseteq V}} V \subseteq V_v$, which proves the first claim.
– For each $n \in \mathbb{N}$, $S^{(n)} = (A_n)_{\mathfrak{m}_v \cap A_n}$ ($n \in \mathbb{N}$) for some finitely generated $k$-subalgebra $A_n \subset K$, because:
  * $\{S^{(n)}\}_{n \in \mathbb{Z}_{\geqq 0}}$ is inductively constructed (57) by iterating local blowing ups with respect to $v$, and such an iteration can be realized by a single local blowing up of $S$ by Proposition 3.2.
  * Thus, we can express $S^{(n)}$ as (36):

$$S^{(n)} = \left(S\left[\frac{a_1}{b_1}, \cdots, \frac{a_s}{b_s}\right]\right)_{\mathfrak{m}_v \cap S\left[\frac{a_1}{b_1}, \cdots, \frac{a_s}{b_s}\right]}$$

  for some $a_i, b_i \in S, i = 1, \ldots, r$ such that $v(b_i) \leq v(a_i)$ for all $i = 1, \ldots, r$.
  * Since $S$ is a finitely generated $k$-algebra by assumption, we only have to set $A_n := S\left[\frac{a_1}{b_1}, \cdots, \frac{a_s}{b_s}\right]$ $(\subset K)$.

- $\left(\widetilde{S^{(n_0)}}\right)_{\mathfrak{m}_v\cap\widetilde{S^{(n_0)}}} = \widetilde{(A_{n_0})}_{\mathfrak{m}_v\cap\widetilde{(A_{n_0})}}$, because:

$$\begin{aligned}
\left(\widetilde{S^{(n_0)}}\right)_{\mathfrak{m}_v\cap\widetilde{S^{(n_0)}}} &= \left(\overbrace{(A_{n_0})_{\mathfrak{m}_v\cap A_{n_0}}}\right)_{\mathfrak{m}_v\cap\widetilde{S^{(n_0)}}} \\
&= \left(\overbrace{(A_{n_0}\setminus(\mathfrak{m}_v\cap A_{n_0}))^{-1}(A_{n_0})}\right)_{\mathfrak{m}_v\cap\widetilde{S^{(n_0)}}} \\
&\overset{(56)}{=} \left(\underbrace{(A_{n_0}\setminus(\mathfrak{m}_v\cap A_{n_0}))^{-1}\widetilde{(A_{n_0})}}_{=\widetilde{S^{(n_0)}}}\right)_{\mathfrak{m}_v\cap\widetilde{S^{(n_0)}}} \\
&= \left(\widetilde{S^{(n_0)}}\setminus\left(\mathfrak{m}_v\cap\widetilde{S^{(n_0)}}\right)\right)^{-1}\left(\underbrace{(A_{n_0}\setminus(\mathfrak{m}_v\cap A_{n_0}))^{-1}\widetilde{(A_{n_0})}}_{=\widetilde{S^{(n_0)}}}\right) \\
&= \left(\widetilde{(A_{n_0})}\setminus\left(\mathfrak{m}_v\cap\widetilde{(A_{n_0})}\right)\right)^{-1}\widetilde{(A_{n_0})} = \widetilde{(A_{n_0})}_{\mathfrak{m}_v\cap\widetilde{(A_{n_0})}}
\end{aligned}$$

- Since $A_{n_0}$ is a finitely generated $K$-algebra, $\widetilde{(A_{n_0})}$ is also a finitely generated $k$-algebra by [ZS60a, p.267, Chap V, Theorem 9] [S00, p.46, Proposition 16] [M89, p.262–263] recalled earlier.
- Now that $S^{(n_0)} \subseteq \left(\widetilde{S^{(n_0)}}\right)_{\mathfrak{m}_v\cap\widetilde{S^{(n_0)}}} = \widetilde{(A_{n_0})}_{\mathfrak{m}_v\cap\widetilde{(A_{n_0})}} \subseteq V_v$ for a finitely generated $k$-subalgebra $\widetilde{(A_{n_0})} \subseteq K$, the previous discussions imply $\dim\left(\widetilde{S^{(n_0)}}\right)_{\mathfrak{m}_v\cap\widetilde{S^{(n_0)}}} = 1$, as was desired.

<u>$\left(\widetilde{S^{(n_0)}}\right)_{\mathfrak{m}_v\cap\widetilde{S^{(n_0)}}} = \widetilde{(A_{n_0})}_{\mathfrak{m}_v\cap\widetilde{(A_{n_0})}}$ **is a discrete valuation ring.**</u> **:** In fact,

- $\widetilde{(A_{n_0})}_{\mathfrak{m}_v\cap\widetilde{(A_{n_0})}}$ is Noetherian because of the Noetherian property of $\widetilde{(A_{n_0})}$, by applying [AM69, p.80, Proposition 7.3] [sp18, 00FN].
- $\widetilde{(A_{n_0})}_{\mathfrak{m}_v\cap\widetilde{(A_{n_0})}}$ is normal because of the normality of $\widetilde{(A_{n_0})}$, by applying (56) "integral closure commutes with localization ":

$$\begin{aligned}
\overbrace{\widetilde{(A_{n_0})}_{\mathfrak{m}_v\cap\widetilde{(A_{n_0})}}} &= \overbrace{\left(\widetilde{(A_{n_0})}\setminus\left(\mathfrak{m}_v\cap\widetilde{(A_{n_0})}\right)\right)^{-1}\widetilde{(A_{n_0})}} \\
&\overset{(56)}{=} \left(\widetilde{(A_{n_0})}\setminus\left(\mathfrak{m}_v\cap\widetilde{(A_{n_0})}\right)\right)^{-1}\overbrace{\widetilde{(A_{n_0})}} \\
&= \left(\widetilde{(A_{n_0})}\setminus\left(\mathfrak{m}_v\cap\widetilde{(A_{n_0})}\right)\right)^{-1}\widetilde{(A_{n_0})} = \widetilde{(A_{n_0})}_{\mathfrak{m}_v\cap\widetilde{(A_{n_0})}}
\end{aligned}$$

- Now that $\left(\widetilde{S^{(n_0)}}\right)_{\mathfrak{m}_v\cap\widetilde{S^{(n_0)}}} = \widetilde{(A_{n_0})}_{\mathfrak{m}_v\cap\widetilde{(A_{n_0})}}$ is a Noetherian local normal domain of dim 1, it is a discrete valuation ring by [S00, p.37, Proposition 8] [AM69, p.94, Proposition 9.2] [M89, p.79, Theorem 11.2], as desired.

$\left(\widetilde{S^{(n_0)}}\right)_{\mathfrak{m}_v \cap \widetilde{S^{(n_0)}}} = \widetilde{(A_{n_0})}_{\mathfrak{m}_v \cap \widetilde{(A_{n_0})}} = V_v$ **: In fact,**

- Obviously, we have an inclusion of discrete valuation rings:

$$\left(\widetilde{S^{(n_0)}}\right)_{\mathfrak{m}_v \cap \widetilde{S^{(n_0)}}} = \widetilde{(A_{n_0})}_{\mathfrak{m}_v \cap \widetilde{(A_{n_0})}} \subseteq V_v \subseteq K = Frac\left(\widetilde{(A_{n_0})}_{\mathfrak{m}_v \cap \widetilde{(A_{n_0})}}\right)$$

- Thus, any $r \in V_v = \{r \in K \mid v(r) \geq 0\}$ can be expressed as $r = \frac{u_1 t^{p_1}}{u_2 t^{p_2}}$, where $t$ is a uniformizer of the discrete valuation ring $\widetilde{(A_{n_0})}$ and $u_i \in \widetilde{(A_{n_0})} \setminus \left(\mathfrak{m}_v \cap \widetilde{(A_{n_0})}\right)$, $p_i \in \mathbb{Z}_{\geq 0}$ for $i = 1.2$.
- $v(r) = (p_1 - p_2)v(t)$ with $v(t) > 0$, because

$$v(r) = v\left(\frac{u_1 t^{p_1}}{u_2 t^{p_2}}\right) = v(u_1) - v(u_2) + (p_1 - p_2)v(t) \quad \text{and}$$

$$t \in \mathfrak{m}_v \cap \left(\widetilde{(A_{n_0})}_{\mathfrak{m}_v \cap \widetilde{(A_{n_0})}}\right) \subseteq \mathfrak{m}_v = \{k \in K \mid v(k) > 0\},$$

$$\begin{aligned} u_i &\in \left(\widetilde{(A_{n_0})}_{\mathfrak{m}_v \cap \widetilde{(A_{n_0})}}\right) \setminus \left(\mathfrak{m}_v \cap \left(\widetilde{(A_{n_0})}_{\mathfrak{m}_v \cap \widetilde{(A_{n_0})}}\right)\right) \\ &\subseteq R_v \setminus (\mathfrak{m}_v \cap R_v) = \{k \in K \mid v(k) = 0\} \end{aligned}$$

- We now see $p_1 \geq p_2$ becuase $0 \leq v(r) = (p_1 - p_2)v(t)$ with $v(t) > 0$.
- Consequently, $r = \frac{u_1 t^{p_1}}{u_2 t^{p_2}} = \frac{u_1}{u_2} t^{p_1 - p_2} \in \widetilde{(A_{n_0})}_{\mathfrak{m}_v \cap \widetilde{(A_{n_0})}}$, as desired.

**For large enough $N \gg n_0$, $S^{(N)} = V_v$. :**

- As we mentioned before, $\widetilde{(A_{n_0})}$ is a finite $A_{n_0}$-module by [ZS60a, p.267, Chap V, Theorem 9][S00, p.46, Proposition 16] [M89, p.262–263]. Since $A_{n_0} \subseteq S^{(n_0)}$ and $\widetilde{(A_{n_0})} \subseteq V_v = \cup_n S^{(n)}$, this means

$$\widetilde{(A_{n_0})} \subseteq S^{(N)} \ (\subseteq V_v) \quad \text{for } N \gg n_0.$$

- Consequently,

$$\begin{array}{ccccc} V_v = \widetilde{(A_{n_0})}_{\mathfrak{m}_v \cap \widetilde{(A_{n_0})}} & \hookrightarrow & \left(S^{(N)}\right)_{\mathfrak{m}_v \cap S^{(N)}} & \hookrightarrow & (V_v)_{\mathfrak{m}_v \cap V_v}, \\ & & \| & & \| \\ & & S^{(N)} & \hookrightarrow & V_v \end{array}$$

which proves the claim: $S^{(N)} = V_v$.

**Completion of the proof:**

By the above, $\left(S^{(N)}, \mathfrak{m}_v \cap S^{(N)}\right) = (V_v, \mathfrak{m}_v)$. On the other hand, $\left(S^{(N)}, \mathfrak{m}_v \cap S^{(N)}\right)$ was constructed, out of $S$, by repeatedly applying the local uniformization as in (57), for which we can realize by a single blowing up by Proposition 3.2. Now, $S^{(N)} = V_v = \widetilde{(A_{n_0})}_{\mathfrak{m}_v \cap \widetilde{(A_{n_0})}}$, where $\widetilde{(A_{n_0})} \subset V_v$ is a normal integral subdomain, which is finitely generated over $k$ with $\mathfrak{m}_v \cap \widetilde{(A_{n_0})} \subset \widetilde{(A_{n_0})}$ its height 1 prime ideal.
Thus, we have finished the proof.

□

- For geometric valuations, the unique composition (32)

$$v = v_1 \circ v_2 \circ \cdots \circ v_{rank(v)},$$

  is given so that each $v_i$ $(1 \leq i \leq rank(v))$ is a rank 1 geometric valuation, a.k.a. a divisorial valuation, of a finitely generated field extension $\kappa(v_{i-1})/k$ with $V_{v_i} \subsetneq \kappa(v_{i-1})$ a discrete valuation ring.
- For this restricted class of geometric variations, we now prove the following strong form of Theorem 3.4 (Although Theorem 3.4 appears to be slightly weaker than the following Theorem 3.6, Theorem 3.4 is of course valid for the much wider class of Abhyankar valuations):

**Theorem 3.6.** (i) *Any rank $r$ geometric valuation $v = v_1 \circ \cdots \circ v_r$ of a finitely generated field extension $K/k$ with the valuation ring $V_v = V_{v_1 \circ \cdots \circ v_r}$, given by $r$ divisorial valuations $v_i$ on $\kappa(v_{i-1})/k$ $(1 \leq i \leq r)$ (where $\kappa(v_0)$ is interpreted to be $K$), can be "geometrically realized."*

*More precisely, there exists a normal integral subdomain $R \subseteq V_v$ with a naturally induced map of systems of ring homomorphisms for $1 \leq i \leq j \leq r$, which we may express via the simplified notations presented in* (29) *as follows:*

$$\begin{array}{ccccc}
\dfrac{\left(\dfrac{R}{\left(\mathfrak{m}_{v_1 \circ \cdots \circ v_{i-1}}\right)}\right)}{\left(\mathfrak{m}_{v_1 \circ \cdots \circ v_{j-1}}\right)} & \hookleftarrow & \dfrac{\left(\dfrac{R}{\left(\mathfrak{m}_{v_1 \circ \cdots \circ v_{i-1}}\right)}\right)}{\left(\mathfrak{m}_{v_1 \circ \cdots \circ v_j}\right)} & \twoheadrightarrow & \dfrac{\left(\dfrac{R}{\left(\mathfrak{m}_{v_1 \circ \cdots \circ v_i}\right)}\right)}{\left(\mathfrak{m}_{v_1 \circ \cdots \circ v_j}\right)}, \\
\downarrow & & \downarrow & & \downarrow \\
V_{v_i \circ \cdots \circ v_{j-1}} & \hookleftarrow & V_{v_i \circ \cdots \circ v_j} & \twoheadrightarrow & V_{v_{i+1} \circ \cdots \circ v_j}
\end{array}$$

*or diagrammatically, a map from the system of local rings*

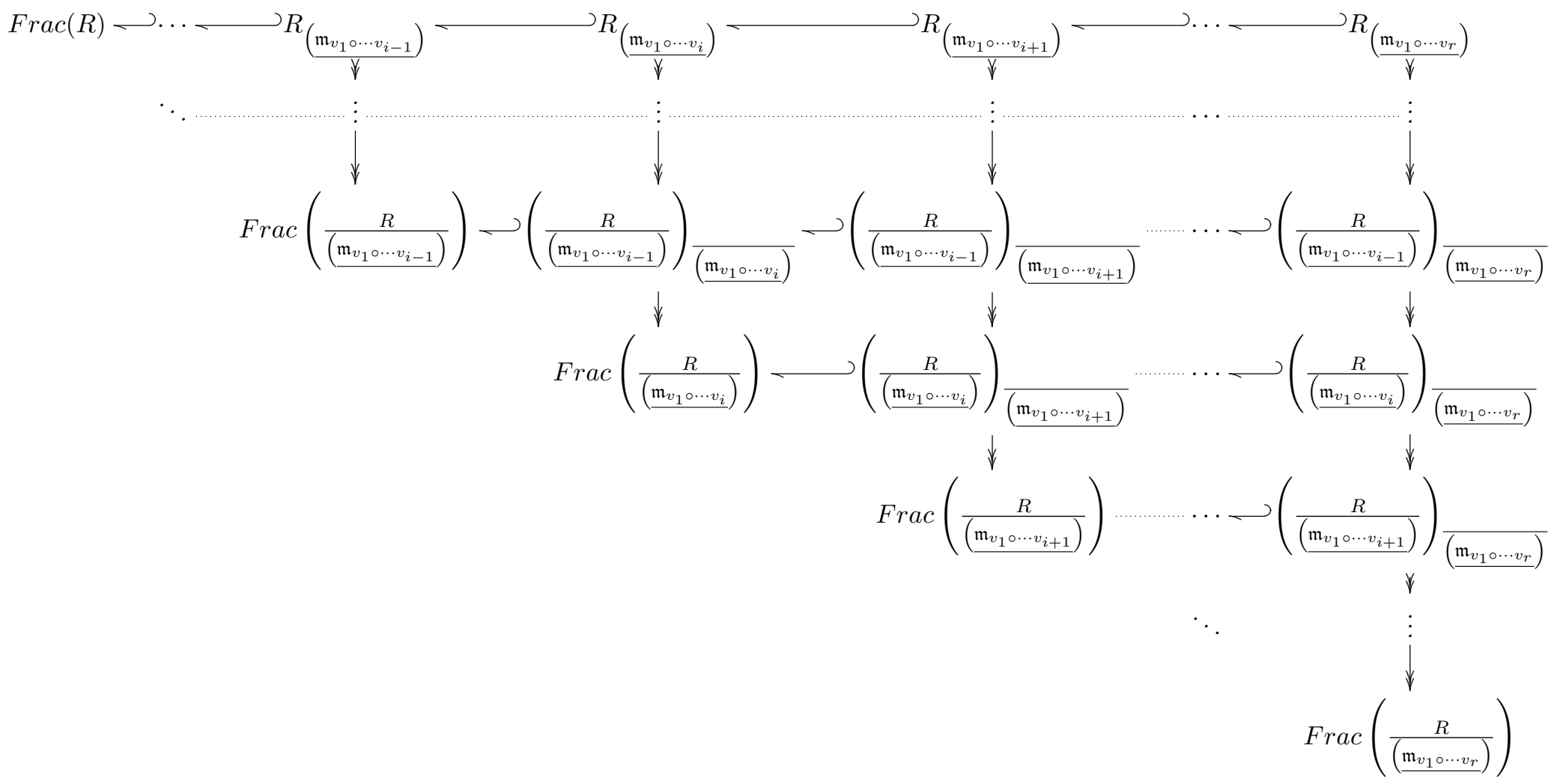

*to the system of local rings*

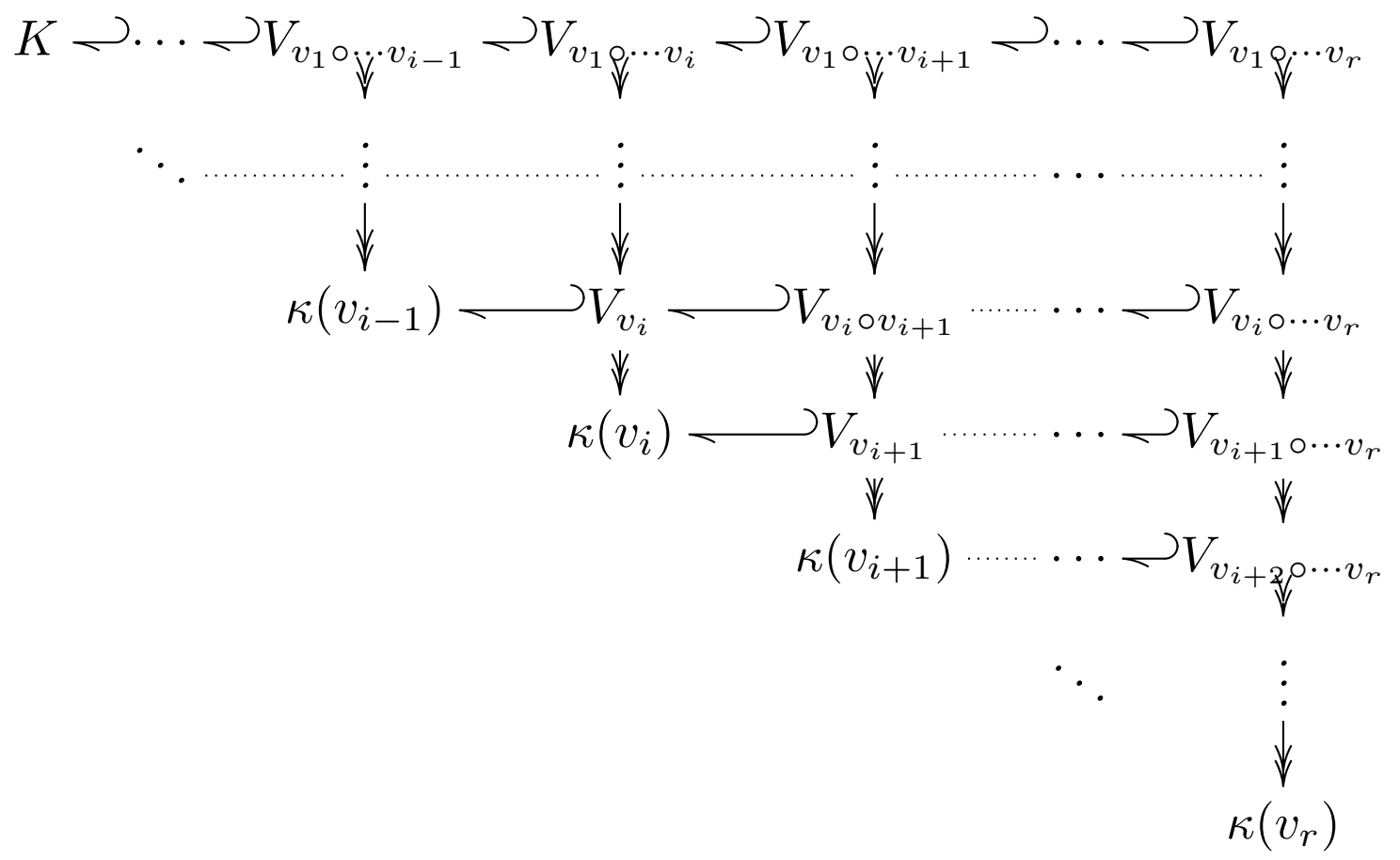

*such that the following conditions are satisfied:*

– $R$ *is finitely generated over* $k$. *In particular, for any* $1 \leq i \leq j \leq r$, $\left(\frac{R}{\left(\mathfrak{m}_{v_1\circ\cdots\circ v_{i-1}}\right)}\right)_{\overline{\left(\mathfrak{m}_{v_1\circ\cdots\circ v_j}\right)}}$ *is Noetherian.*

– $R_{\overline{\left(\mathfrak{m}_{v_1\circ\cdots\circ v_r}\right)}} = R_{R\cap\mathfrak{m}_{v_1\circ\cdots\circ v_r}}$ *is a regular local ring with a regular system of generators* $x_1, \ldots, x_r$ *such that for each* $1 \leq i \leq j \leq r$, :

  - $\{x_1, \ldots, x_{i-1}\} \in R \cap \mathfrak{m}_{v_1\circ\cdots\circ v_{i-1}} = \mathrm{Ker}\ \pi_{v_1\circ\cdots\circ v_{i-1}}\big|_R$,
  - *the set* $\left\{(v_1 \circ \cdots \circ v_{i-1})\big|_R(x_i), \ldots, (v_1 \circ \cdots \circ v_{i-1})\big|_R(x_j)\right\}$ *forms a* $\mathbb{Z}$*-basis of* $\Gamma_{i-1}/\Gamma_j \cong \mathbb{Z}^{j-i+1}$.
  - *the set* $\left\{\pi_{v_1\circ\cdots\circ v_{i-1}}\big|_R(x_i), \ldots, \pi_{v_1\circ\cdots\circ v_{i-1}}\big|_R(x_j)\right\}$ *forms a regular system of generators of* $\left(\frac{R}{\left(\mathfrak{m}_{v_1\circ\cdots\circ v_{i-1}}\right)}\right)_{\overline{\left(\mathfrak{m}_{v_1\circ\cdots\circ v_j}\right)}}$.

*Consequently, for any* $1 \leq i \leq j \leq r$, $\left(\frac{R}{\left(\underline{\mathfrak{m}_{v_1\circ\cdots\circ v_{i-1}}}\right)}\right)_{\overline{\left(\underline{\mathfrak{m}_{v_1\circ\cdots\circ v_j}}\right)}}$ *is a regular local ring. Furthermore, regarding the prime ideals* $\underline{\mathfrak{m}_{v_1\circ\cdots v_j}} = \mathfrak{m}_{v_1\circ\cdots v_j} \cap R$ $(1 \leq j \leq r)$ *as schematic points of* Spec$R$*, we therefore find they all belong to the smooth locus* $L$ *of* Spec$R$*, where* $L$ *is smooth and finite type over* $k$.[16] *Consequently, we have equivalences of local rings:*

$$R_{\left(\underline{\mathfrak{m}_{v_1\circ\cdots v_j}}\right)} = R_{\mathfrak{m}_{v_1\circ\cdots v_j}\cap R} \overset{\cong}{\longrightarrow} \mathcal{O}_{L,\mathfrak{m}_{v_1\circ\cdots v_j}\cap R} \qquad (1 \leq j \leq r). \tag{58}$$

– *The above map between systems of local rings restricts to equalities on the diagonal and the next to the diagonal part (where the equalities are commonly called of the* first kind *):*

$$\begin{array}{ccccc} Frac\left(\frac{R}{\left(\underline{\mathfrak{m}_{v_1\circ\cdots\circ v_{i-1}}}\right)}\right) & \hookleftarrow & \left(\frac{R}{\left(\underline{\mathfrak{m}_{v_1\circ\cdots\circ v_{i-1}}}\right)}\right)_{\overline{\left(\underline{\mathfrak{m}_{v_1\circ\cdots\circ v_i}}\right)}} & \twoheadrightarrow & Frac\left(\frac{R}{\left(\underline{\mathfrak{m}_{v_1\circ\cdots\circ v_i}}\right)}\right) \\ \| & & \| & & \| \\ \kappa(v_{i-1}) & \hookleftarrow & V_{v_i} & \twoheadrightarrow & \kappa(v_i) \end{array} \tag{59}$$

(ii) *Let* $v = v_1 \circ \cdots \circ v_r$ *be a rank* $r$ *geometric valuation of a finitely generated field extension* $K/k$*, given by* $r$ *divisorial valuations* $v_i$ *on* $\kappa(v_{i-1})/k$ $(1 \leq i \leq r)$ *(where* $\kappa(v_0)$ *is interpreted to be* $K$*), as in (i).*

*Then, for any finitely generated* $k$*-algebra* $S$ *such that*

$$S \subset V_{v_1\circ\cdots\circ v_r} \subset K = Frac(S),$$

*there is a* local blowup of $S$ with respect to $v_1 \circ \cdots \circ v_r$, *i.e. an obviously defined local ring homomorphisms*

$$\left(S_{\mathfrak{m}_{v_1\circ\cdots\circ v_r}\cap S}, \mathfrak{m}_{v_1\circ\cdots\circ v_r} \cap S_{\mathfrak{m}_{v_1\circ\cdots\circ v_r}\cap S}\right)$$
$$\to \left(S\left[\frac{a_1}{b_1}, \ldots, \frac{a_r}{b_r}\right]_{\mathfrak{m}_v\cap S\left[\frac{a_1}{b_1},\ldots,\frac{a_r}{b_r}\right]}, \mathfrak{m}_{v_1\circ\cdots\circ v_r} \cap S\left[\frac{a_1}{b_1}, \ldots, \frac{a_r}{b_r}\right]_{\mathfrak{m}_v\cap S\left[\frac{a_1}{b_1},\ldots,\frac{a_r}{b_r}\right]}\right)$$
$$\to \left(V_{v_1\circ\cdots\circ v_r}, \mathfrak{m}_{v_1\circ\cdots\circ v_r}\right)$$

*for some* $a_i, b_i \in S$ $(i = 1, \ldots, r)$ *with*

$$(v_1 \circ \cdots \circ v_r)(a_i) \geq (v_1 \circ \cdots \circ v_r)(b_i) \quad (i = 1, \ldots, r),$$

*such that the following conclusions hold:*

(a) $S\left[\frac{a_1}{b_1}, \ldots, \frac{a_r}{b_r}\right]_{\mathfrak{m}_v\cap S\left[\frac{a_1}{b_1},\ldots,\frac{a_r}{b_r}\right]}$ *contains a normal integral subdomain*

$$R_f \subset S\left[\frac{a_1}{b_1}, \ldots, \frac{a_r}{b_r}\right]_{\mathfrak{m}_v\cap S\left[\frac{a_1}{b_1},\ldots,\frac{a_r}{b_r}\right]} \subset V_{v_1\circ\cdots\circ v_r},$$

*which is finitely generated over* $k$ *with* $\mathfrak{m}_{v_1\circ\cdots\circ v_r} \cap R_f \subset R_f$ *its height* $r$ *prime ideal.*

(b) *For* $R = R_f$, $S\left[\frac{a_1}{b_1}, \ldots, \frac{a_r}{b_r}\right]_{\mathfrak{m}_v\cap S\left[\frac{a_1}{b_1},\ldots,\frac{a_r}{b_r}\right]}$, *properties* (58) *and* (59)*, stated in* (i)*, are satisfied.*

[16]Here, we have applied Proposition 1.8.

Restricting to the case of rank 2 geometric valuations, Theorem 3.6, especially (58) and (59), immediately implies the following:

**Corollary 3.7.** *Let $v = v_1 \circ v_2$ be a rank 2 geometric valuation of a finitely generated field extension $K/k$, given by a divisorial valuation $v_1$ on $K/k$ and a divisorial valuation $v_2$ on $\kappa(v_1)/k$.*

*Then there exist a normal integral subdomain $R \subseteq V_v$, which is finitely generated over the base field $k$, with the smooth locus $L$ of the affine scheme* Spec $R$*, containing the schematic points $R \cap \mathfrak{m}_{v_1}, R \cap \mathfrak{m}_{v_1 \circ v_2}$, such that the following commutative diagram, where all the morphisms except those connecting $V_{v_1 \circ v_2}$ are in $Pr\,(Sm_k)$, exists:*

(60)
$$\begin{array}{ccccc}
 & & \mathcal{O}_{L,\mathfrak{m}_{v_1}\cap R} & \hookleftarrow & \mathcal{O}_{L,\mathfrak{m}_{v_1\circ v_2}\cap R} \\
 & & \| & & \| \\
Frac(R) & \hookleftarrow & R_{\mathfrak{m}_{v_1}\cap R} & \hookleftarrow & R_{\mathfrak{m}_{v_1\circ v_2}\cap R} \\
\| & & \| & \swarrow & \downarrow \\
K & \hookleftarrow & V_{v_1} & \hookleftarrow V_{v_1\circ v_2} & \\
 & & \twoheaddownarrow & \twoheaddownarrow & \\
 & & \kappa(v_1) = Frac\left(\frac{R}{\mathfrak{m}_{v_1}\cap R}\right) & \hookleftarrow V_{v_2} = & \left(\frac{R}{\mathfrak{m}_{v_1}\cap R}\right)_{\left(\mathfrak{m}_{v_1\circ v_2}\cap R\right)}
\end{array}$$

□

**Remark 3.8.** Essentially, the only but the truly critical improvement of Theorem 3.6 over Theorem 3.4 in this special case of geometric valuations are the first kind equalities next to the diagonal (59), i.e. realizations by localizations of finitely generated algebras over the base field $k$ of: $V_v$ in (52) of Theorem 3.5; $V_{v_i}$ in (59) of Theorem 3.6; and $V_{v_1}$ and $V_{v_2}$ in (60) of Corollary 3.7.

We stress that our realization of $V_{v_i}$ in (59) of Theorem 3.6 does not hold for non-geometric Abhyankar valuation $v = v_1 \circ \cdots \circ v_r$. In fact, if an Abhyankar valuation $v = v_1 \circ \cdots \circ v_r$ is not geometric, then at least one of $v_i$'s has its rational rank $> 1$ by (22) and (31). And for this $v_i$,

$$\Gamma_{v_i} \not\cong \mathbb{Z} \quad \Longrightarrow \quad V_{v_i} \text{ is not Noetherian}$$
$$\Longrightarrow \quad V_{v_i} \text{ is not a localization of a finitely generated } k\text{-algebra.}$$

We also note that all of the proofs [KK05, p.835, Theorem 1.1] [T13, p.114, Theorem 5.5.2] [C21, Theorem 1.3] of Theorem 3.4 are somewhat involved. Fortunately, for this special case of geometric valuations, we may present a more direct and transparent inductive proof with respect to the rank of a geometric valuation $v$ as follows:

*Proof of Theorem 3.6.* Just like the proof of Theorem 3.5, (i) is reduced to (ii). Although the claim about the smooth locus is not stated in Theorem 3.5, this is clear. Thus, it suffices to show (ii), which we only have to prove the output $S\left[\frac{a_1}{b_1}, \ldots, \frac{a_s}{b_s}\right]_{\mathfrak{m}_v \cap S\left[\frac{a_1}{b_1}, \ldots, \frac{a_s}{b_s}\right]}$ enjoys the stated properties of $R$ in (i), in view of Proposition 3.2.

For this, we shall prove by induction on $rank(v)$, and then the case $rank(v) = 1$ is already taken care of by Theorem 3.5. Thus, we shall concentrate on proving the following:

The case $rank < r \implies$ The case $rank = r$ for an inductive local blowup construction of $R = S\left[\frac{a_1}{b_1}, \ldots, \frac{a_s}{b_s}\right]_{\mathfrak{m}_v \cap S\left[\frac{a_1}{b_1}, \ldots, \frac{a_s}{b_s}\right]}$, which enjoys the stated properties in Theorem 3.6 (i).

Now our proof is a pushout of Theorem 3.5 and the proof of [NS14, Theorem 3.1]:

- **Step 1: $\exists$ Local blowup $\widetilde{b} : S \to \widetilde{S}^{(1)}$ w.r.t. $v_1$ s.t. $\widetilde{S}^{(1)} = V_{v_1}$:**
  This procedure is nothing but the the case $rank(v) = rank(v_1) = 1$ of our inductive proof, or the strong form of a local uniformization with respect to divisorial valuations due to Zariski [Z39] (whose proof was given in later literatures like [ZS60b, p.89, VI, §14, Theorem 31] [A86a, p.230, Theorem 5.2] [KM98, p. 61, Lemma 2.45]), which we recalled as Theorem 3.5.
- **Step 2: $\exists$ Local blowup $b : S \to S^{(1)}$ w.r.t. $v_1 \circ \cdots \circ v_r$, inducing $\widetilde{b}$ by the composite $S \xrightarrow{b} S^{(1)} \to S^{(1)}_{\mathfrak{m}_{v_1} \cap S^{(1)}} = \widetilde{S}^{(1)} \overset{\text{Step 1}}{=} V_{v_1}$:**
  Since the local blowup $\widetilde{b} : S \to \widetilde{S}^{(1)}$ w.r.t. $v_1$ is a composition of simple local blowups by Proposition 3.1, we suppose $\widetilde{b} : S \to \widetilde{S}^{(1)}$ is a simple local blowing up as in (39) w.r.t. $v_1$ :

(61)
$$\widetilde{b} : (S, \mathfrak{m}_{v_1} \cap S) \to \left( S\left[\frac{a}{b}\right]_{\mathfrak{m}_{v_1} \cap S\left[\frac{a}{b}\right]}, \mathfrak{m}_{v_1} \cap S\left[\frac{a}{b}\right]_{\mathfrak{m}_{v_1} \cap S\left[\frac{a}{b}\right]} \right), \quad v_1(a) \geq v_1(b),$$

  and we shall search after a local blowup $b : S \to S^{(1)}$ w.r.t. $v_1 \circ \cdots \circ v_r$, inducing $\widetilde{b}$ :
  * **Case 1: $(v_1 \circ \cdots \circ v_r)(a) \geq (v_1 \circ \cdots \circ v_r)(b)$:** In this case, as is shown in [NS14, Lemma 2.15], (61) is obtained by applying the localization w.r.t. $v_1$ to the following composite local blowup w.r.t. $v_1 \circ \cdots \circ v_r$, as desired:
  $$S \to S' := (S[a, b])_{\mathfrak{m}_{v_1 \circ \cdots \circ v_r} \cap S[a,b]} \to S'' := \left( S'\left[\frac{a}{b}\right] \right)_{\mathfrak{m}_{v_1 \circ \cdots \circ v_r} \cap S'\left[\frac{a}{b}\right]}$$
  * **Case 2: $(v_1 \circ \cdots \circ v_r)(a) < (v_1 \circ \cdots \circ v_r)(b)$:** In this case, since $v_1(a) \geq v_1(b)$, we find $v_1(a) = v_1(b)$. This means $\frac{b}{a} = \left(\frac{a}{b}\right)^{-1} \in S\left[\frac{a}{b}\right]_{\mathfrak{m}_{v_1} \cap S\left[\frac{a}{b}\right]}$ and (61) can be identified with
  $$(S, \mathfrak{m}_{v_1} \cap S) \to \left( S\left[\frac{b}{a}\right]_{\mathfrak{m}_{v_1} \cap S\left[\frac{a}{b}\right]}, \mathfrak{m}_{v_1} \cap S\left[\frac{b}{a}\right]_{\mathfrak{m}_{v_1} \cap S\left[\frac{b}{a}\right]} \right), \quad v_1(b) = v_1(a)$$
  Thus, we are in Case 1 with the roles of $a$ and $b$ exchanged. So the claim follows from Case 1.

  Now our desired local blowup $b : S \to S^{(1)}$ is obtained by the iterate of above local blowups as in [NS14, Corollary 2.17].

– **Step 3: $\exists$ Local blowup $\overline{b}^{(1)} : \overline{S}^{(1)} := \frac{S^{(1)}}{\mathfrak{m}_{v_1} \cap S^{(1)}} \to \overline{S}^{(2)}$ w.r.t. $v_2 \circ \cdots \circ v_r$ s.t. $\overline{S}^{(2)}$ enjoys the stated properties in Theorem 3.6 (i):**
This procedure is nothing but the the case $rank(v) = rank(v_2 \circ \cdots \circ v_r) = r-1$ of our inductive proof.

– **Step 4: $\exists$ Local blowup $b^{(1)} : S^{(1)} \to S^{(2)}$ w.r.t. $v_1 \circ v_2 \circ \cdots \circ v_r$ s.t. $S^{(2)}_{\mathfrak{m}_{v_1} \cap S^{(2)}} = S^{(1)}_{\mathfrak{m}_{v_1} \cap S^{(1)}} = V_{v_1}$ and $\frac{S^{(2)}}{\mathfrak{m}_{v_1} \cap S^{(2)}} = \overline{S}^{(2)}$ enjoys the stated properties in Theorem 3.6 (i):**
Since the local blowup $\overline{b}^{(1)} : \overline{S}^{(1)} := \frac{S^{(1)}}{\mathfrak{m}_{v_1} \cap S^{(1)}} \to \overline{S}^{(2)}$ w.r.t. $v_2 \circ \cdots \circ v_r$ is a composition of simple local blowups by Proposition 3.1, we suppose $\overline{b}^{(1)} : \overline{S}^{(1)} \to \overline{S}^{(2)}$ is a simple local blowing up (39) of the form:

(62)

$$\left(\overline{S}^{(1)}, \mathfrak{m}_{v_2 \circ \cdots \circ v_r} \cap \overline{S}^{(1)}\right) \to \left(\overline{S}^{(1)} \left[\frac{\overline{a}}{\overline{b}}\right]_{\mathfrak{m}_{v_2 \circ \cdots \circ v_r} \cap \overline{S}^{(1)}\left[\frac{\overline{a}}{\overline{b}}\right]}, \mathfrak{m}_{v_2 \circ \cdots \circ v_r} \cap \overline{S}^{(1)} \left[\frac{\overline{a}}{\overline{b}}\right]_{\mathfrak{m}_{v_2 \circ \cdots \circ v_r} \cap \overline{S}^{(1)}\left[\frac{\overline{a}}{\overline{b}}\right]}\right),$$

$$\overline{a}, \overline{b} \in \overline{S}^{(1)} := \frac{S^{(1)}}{\mathfrak{m}_{v_1} \cap S^{(1)}}, \text{ s.t. } (v_2 \circ \cdots \circ v_r)(\overline{a}) \geq (v_2 \circ \cdots \circ v_r)(\overline{b})$$

Now, take $a, b \in S^{(1)} \setminus \left(\mathfrak{m}_{v_1} \cap S^{(1)}\right)$ so that their images in $\overline{S}^{(1)} := \frac{S^{(1)}}{\mathfrak{m}_{v_1} \cap S^{(1)}}$ are $\overline{a}, \overline{b}$, respectively. Then, we immediately see $(v_1 \circ v_2 \circ \cdots \circ v_r)(a) \geq (v_1 \circ v_2 \circ \cdots \circ v_r)(b)$ and the following local blowup w.r.t. $v_1 \circ v_2 \circ \cdots \circ v_r$

$$S^{(1)} \to S^{(1)} \left[\frac{a}{b}\right]_{\mathfrak{m}_{v_1 \circ v_2 \circ \cdots \circ v_r} \cap S^{(1)}}$$

reduces to (62), i.e. $\frac{S^{(1)}\left[\frac{a}{b}\right]_{\mathfrak{m}_{v_1 \circ v_2 \circ \cdots \circ v_r} \cap S^{(1)}}}{\mathfrak{m}_{v_1} \cap S^{(1)}\left[\frac{a}{b}\right]_{\mathfrak{m}_{v_1 \circ v_2 \circ \cdots \circ v_r} \cap S^{(1)}}} = \overline{S}^{(1)} \left[\frac{\overline{a}}{\overline{b}}\right]_{\mathfrak{m}_{v_2 \circ \cdots \circ v_r} \cap \overline{S}^{(1)}\left[\frac{\overline{a}}{\overline{b}}\right]}$
while preserving the $v_1$-localization, i.e.
$S^{(1)}_{\mathfrak{m}_{v_1} \cap S^{(1)}} \xrightarrow{\cong} \left(S^{(1)} \left[\frac{a}{b}\right]_{\mathfrak{m}_{v_1 \circ v_2 \circ \cdots \circ v_r} \cap S^{(1)}}\right)_{\mathfrak{m}_{v_1} \cap S^{(1)}\left[\frac{a}{b}\right]_{\mathfrak{m}_{v_1 \circ v_2 \circ \cdots \circ v_r} \cap S^{(1)}}}$
because of our choices of $a, b \in S^{(1)} \setminus \left(\mathfrak{m}_{v_1} \cap S^{(1)}\right)$ [NS14, Lemma 2.18].
Then our desired local blowup $b^{(1)} : S^{(1)} \to S^{(2)}$ is obtained by the iterate of the above local blowups as in [NS14, Corollary 2.20].
In fact, from our construction of $b^{(1)} : S^{(1)} \to S^{(2)}$ in this Step 4 and $b : S \to S^{(1)}$ in Step 2, we recognize the following properties are satisfied:

* $S^{(2)}_{\mathfrak{m}_{v_1} \cap S^{(2)}} = S^{(1)}_{\mathfrak{m}_{v_1} \cap S^{(1)}} = V_{v_1}$, a rank 1 discrete valuation ring.
* By [NS14, Lemma 2.21], there exists some $y \in \mathfrak{m}_{v_1} \cap S^{(2)} \subseteq \mathfrak{m}_{v_1} \cap S^{(2)}_{\mathfrak{m}_{v_1} \cap S^{(2)}}$, which becomes a uniformizing parameter of $S^{(2)}_{\mathfrak{m}_{v_1} \cap S^{(2)}} = S^{(1)}_{\mathfrak{m}_{v_1} \cap S^{(1)}} = V_{v_1}$.
* $\frac{S^{(2)}}{\mathfrak{m}_{v_1} \cap S^{(2)}} = \overline{S}^{(2)}$ such that, for any $2 \leq i \leq r$,

$$\left(\frac{S^{(2)}}{\mathfrak{m}_{v_1 \circ \cdots \circ v_{i-1}} \cap S^{(2)}}\right)_{\left(\mathfrak{m}_{v_1 \circ \cdots \circ v_i} \cap S^{(2)}\right)} = \left(\frac{\overline{S}^{(2)}}{\mathfrak{m}_{v_2 \circ \cdots \circ v_{i-1}} \cap \overline{S}^{(2)}}\right)_{\left(\mathfrak{m}_{v_2 \circ \cdots \circ v_i} \cap \overline{S}^{(2)}\right)} = V_{v_i}.$$

* There exist some $x_2, \dots, x_r \in (\mathfrak{m}_{v_1 \circ \cdots \circ v_r} \setminus \mathfrak{m}_{v_1}) \cap S^{(2)}$, which become a regular system of generators of

$$\frac{\mathfrak{m}_{v_1 \circ \cdots \circ v_r} \cap S^{(2)}}{\mathfrak{m}_{v_1} \cap S^{(2)}} = \mathfrak{m}_{v_2 \circ \cdots \circ v_r} \cap \overline{S}^{(2)} \subset \overline{S}^{(2)} = \frac{S^{(2)}}{\mathfrak{m}_{v_1} \cap S^{(2)}}$$

– **Step 5: $\exists$ Local blowup $b^{(2)} : S^{(2)} \to S^{(3)}$ w.r.t. $v_1 \circ v_2 \circ \cdots \circ v_r$ s.t. $S^{(3)}$ enjoys the requirement of Theorem 3.6:**
In Step 4, we observed there exists some

$$y \in \mathfrak{m}_{v_1} \cap S^{(2)} \subseteq \mathfrak{m}_{v_1} \cap S^{(2)}_{\mathfrak{m}_{v_1} \cap S^{(2)}},$$

which becomes a uniformizing parameter of a rank 1 discrete valuation ring $S^{(2)}_{\mathfrak{m}_{v_1} \cap S^{(2)}} = S^{(1)}_{\mathfrak{m}_{v_1} \cap S^{(1)}} = V_{v_1}$.

* **Case 1: $y \in \mathfrak{m}_{v_1} \cap S^{(2)}$ generates $\mathfrak{m}_{v_1} \cap S^{(2)}$:**
In this case, $\left(S^{(2)}, \mathfrak{m}_{v_1 \circ \cdots \circ v_r} \cap S^{(2)}\right)$ is a regular local ring, because, by [NS14, in the proof of Theorem 3.1], we see

$$\dim S^{(2)} = ht\left(\mathfrak{m}_{v_1 \circ \cdots \circ v_r} \cap S^{(2)}\right)$$
$$\geq ht\left(\mathfrak{m}_{v_1} \cap S^{(2)}\right) + ht\left(\frac{\mathfrak{m}_{v_1 \circ \cdots \circ v_r} \cap S^{(2)}}{\mathfrak{m}_{v_1} \cap S^{(2)}}\right) = 1 + (r-1) \overset{(\bigstar)}{\geq} \dim S^{(2)},$$

where $(\bigstar)$ follows because $y, x_2, \dots, x_r \in \mathfrak{m}_{v_1 \circ \cdots \circ v_r} \cap S^{(2)}$ generate $\mathfrak{m}_{v_1 \circ \cdots \circ v_r} \cap S^{(2)}$ by the assumption. Then, together with Step 4, we find $S^{(2)}$ enjoys the stated properties in Theorem 3.6 (i), and so, we may simply set $S^{(3)} := S^{(2)}$.

* **Case 2: $y \in \mathfrak{m}_{v_1} \cap S^{(2)}$ does not generate $\mathfrak{m}_{v_1} \cap S^{(2)}$:**
In this case, since $S^{(2)}$ is Noetherian, we may still find some $y_1, \dots, y_{r'} \in \mathfrak{m}_{v_1} \cap S^{(2)}$ such that $y, y_1, \dots, y_{r'} \in \mathfrak{m}_{v_1} \cap S^{(2)}$ generate $\mathfrak{m}_{v_1} \cap S^{(2)}$ :

$$\mathfrak{m}_{v_1} \cap S^{(2)} = (y, y_1, \dots, y_{r'}) \subset S^{(2)}. \tag{63}$$

Then, we may apply a single local blowup with respect to $v_1 \circ \cdots \circ v_r$ to reduce to the case $r' = 0$, i.e. to Case 1, as follows:
  · For each $1 \leq k \leq r'$, since

$$y_k \in \mathfrak{m}_{v_1} \cap S^{(2)} \subseteq \mathfrak{m}_{v_1} \cap S^{(2)}_{\mathfrak{m}_{v_1} \cap S^{(2)}} = (y) \subset S^{(2)}_{\mathfrak{m}_{v_1} \cap S^{(2)}},$$

there exists some $a_k \in S^{(2)} \setminus \left(\mathfrak{m}_{v_1} \cap S^{(2)}\right), b_k \in S^{(2)}$ such that

$$y_k = \frac{b_k}{a_k} y \qquad (1 \leq k \leq r'). \tag{64}$$

Setting

$$\begin{aligned} a &:= \prod_{1 \leq k \leq r'} a_k \in S^{(2)} \setminus \left(\mathfrak{m}_{v_1} \cap S^{(2)}\right), \\ y^{(1)} &:= \frac{y}{a} \in \mathfrak{m}_{v_1} \cap S^{(2)}_{\mathfrak{m}_{v_1} \cap S^{(2)}}, \end{aligned} \tag{65}$$

we may rewrite (64) in the following form:

(66)
$$y_k \overset{(64)}{=} \frac{b_k}{a_k} y = \frac{b_k \prod_{1\le j\ne k\le r'} a_j}{a_k \prod_{1\le j\ne k\le r'} a_j} y \overset{(65)}{=} b'_k y^{(1)}$$
$$\text{where } b'_k := b_k \prod_{1\le j\ne k\le r'} a_j \in S^{(2)}, \quad (1\le k\le r').$$

· Now, let us blowup $S^{(2)}$ with respect to $v_1 \circ \cdots \circ v_r$ along the ideal $(a, y)$, where $(v_1 \circ \cdots \circ v_r)(a) < (v_1 \circ \cdots \circ v_r)(y)$ for $v_1(y) > 0 = v_1(a)$. So, we perform the blowup of $S^{(2)}$ given by

(67)
$$\left(S^{(2)}\right)^{(1)} := \left(S^{(2)}\left[y^{(1)}\right]\right)_{\mathfrak{m}_{v_1\circ\cdots\circ v_r}\cap S^{(2)}\left[y^{(1)}\right]} \quad \left(y^{(1)} := \frac{y}{a}\right).$$

· Here, we can show

(68)
$$\mathfrak{m}_{v_1} \cap \left(S^{(2)}\right)^{(1)} = \left(y^{(1)}\right).$$

In fact, $\left(y^{(1)}\right) \subseteq \mathfrak{m}_{v_1} \cap \left(S^{(2)}\right)^{(1)}$ follows immediately from the definitions (65) (67). For $\mathfrak{m}_{v_1} \cap \left(S^{(2)}\right)^{(1)} \subseteq \left(y^{(1)}\right)$, write a general element $z$ of $\mathfrak{m}_{v_1} \cap \left(S^{(2)}\right)^{(1)} = \mathfrak{m}_{v_1} \cap \left(S^{(2)}\left[y^{(1)}\right]\right)_{\mathfrak{m}_{v_1\circ\cdots\circ v_r}\cap S^{(2)}\left[y^{(1)}\right]}$ in the form

$$z = \frac{f(0) + y^{(1)}\overline{f}}{g}, \qquad \text{where}$$
$$f(0) \in S^{(2)}, \overline{f} \in S^{(2)}\left[y^{(1)}\right], g \in S^{(2)}\left[y^{(1)}\right] \setminus \left(\mathfrak{m}_{v_1\circ\cdots\circ v_r} \cap S^{(2)}\left[y^{(1)}\right]\right).$$

Then, since $v_1(r) = v_1(g) = 0, v_1\left(y^{(1)}\overline{f}\right) > 0$, we conclude $v_1(f(0)) = 0$. This means

$$f(0) \overset{(63)}{\in} (y, y_1, \ldots, y_{r'}) \overset{(65)(64)}{\subset} \left(y^{(1)}\right) \subset \left(S^{(2)}\right)^{(1)},$$

which implies $z = \frac{f(0)+y^{(1)}\overline{f}}{g} \in \left(y^{(1)}\right) \subset \left(S^{(2)}\right)^{(1)}$, as desired.

· As in [NS14, Lemma 2.22 (ii)], we can show

(69)
$$\frac{\left(S^{(2)}\right)^{(1)}}{\mathfrak{m}_{v_1} \cap \left(S^{(2)}\right)^{(1)}} = \frac{S^{(2)}}{\mathfrak{m}_{v_1} \cap S^{(2)}} \overset{\text{Step 4}}{=} \overline{S}^{(2)}.$$

For this, write a general element $w$ of $\left(S^{(2)}\right)^{(1)} = \left(S^{(2)}\left[y^{(1)}\right]\right)_{\mathfrak{m}_{v_1\circ\cdots\circ v_r}\cap S^{(2)}\left[y^{(1)}\right]}$ in the form

$$w = \frac{p(0) + y^{(1)}\overline{p}}{q(0) + y^{(1)}\overline{q}} \in \left(S^{(2)}\right)^{(1)} = \left(S^{(2)}\left[y^{(1)}\right]\right)_{\mathfrak{m}_{v_1\circ\cdots\circ v_r}\cap S^{(2)}\left[y^{(1)}\right]},$$

where $p(0), q(0) \in S^{(2)}, \overline{p}, \overline{q} \in S^{(2)}\left[y^{(1)}\right]$, $q(0)+y^{(1)}\overline{q} \in S^{(2)}\left[y^{(1)}\right] \setminus \left(\mathfrak{m}_{v_1\circ\cdots\circ v_r} \cap S^{(2)}\left[y^{(1)}\right]\right)$. Then, since $(v_1 \circ \cdots \circ v_r)\left(q(0) + y^{(1)}\overline{q}\right) = 0, (v_1 \circ \cdots \circ v_r)\left(y^{(1)}\overline{q}\right) > 0$, we see $(v_1 \circ \cdots \circ v_r)(q(0)) = 0$. These in particular imply

$$\begin{aligned} \frac{1}{q(0)} &\in \left(S^{(2)}\right)_{\mathfrak{m}_{v_1 \circ \cdots \circ v_r} \cap S^{(2)}} = S^{(2)}. \\ \frac{1}{q(0)\left(q(0)+y^{(1)}\overline{q}\right)} &\in \left(S^{(2)}\left[y^{(1)}\right]\right)_{\mathfrak{m}_{v_1 \circ \cdots \circ v_r} \cap S^{(2)}\left[y^{(1)}\right]} = \left(S^{(2)}\right)^{(1)}, \end{aligned} \tag{70}$$

and the following (obvious) equality in $\left(S^{(2)}\right)^{(1)}$ :

$$w = \frac{p(0) + y^{(1)}\overline{p}}{q(0) + y^{(1)}\overline{q}} = \frac{p(0)}{q(0)} + \left(\frac{p(0)\overline{q} - q(0)\overline{p}}{q(0)\left(q(0)+y^{(1)}\overline{q}\right)}\right) y^{(1)}.$$

By (65), this means the image of $w = \frac{p(0)+y^{(1)}\overline{p}}{q(0)+y^{(1)}\overline{q}} \in \left(S^{(2)}\right)^{(1)}$ in $\frac{\left(S^{(2)}\right)^{(1)}}{\mathfrak{m}_{v_1} \cap \left(S^{(2)}\right)^{(1)}}$ is the same as the image of $\frac{p(0)}{q(0)} \overset{(70)}{\in} S^{(2)}$. This implies $\frac{\left(S^{(2)}\right)^{(1)}}{\mathfrak{m}_{v_1} \cap \left(S^{(2)}\right)^{(1)}} \cong \frac{S^{(2)}}{\mathfrak{m}_{v_1} \cap S^{(2)}}$, as desired.

· From (68), (69) and Step 4, we see $y^{(1)}, x_2, \ldots, x_r \in \mathfrak{m}_{v_1 \circ \cdots \circ v_r} \cap \left(S^{(2)}\right)^{(1)}$ generate $\mathfrak{m}_{v_1 \circ \cdots \circ v_r} \cap \left(S^{(2)}\right)^{(1)}$. Then, just like Step 5 Case 1, we find $\left(S^{(2)}\right)^{(1)}$ enjoys the stated properties in Theorem 3.6 (i), and so, we may simply set $S^{(3)} := \left(S^{(2)}\right)^{(1)}$.

– **Step 6: Completion of the proof:** For the inductive procedure, take the composition of the local blowups constructed in Step 2, Step 4 and Step 5:

$$S \xrightarrow{b} S^{(1)} \xrightarrow{b^{(1)}} S^{(2)} \xrightarrow{b^{(2)}} S^{(3)},$$

which enjoys the stated properties in Theorem 3.6 (i) by its construction and Proposition 3.2

This concludes the inductive procedure and the proof is complete.

□

## 4. Unramified sheaf, SBNR and Main theorem

We first review our conventions and set up some notations:

- Let $k$ be a perfect field, which gives our base scheme $S = \text{Spec } k$.
- Let $\mathcal{F}_k$ be the category of field extensions $k \subset F$ of $k$ such that $F$ is of finite transcendence degree over $k$. Since an arbitrary element $F$ of $\mathcal{F}_k$ can be written as a filtered colimit of finitely generated $k$-algebra, our perfect assumption of $k$ allows us to conclude $\mathcal{F}_k \subset Pro^{Aff}\left(Sm_k^{ft}\right)$.
- Given $X \in Sm_k^{ft}$, we denote by $X^{(c)}$ $(0 \leq c \leq \dim X)$ the set of codimension $c$ schematic points of $X$, and, for each $x \in X^{(c)}$, $\mathcal{O}_{X,x}$ is a regular local ring of dimension $c$ and Spec $\mathcal{O}_{X,x} \in Pro^{Aff}\left(Sm_k^{ft}\right)$.
- When $c = 0$, each $\alpha \in X^{(0)}$ corresponds to an irreducible component $X_\alpha$ and $\mathcal{O}_{X,\alpha} = k(X_\alpha) \in Pro^{Aff}\left(Sm_k^{ft}\right)$.

- When $c = 1$, for each $x \in X^{(1)}$, the regular local ring $\mathcal{O}_{X,x} \in Pro^{Aff}\left(Sm_k^{ft}\right)$ of dimension 1 is identified with the rank 1 discrete valuation ring $V_{v_x}$ of a uniquely determined rank 1 geometric valuation, a.k.a. divisorial valuation, $v_x$ of $k(X_\alpha)/k$ s.t. $x \in X_\alpha$.
- We denote by $\widetilde{Sm_k^{ft}}$ the subcategory of $Sm_k^{ft}$ consisting of the same objects but only of smooth morphisms. We shall also consider the resulting subcategory $Pro^{Aff}\left(\widetilde{Sm_k^{ft}}\right) \subset Pro^{Aff}\left(Sm_k^{ft}\right)$.

Now, in view of Theorem 1.11, we are urged to search after birational presheaves. For this purpose, let us turn our attention to Morel's unramified presheaves:

**Definition 4.1.** ([M12, Definition 2.1, Remark 2.2]) *An* **unramified presheaf** $S$ *of* $\mathcal{C}$ *(*$\mathcal{C}$ *= sets. groups, or abelian groups) on* $Sm_k^{ft}$ *(resp. on* $\widetilde{Sm_k^{ft}}$*) is a presheaf of* $S$ *of* $\mathcal{C}$ *(*$\mathcal{C}$ *= sets. groups, or abelian groups) on* $Sm_k^{ft}$ *(resp. on* $\widetilde{Sm_k^{ft}}$*), which we uniquely extend to a presheaf on* $Pro^{Aff}\left(Sm_k^{ft}\right)$ *(resp. on* $Pro^{Aff}\left(\widetilde{Sm_k^{ft}}\right)$*)* [17] *by the left Kan extension guaranteed by Corollary 2.2:*

$$S : Pro^{Aff}\left(Sm_k^{ft}\right)^{op} \to \mathcal{C} \quad \left(resp.\ S : Pro^{Aff}\left(\widetilde{Sm_k^{ft}}\right)^{op} \to \mathcal{C}\right)$$

$$\varprojlim_\alpha X_\alpha \mapsto \varinjlim_\alpha S(X_\alpha),$$

*such that the following three conditions hold (Here, for an affine scheme* Spec $A \in Pro^{Aff}\left(Sm_k^{ft}\right)$*, we have abbreviated* $S(\text{Spec } A)$ *simply by* $S(A)$*.):*

**(U0):** *For any* $X \in Sm_k^{ft}$*, the obvious map*

$$S(X) \to \prod_{\alpha \in X^{(0)}} S(X_\alpha)$$

*is a bijection;*

**(U1):** *For any* $X \in Sm_k^{ft}$ *and any everywhere dense open subscheme* $U \subset X$ *the restriction map*

$$S(X) \to S(U)$$

*is injective;*

**(U2):** *With* **(U0)** *at hand, let us suppose* $X \in Sm_k^{ft}$ *is irreducible. Then the injective map (where the injectivity is guaranteed by* **(U1)** *)*

$$S(X) \to \cap_{x \in X^{(1)}} S\left(\mathcal{O}_{X,x}\right)$$

*is a bijection, where* $\cap_{x \in X^{(1)}} S\left(\mathcal{O}_{X,x}\right)$ *is computed in* $S(k(X))$*.*

Morel observed the essence of unramified (pre)sheaves lies upon their values on $\mathcal{F}_k$, which Morel axiomatized as the unramified data [M12, Definition 2.6, Definition 2.9] (see also [C08, 3.2] [F21b, 2.1]):

**Definition 4.2.** *An* ***unramified*** $\widetilde{\mathcal{F}_k}$***-datum of sets (resp. groups)*** *consists of:*

[17] This extension of $S$ allows us to define $S(k(X))$ and $S(\mathcal{O}_{X,x})$ $(x \in X)$ for $X \in Sm_k$, $S(V_v)$ for a divisorial valuation $v$ by Theorem 3.5, and, more generally, $S(R)$ for an arbitrary regular local ring $R$ over a perfect field $k$ by Popescu's theorem (Theorem 2.3, Corollary 2.5), amongst of all.

**D1**: *A continous functor*

$$S : \mathcal{F}_k \to Set \ (resp.\ Grp,\ Ab),$$

*i.e.* $S(F = \varinjlim_\alpha F_\alpha) = \varinjlim_\alpha S(F_\alpha)$ *where* $F_\alpha$ *run over the set of subfields of* $F$ *of finite type over* $k$;

**D2**: *For any* $F \in \mathcal{F}_k$ *and any divisorial valuation* $v$ *on* $F$, [18] *a subset (resp. subgroup)*

$$S(V_v) \ \subset \ S(F)$$

*These data should satisfy the following axioms:*

**A1**:

- *If* $i : E \subseteq F$ *is a separable extension in* $\mathcal{F}_k$, *and* $v$ *is a divisorial valuation on* $F$ *which restricts to a divisorial valuation* $w$ *on* $E$ *with ramification index* $e(v/w) = 1$ *then* $S(i)$ *maps* $S(V_w)$ *into* $S(V_v)$;

$$\begin{array}{ccc} S(V_w) & \overset{\exists 1}{\dashrightarrow} & S(V_v) \\ \boxed{D2}\ \cap\downarrow & & \cap\downarrow\ \boxed{D2} \\ S(E) & \xrightarrow[\boxed{D1}]{S(i)} & S(F) \end{array}$$

- *moreover if the induced extension* $i : \kappa(w) \to \kappa(v)$ *is an isomorphism, then the above square of sets (resp. groups) is cartesian (which is induced from a Nisnevich elementary distinguished square* [19] *over* $\mathrm{Spec}\,(V_w)$ : *):*

$$\begin{array}{ccc} \mathrm{Spec}(F) & \longrightarrow & \mathrm{Spec}\,(V_v) \\ \downarrow & & \downarrow \\ \mathrm{Spec}(E) & \longrightarrow & \mathrm{Spec}\,(V_w) \end{array}$$

**A2**: *Let* $X \in Sm_k$ *be irreducible with function field* $F$. *Then,*

$$S(F) = \cup_{\substack{\Phi \subset X^{(1)} \\ \Phi,\ finite\ subset}} \cap_{x \in X^{(1)} \setminus \Phi} S(\mathcal{O}_{X,x})$$

*In other words, any* $\sigma \in S(F)$ *is contained in all but a finite number of* $S(\mathcal{O}_{X,x})$*'s, where* $x$ *runs over the set* $X^{(1)}$ *of points of codimension one.*

---

[18] Since $v$ is a rank 1 discrete valuation, its valuation ring $V_v$ is well known to be regular (see e.g. [AM69, Prop.9.2;Th.11.22]). However, we had already observed directly $V_v \in Pro^{Aff}\left(Sm_k^{ft}\right)$ in Theorem 3.5, even without using Popescu's theorem.

[19] (See e.g. [MV99, p.96]) An (Nisnevich) elementary distinguished square is a cartesian square of the form

$$\begin{array}{ccc} U \times_X V & \longrightarrow & V \\ \downarrow & & \downarrow p \\ U & \xrightarrow{j} & X \end{array}$$

s.t. $p$ is étale , $j$ an open embedding with $p^{-1}((X \setminus U)_{red}) \to (X \setminus U)_{red}$ an isomorphism.

This relevance to the Nisnevich distinguished square was featured in [C19, Proof of Proposition 1.2], noting that Spec $(V_v) \to$ Spec $(V_w)$ is étale in this particular case (see e.g. [sp18, 09E4,09E7,07NP,00U6] for this).

*An* <u>*unramified $\mathcal{F}_k$-datum*</u> *$S$ is an unramified $\widetilde{\mathcal{F}_k}$-data together with the following additional data:*

$\boxed{\textbf{D3}}$ **:** *For any $F \in \mathcal{F}_k$ and any divisorial valuation $v$ on $F$, a map $s_v : S(V_v) \to S(\kappa(v))$, called the* <u>*specialization map associated to $v$*</u>.

*Furthermore, these data should satisfy the following additional axioms:*

$\boxed{\textbf{A3}}$ **:** *If $i : E \subset F$ is an extension in $\mathcal{F}_k$, and $v$ is a divisorial valuation on $F$.*

(i) *If $v$ restricts to a divisorial valuation $w$ on $E$, then $S(i)$ maps $S(V_w)$ to $S(V_v)$ and the following diagram is commutative:*

$$\begin{array}{ccccc} S(E) & \underset{\boxed{D2}}{\hookleftarrow} & S(V_w) & \xrightarrow[\boxed{D3}]{s_w} & S(\kappa(w)) \\ {\scriptstyle S(i)}\downarrow \boxed{D1} & & \exists 1 \, \vdots & & \boxed{D1} \downarrow {\scriptstyle S(\bar{i})} \\ S(F) & \underset{\boxed{D2}}{\hookleftarrow} & S(V_v) & \xrightarrow[\boxed{D3}]{s_v} & S(\kappa(v)) \end{array}$$

(ii) *If $v$ restricts to $0$ on $E$ and denote by $j : E \subset \kappa(v)$ the resulting field extension, then $S(i)$ maps $S(E)$ to $S(V_v)$ and the following diagram is commutative:*

$$\begin{array}{ccccc} S(E) & & \xrightarrow{\quad S(j)\ \boxed{D1} \quad} & & S(\kappa(v)) \\ \boxed{D1}\ \downarrow {\scriptstyle S(i)} & \overset{\exists 1}{\searrow} & & & \\ S(F) & \underset{\boxed{D2}}{\hookleftarrow} & S(V_v) & \xrightarrow[\boxed{D3}]{s_v} & S(\kappa(v)) \end{array}$$

$\boxed{\textbf{A4}}$ **:**

(i) *Let $X = \mathrm{Spec}A \in Pro^{Aff}\left(Sm_k^{ft}\right)$ be an arbitrary affine scheme of a dimension $2$ local ring $A$ with its unique closed point $z$. Then, for any $y_0 \in X^{(1)}$, which corresponds to a uniquely determined divisorial valuation $v_{y_0}$, such that $\overline{y_0} \in Pro^{Aff}\left(Sm_k^{ft}\right)$, the specialization map $s_{v_{y_0}} : S\left(\mathcal{O}_{X,y_0}\right) = S\left(V_{v_{y_0}}\right) \to S\left(\kappa(v_{y_0})\right) = S\left(\kappa(y_0)\right)$ maps $\cap_{y \in X^{(1)}} S\left(\mathcal{O}_{X,y}\right)$ into $S\left(\mathcal{O}_{\overline{y_0},z}\right)$.*

$$\begin{array}{ccc} S\left(\mathcal{O}_{X,y_0}\right) & \xrightarrow[\boxed{D3}]{s_{y_0}} & S\left(\kappa(y_0)\right) = S\left(k\left(\overline{y_0}\right)\right) \\ \uparrow & & \uparrow \boxed{D2} \\ \cap_{y \in X^{(1)}} S\left(\mathcal{O}_{X,y}\right) & \overset{\exists 1}{\dashrightarrow} & S\left(\mathcal{O}_{\overline{y_0},z}\right) \end{array}$$

(ii) *The composition*

$$\cap_{y \in X^{(1)}} S\left(\mathcal{O}_{X,y}\right) \overset{(i)}{\longrightarrow} S\left(\mathcal{O}_{\overline{y_0},z}\right) \xrightarrow[\boxed{(D3)}]{s_{v_z}} S\left(\kappa(z)\right),$$

*where $v_z$ is the uniquely deternined divisorial valuation of $k(\overline{y_0})/k$ whose valuation ring is $\mathcal{O}_{\overline{y_0},z}$, doesn't depend on the choice of $y_0$ such that $\overline{y_0} \in Pro^{Aff}\left(Sm_k^{ft}\right)$.*

In fact, Morel [M12] showed unramified sheaves and unramified $\mathcal{F}_k$-data are equivalent as follows:

**Theorem 4.3.** (i) ([M12, Proposition 2.8]) *The category of unramified* $\boxed{\widetilde{\mathcal{F}_k}}$*-data of sets (resp. groups, abelian groups) is equivalent to the category of unramified sheaves of sets (resp. groups, abelian groups) for the Nisnevich topology on* $\boxed{\widetilde{Sm_k^{ft}}.}$

(ii) ([M12, Theorem 2.11]) *Given an unramified* $\boxed{\mathcal{F}_k}$*-datum* $S$, *there is a unique way to extend the above unramified sheaf defined on* $\boxed{\widetilde{Sm_k^{ft}}}$ *to a sheaf on* $\boxed{Sm_k^{ft},}$ *such that*

- *for any divisorial valuation* $v$ *on* $F/k$ *such that* $F \in \mathcal{F}_k$, *inducing the map* $\pi(v) : V_v \twoheadrightarrow \kappa(v) = V_v/\mathfrak{m}_v$ *in* $Pro^{Aff}\left(Sm_k^{ft}\right)$, [20] *the map* $S(\pi(v)) : S(V_v) \to S(\kappa(v))$ *induced by the sheaf structure is the specialization map* $s_v$ *of* $S$,
- *and the extended sheaf is unramified.*

For our purpose, it is useful to review:

*An outline of Morel's proof of: unramified datum $\Longrightarrow$ unramified sheaf:* .

(1) Given an unramified datum $S$, we first extend $S$ to a presheaf on $Sm_k^{ft}$ at object level, by setting $S(X)$ for $X \in Sm_k^{ft}$ as **(U0)** and **(U2)** in Definition 4.1:

$$S(X) := \prod_{\alpha \in X^{(0)}} S(X_\alpha) := \prod_{\alpha \in X^{(0)}} \left( \cap_{x \in X_\alpha^{(1)}} S(\mathcal{O}_{X_\alpha, x}) \right) \left( \subset \prod_{\alpha \in X^{(0)}} k(X_\alpha) \right)$$

(2) To extend $S$ to a presheaf on $Sm_k^{ft}$ at morphism level, observe that an arbitrary morphism $f : X \to Y$ in $Sm_k^{ft}$ decomposes as

$$X \overset{\Gamma_f}{\hookrightarrow} X \times Y \overset{\pi_Y}{\twoheadrightarrow} Y; \qquad x \overset{\Gamma_f}{\mapsto} (x, y = f(x)) \overset{\pi_Y}{\mapsto} y = f(x), \tag{71}$$

where $\Gamma_f$ is a regular immersion and $\pi_Y$ is a smooth morphism (see e.g. [F98, B7] [sp18, 0E9K]).

(3) To define $S$ for a regular immersion $i : Z \hookrightarrow X$, observe that we can cover $X$ by Zariski opens $U$ s.t. each induced regular immersion $i_{Z\cap U} : Z \cap U \hookrightarrow U$ is a composition

$$i_{Z\cap U} : Z \cap U = Y_0 \overset{j_0}{\hookrightarrow} Y_1 \overset{j_1}{\hookrightarrow} Y_2 \overset{j_2}{\hookrightarrow} \cdots \overset{j_{d-1}}{\hookrightarrow} Y_d = U \tag{72}$$

of codimension one regular immersions $j_k : Y_k \hookrightarrow Y_{k+1}$ $(0 \le k \le d-1)$. (For an existence of a decomposition as in (72), see [G67, Corollary 17.12.2] and [C19, Lemma 1.5].)

(4) Any smooth morphism $f : Y \to X$ induces the corresponding extension of function fields

$$k(f) : \; E := k(X) \; \subset \; k(Y) =: F.$$

[20] In the original exposition of [M12, Lemma 2.12] (and [C19, Theorem 1.3]), $\kappa(v)$ is required to be separable. However, this is automatically satisfied because our base field is assumed to be perfect [sp18, 030Y,030Z].

Now, $k(f)$ induces $S(k(f))$, which Morel observed to restrict to $S(f)$ because of **(U0)**, **(U2)** and $\boxed{\textbf{A1}}$ [M12, Proof of Proposition 2.8] (see also [C19, Proof of Proposition 1.2] ) :

$$\begin{array}{ccc} S\left(k(X)\right) & \hookleftarrow & S(X) = \bigcap_{x \in X^{(1)}} S\left(\mathcal{O}_{X,x}\right) \\ \Big\downarrow{\scriptstyle S(k(f))} & & \vdots\ {\scriptstyle \exists S(f)} \\ S\left(k(Y)\right) & \hookleftarrow & S(Y) = \bigcap_{y \in Y^{(1)}} S\left(\mathcal{O}_{Y,y}\right) \end{array} \tag{73}$$

(5) Any codimension 1 regular immersion $i : \overline{\{y\}} = Y \hookrightarrow X$ $(y \in X^{(1)})$ induces the canonical quotient

$$\pi_y : \mathcal{O}_{X,y} \twoheadrightarrow \kappa(y) := \mathcal{O}_{X,y}/\mathfrak{m}_{X,y}, \quad \text{a morphism in } Pro^{Aff}\left(Sm_k^{ft}\right).$$

Now $\pi_y$ induces the specialization map $s_{v_y}$, where $v_y$ is the rank 1 geometric valuation corresponding to the discrete valuation ring $\mathcal{O}_{X,y}$, and Morel observed $s_{v_y}$ restrict to $S(i)$ because of **(U0)**, **(U2)** and $\boxed{\textbf{A4}}$ [M12, Lemma 2.12] (see also [C19, Theorem 1.3] ) :

$$\begin{array}{ccc} S\left(\mathcal{O}_{X,y}\right) & \hookleftarrow & S(X) = \bigcap_{y' \in X^{(1)}} S\left(\mathcal{O}_{X,y'}\right) \\ \Big\downarrow{\scriptstyle s_{v_y}} & & \vdots\ {\scriptstyle \exists S(i)} \\ S(\kappa(y)) = S\left(k(Y)\right) & \hookleftarrow & S(Y) = \bigcap_{z \in Y^{(1)}} S\left(\mathcal{O}_{Y,z}\right) \end{array} \tag{74}$$

(6) For a general regular immersion $i : Z \hookrightarrow X$, set $S(i)$ by patching together $S\left(i_{Z\cap U}\right)$'s, which we set to be $S\left(j_0\right) \circ S\left(j_1\right) \circ S\left(j_2\right) \circ \cdots \circ S\left(j_{d-1}\right)$ as the composition of $S$'s of the codimsion one regular immersions in (72). Of course, we must check its well-definedness, but this is what Morel did in [M12, Lemma 2.13] (see also [C19, Lemma 1.4]).

(7) For a general morphism $f : X \to Y$ in $Sm_k$, set $S(f)$ to be the composition $S\left(\Gamma_f\right) \circ S\left(\pi_Y\right)$ of $S$'s of the regular immersion $\Gamma_f$ and the smooth morphism $\pi_Y$ in (71).

We are then left to show the this correspondence $S$ of smooth morphisms is functorial. However, this is what Morel did in [M12, just before Remark 2.14] (see also [C19, p.10–11]).

□

Actually, there are considerably many examples of unramified sheaves:

**Example 4.4.** (i) ([M12, Remark 6.10] [K21a, Corollary 2.8]) Any <u>strongly $\mathbb{A}^1$ invariant Nisnevich sheaf of group</u> $G$, i.e. $U \mapsto H^n_{Nis}(U, G)$ is $\mathbb{A}^1$ invariant for $n = 0, 1$, is unramified.

(ii) Also, familiar $\mathbb{A}^1$ invariant Nisnevich sheaves with transfer are unramified. In fact, these are special cases of (i), because, by [MVW06, Theorem 13.8], any $\mathbb{A}^1$ invariant Nisnevich sheaf of abelian groups with transfer $F$ is <u>*strictly* $\mathbb{A}^1$ *invariant*</u>, i.e. $U \mapsto H^n_{Nis}(U, F)$ is $\mathbb{A}^1$ invariant for any $n \in \mathbb{Z}_{\geq 0}$ under our perfect base field assumption.

Déglise [D06, Proposition 6.9] [D11, Théorème 3.7] showed an $\mathbb{A}^1$ invariant Nisnevich sheaves of abelian groups with transfer $F$ is essentially equivalent to Rost's cycle module [R96]. Rost's cycle module is a data of the following form:

$$M = \Big( M_* : \mathcal{F}_k \to \mathcal{A}b_*, \{\phi^* : M_*(F) \to M_*(F) \mid \text{finite extension } E \subset F \text{ in } \mathcal{F}_k\}, \{\partial_v : M_*(F) \to M_*(\kappa(v)) \mid v, \text{ a geometric discrete valuation on } F|k \in \mathcal{F}_k\} \Big),$$

and the desired $\mathbb{A}^1$ invariant Nisnevich sheaf of abelian groups with transfer is given by the associated Rost Chow groups with coefficients in the given cycle module $M$:
(75)

$$U \mapsto A^0(U, M_*) := \mathrm{Ker}\left( \oplus_{x\in U^{(0)}} M_*(\kappa(x)) \xrightarrow{\oplus_{x\in U^{(0)}}\left(\oplus_{y\in U^{(1)}}\partial^x_y\right)} \oplus_{y\in U^{(1)}} M_{*-1}(\kappa(y)) \right),$$

where $U^{(c)}$ is the set of codimension $c$ schematic points of $U$. [21] While we only consider the degree 0 part in (75), the higher degree parts also emerges as the Bloch-Ogus [BO74] type theorem for Rost's cycle modules [R96, Theorem (6.1), Corollary (6.5)], which gives us a conceptual transparent description of $A^p(X, M_*)$ for a smooth $k$-variety $X$:

$$A^p(X, M_*) \cong H^p_{Zar}(X, \mathcal{M}_*), \tag{76}$$

where $\mathcal{M}_*$ is the Zariski sheaf given by (75).

For instance, starting with the Galois cohomology as Rost's cycle module

$$\left( H^j_{Gal}(-, D \otimes \mu_r^{\otimes j}), \{\phi^*\}, \{\partial_v\} \right),$$

where $r \in \mathbb{N}$ is taken to be coprime to *char* $k$ with $\mu_r \subset \overline{k}^*$ the group of $r$-th roots of unity and $D$ a finite continuous *Gal* $(\overline{k}/k)$-module of exponent $r$ [R96, p.335, Remark (1.11)], we arrive at the familiar *unramified cohomology*

$$\begin{aligned} H^j_{nr}\left(X, D \otimes \mu_r^{\otimes j}\right) &= A^0\left(X, (H^j_{Gal}(-, D \otimes \mu_r^{\otimes j}), \{\phi^*\}, \{\partial_v\})\right) \\ &\cong H^0_{Zar}\left(X, \mathcal{H}^j_{ét}(D \otimes \mu_r^{\otimes j})\right). \end{aligned}$$

as un unramified Zariski sheaf. Again, see [R96] [M05a, 2.2] [AB17, 3] for more details including the notation.

(iii) ([F20][F21a][F21b]) The concept of Rost's cycle module was generalized by Feld [F20] to his *Milnor-Witt cycle modules.* Whereas Rost's cycle module admits a graded action of the Milnor $K$-theory $\mathbf{K}^M_*$ [M70], Feld's Milnor-Witt cycle module admits a graded action of the Morel's Milnor-Witt $K$-theory $\mathbf{K}^{MW}_*$ [M12, Definition 3.1], so that the Rost cycle module is a special case of the Milnor-Witt cycle module with the trivial action of $\eta \in \mathbf{K}^{MW}_{-1}$. (Recall $\mathbf{K}^{MW}_*/(\eta) \cong \mathbf{K}^M_*$ [M12, Remark 3.2].)

For any Milnor-Witt cycle module $M$ and any smooth $k$-scheme $X$, just like the case of Rost's cycle module, a $\mathbb{Z}$-graded pair $\left\{F^M_n, \omega_n : F^M_{n-1} \to (F^M_n)_{-1}\right\}_{n\in\mathbb{Z}}$ are produced so that the $\mathbb{Z}$-graded Zariski sheaves of abelian groups $\left\{F_n^M\right\}_{n\in\mathbb{Z}}$ are of the form

$$\begin{cases} U \mapsto F_n^M(U) := A^0\left(U, M, -\Omega_{U/k} + \langle n\rangle\right) \cong H^0_{Zar}\left(U, \mathcal{M}_{U,n}\right) \\ A^p\left(X, M, -\Omega_{X/k} + \langle n\rangle\right) \cong H^p_{Zar}\left(X, \mathcal{M}_{X,n}\right) \end{cases}, \tag{77}$$

making $\left\{F_n^M, \omega_n : F_{n-1}^M \to (F_n^M)_{-1}\right\}_{n\in\mathbb{Z}}$ into a *homotopy module* [F20, Corollary 8.5] [F21a, Theorem 4.1.7], i.e. a pair $\{M_*, \omega_*\}_{n\in\mathbb{Z}}$ of a $\mathbb{Z}$-graded strictly $\mathbb{A}^1$-invariant Nisnevich sheaf of abelian groups $M_*$ on $Sm_k$ and the *desuspension map* $\omega_n : M_{n-1} \to (M_n)_{-1}$ for each $n \in \mathbb{Z}$. See [F20][F21a][F21b] for more details including the notation.

[21]For the precise meanings of these notations in this expression and general facts about Rost's cycle modules, consult surveys [M05a, 2.2] [AB17, 3] or Rost's original article [R96]. For its relation with the Bloch-Ogus theory [BO74], see [CTHK97, 7.3(5)].

Consequently, each $F_n^M$ is *unramified* by (i).

(iv) By Morel [M12, Theorem 6.1, Corollary 6.2], for any pointed space $\mathcal{X} \in Fun\left(\Delta^\bullet, Pr(Sm_k^{ft})\right)$, its homotopy group sheaf $\pi_n^{\mathbb{A}^1}(\mathcal{X})$, which is defined to be the associated Nisnevich sheaf to the preseaf

$$\begin{aligned} Sm_k^{ft} &\to Groups \\ U &\mapsto Hom_{\mathcal{H}_\bullet(k)}\left(\Sigma^n(U_+), \mathcal{X}\right), \end{aligned}$$

is unramified for $n \in \mathbb{Z}_{\geq 1}$.[22] Here, because of the adjunction:

$$Hom_{\mathcal{H}_\bullet(k)}\left(\Sigma^n(U_+), \mathcal{X}\right) \cong Hom_{\mathcal{H}_\bullet(k)}\left(\Sigma(U_+), \Omega_{S^1}^{n-1}\mathcal{X}\right),$$

where $\Omega_{S^1}^{n-1}\mathcal{X}$ is the *(n-fold) simplicial loop space* (c.f. e.g. [VRO07, p.215, Exercise 5.51]), the essence of Morel's proof was the case $n = 1$. In fact, in the course of proving the unramified property for this critical case of $n = 1$, Morel [M12, p,154–155, Proof of Corollary 6.9 2)] actually argued along the following line:

- For a fixed $\mathbb{A}^1$-connected and $\mathbb{A}^1$-local space $\mathcal{B}$, denote by $\pi_1(\mathcal{B})$ the presheaf

  $$\begin{aligned} \pi_1(\mathcal{B}) : Sm_k^{ft} &\to Groups \\ X &\mapsto \pi_1(\mathcal{B})(X) := Hom_{\mathcal{H}_\bullet(k)}\left(\Sigma(X_+), \mathcal{B}\right) \end{aligned}$$

  and set $\pi_1^{Zar}(\mathcal{B})$ and $\pi_1^{\mathbb{A}^1}(\mathcal{B})$ respectively to be the associated sheafications with respect to the Zariski topology and Nisnevich topology, respectively. We extend these (pre)sheaves to $Pro^{Aff}\left(Sm_k^{ft}\right)$ by the left Kan extension as usual. Of course, $\pi_1^{\mathbb{A}^1}(\mathcal{B})$ is the homotopy group sheaf in question (note however that our $\pi_1^{Zar}(\mathcal{B})$ is simply denoted by $\mathcal{G}$ in [M12, p,154–155, Proof of Corollary 6.9 2)]).
- For the affine scheme $L \in Pro^{Aff}\left(Sm_k^{ft}\right)$ of a Zariski local ring of $X \in Sm_k^{ft}$ (which also includes the function field $k(X)$ as the local ring at the generic point), we have the following short exact sequence of groups, as a consequence of Gabber's geometric representation lemma [G94][CTHK97] (c.f. [M12, Lemma 1.15]):

$$1 \to \pi_1(\mathcal{B})(L) \to \pi_1(\mathcal{B})(k(X)) \Rightarrow {\prod}'_{y\in L^{(1)}} \pi_1(\mathcal{B})(k(X))/\pi_1(\mathcal{B})\left(\mathcal{O}_{L,y}\right), \tag{78}$$

  where the double arrow refers the action of $\pi_1(\mathcal{B})(k(X))$ on ${\prod}'_{y\in L^{(1)}} \pi_1(\mathcal{B})(k(X))/\pi_1(\mathcal{B})\left(\mathcal{O}_{L,y}\right)$ with the isotropy subgroup $\pi_1(\mathcal{B})(L)$ included by the middle single arrow.

[22] See also [C08, Corollary 5.2.12].

- Under the same assumption and notation, $\pi_1(\mathcal{B})(L) = \pi_1^{Zar}(\mathcal{B})(L)$, which together with the group short exact sequence (78), we obtain the following group short exact sequence for any irreducible $X \in Sm_k^{ft}$ :

(79)
$$1 \to \pi_1^{Zar}(\mathcal{B})(X) \to \pi_1^{Zar}(\mathcal{B})(k(X)) \Rightarrow {\prod}'_{y \in X^{(1)}} \pi_1^{Zar}(\mathcal{B})(k(X))/\pi_1^{Zar}(\mathcal{B})\left(\mathcal{O}_{X,y}\right),$$

- Considering the definition of the Nisnevich topology and the group short exact sequence (79), we find $\pi_1^{Zar}(\mathcal{B}) = \pi_1^{\mathbb{A}^1}(\mathcal{B})$, the homotopy group sheaf in question. Furthermore, as is show in [M12, p.23; p.26, Remark 2.23],
$$\begin{aligned} \pi_1^{Zar}(\mathcal{B})(k(X))/\pi_1^{Zar}(\mathcal{B})\left(\mathcal{O}_{X,y}\right) &\cong \pi_1^{\mathbb{A}^1}(\mathcal{B})_{-1}\left(\kappa(v_y)\right) \\ &= \pi_1^{Zar}(\mathcal{B})_{-1}\left(\kappa(v_y)\right) = \pi_1(\mathcal{B})_{-1}\left(\kappa(v_y)\right), \end{aligned}$$
where $v_y$ is the divisorial valuation associated with $y \in X^{(1)}$, and, for a presheaf $G$ of groups on $Sm_k^{ft}$, $G_{-1}$ is the **Voevodsky contraction**, defined by
$$\begin{aligned} G_{-1} : Sm_k^{ft} &\to Groups \\ X &\mapsto \operatorname{Ker}\left(G(\mathbb{G}_m \times X) \xrightarrow{ev_1} G(X)\right). \end{aligned}$$

To sum up, the homotopy group sheaf $\pi_1^{\mathbb{A}^1}(\mathcal{B})$ is unramified, and, consequently, we may consider its associated stably birationalized subsheaf $\pi_1^{\mathbb{A}^1}(\mathcal{B})_{sb}$. Furthermore, from the above summary, we have the expressions of
$$\pi_1^{\mathbb{A}^1}(\mathcal{B})_{sb}\left(k(X)\right) \subseteq \pi_1^{\mathbb{A}^1}(\mathcal{B})(X)$$
for arbitrary $X \in Sm_k^{ft}$ as follows:

**homotopy group sheaf is unramified homotopy group**

For $X \in Sm_k^{ft}$, we have the commutative diagram with horizontal group exact sequences, characterizing the embedding: $\pi_1^{\mathbb{A}^1}(\mathcal{B})_{sb}\left(k(X)\right) \subseteq \pi_1^{\mathbb{A}^1}(\mathcal{B})(X)$, together with the canonical homomorphism: $\pi_1(\mathcal{B})(X) \to \pi_1^{\mathbb{A}^1}(\mathcal{B})(X)$ :

$$\begin{array}{ccccccc} 1 & \longrightarrow & \pi_1^{\mathbb{A}^1}(\mathcal{B})_{sb}\left(k(X)\right) & \longrightarrow & \pi_1(\mathcal{B})\left(k(X)\right) & \Longrightarrow & {\prod}'_{\substack{v\,:\,\text{divisorial valuation} \\ \text{on } k(X)/k}} \pi_1(\mathcal{B})_{-1}\left(\kappa(v)\right) \\ & & \big\downarrow & & \big\| & & \big\downarrow \\ 1 & \longrightarrow & \pi_1^{\mathbb{A}^1}(\mathcal{B})\left(X\right) & \longrightarrow & \pi_1(\mathcal{B})\left(k(X)\right) & \Longrightarrow & {\prod}'_{y \in X^{(1)}} \pi_1(\mathcal{B})_{-1}\left(\kappa(v_y)\right) \\ & & & \nwarrow^{\exists 1} & \big\uparrow & & \\ & & & & \pi_1(\mathcal{B})(X) & & \end{array}$$

Of course, our Theorem and Definition 1.12 implies $\pi_1^{\mathbb{A}^1}(\mathcal{B})_{sb}$ is a stable birational invariant Nisnevich sheaf on $Sm_k^{ft}$ and $\pi_1^{\mathbb{A}^1}(\mathcal{B})$, which is equal to $\pi_1^{\mathbb{A}^1}(\mathcal{B})_{sb}$ for smooth proper $k$-schemes, is a stable birational invariant of smooth proper $k$-schemes. Thus, it might be more natural to call $\pi_1^{\mathbb{A}^1}(\mathcal{B})$ an unramified homotopy sheaf and denote by $\pi_1(\mathcal{B})_{nr}$, instead.

Of course, this result is also applicable to generalized motivic cohomology theories. Actually, there are two kinds of motivic stable homotopy categories, where generalized motivic cohomology theories can be defined. The first one is more "unstable" <u>*stable* $\mathbb{A}^1$*-homotopy category of* $S^1$*-spectra*</u> $\mathcal{SH}^{S^1}(k)$ [M04, p.405] [M05b, p.27, Definition 4.1.1.(3)] [VRO07, Definition 2.6], where a generalized motivic cohomology theory $E^s(-)_t$ is represented by a $S^1$-spectrum $E$ as follows:

$$E^s(-)_t : Sm_k^{ft} \to AbelianGroups$$
$$X \mapsto \begin{cases} \mathrm{Hom}_{\mathcal{SH}^{S^1}(k)}\left(\Sigma^\infty X_+, E \wedge \left(\Sigma^\infty S^1\right)^s \wedge \left(\Sigma^\infty \mathbb{G}_m\right)^{\wedge t}\right) & \text{if } t \geq 0 \\ \mathrm{Hom}_{\mathcal{SH}^{S^1}(k)}\left(\Sigma^\infty X_+ \wedge \left(\Sigma^\infty \mathbb{G}_m\right)^{\wedge -t}, E \wedge \left(\Sigma^\infty S^1\right)^s\right) & \text{if } t \leq 0 \end{cases}$$

Then, we shall call its Nisnevich sheafication, which we denote by $E^*_{nr}(-)_t$, the <u>*unramified* $E^*(-)_t$ *thery*</u>, and call its stably birationalized Nisnevich subsheaf, which we denote by $E^*_{sb}(-)_t$, the <u>*stably birationalized* $E^*(-)_t$ *theory*.</u> Then, our analysis of of the homotopy group sheaf, presented above, implies the following:

$$E^s_{sb}\left(k(X)\right)_t \hookrightarrow E^s_{nr}(X)_t \overset{\exists 1}{\longleftarrow} E^s(X)_t$$

For $X \in Sm_k^{ft}$, we have the commutative diagram with horizontal group exact sequences, characterizing the embedding: $E^s_{sb}\left(k(X)\right)_t \subseteq E^s_{nr}(X)_t$, together with the canonical homomorphism: $E^s(X)_t \to E^s_{nr}(X)_t$

$$\begin{array}{ccccccc} 0 & \longrightarrow & E^s_{sb}\left(k(X)\right)_t & \longrightarrow & E^s\left(k(X)\right)_t & \longrightarrow & \oplus_{v\,:\,\text{divisorial valuation on } k(X)/k} E^s\left(\kappa(v)\right)_{t-1} \\ & & \big\downarrow & & \big\| & & \big\downarrow \\ 0 & \longrightarrow & E^s_{nr}(X)_t & \longrightarrow & E^s\left(k(X)\right)_t & \longrightarrow & \oplus_{y \in X^{(1)}} E^s\left(\kappa(v_y)\right)_{t-1} \\ & & & \nwarrow^{\exists 1} & \big\uparrow & & \\ & & & & E^s(X)_t & & \end{array}$$

The second motivic stable homotopy category is the more standard <u>*stable* $\mathbb{A}^1$*-homotopy category of* $\mathbb{P}^1$*-spectra*</u> $\mathcal{SH}^{\mathbb{P}^1}(k)$ [V98, Definition 5.7] [M04, p.417] [VRO07, Definition 2.10], which is also obtained from $\mathcal{SH}^{S^1}(k)$ by "stabilizing" with respect to the Tate sphere $\mathbb{G}_m$ [VRO07, Proposition 2.3]. Since this is more standard, we shall denote $\mathcal{SH}^{\mathbb{P}^1}(k)$ simply by $\mathcal{SH}(k)$, and call it simply the <u>*motivic stable homotopy category*.</u> Then its generalized motivic cohomology theory $\mathbb{M}^{p \cdot q}(-)$ is represented by a $\mathbb{P}^1$-spectrum $M$ as follows:

$$\mathbb{M}^{p \cdot q}(-) : Sm_k^{ft} \to AbelianGroups$$
$$X \mapsto \mathrm{Hom}_{\mathcal{SH}(k)}\left(\Sigma^\infty_{\mathbb{P}^1} X_+, M \wedge \left(\Sigma^\infty_{\mathbb{P}^1} S^1\right)^{p-q} \wedge \left(\Sigma^\infty_{\mathbb{P}^1}\mathbb{G}_m\right)^{\wedge q}\right)$$

Then, we shall call its Nisnevich sheafication, which we denote by $\mathbb{M}^{*,*}_{nr}$, the <u>*unramified* $\mathbb{M}^{*,*}$ *thery*</u>, and call its stably birationalized Nisnevich subsheaf, which we denote by $\mathbb{M}^{*,*}_{sb}$, the <u>*stably birationalized* $\mathbb{M}^{*,*}$ *thery*.</u> Then, our analysis of of the homotopy group sheaf, presented above, implies the following:

$\mathbb{M}^{p.q}_{sb}(k(X)) \hookrightarrow \mathbb{M}^{p.q}_{nr}(X) \xleftarrow{\exists 1} \mathbb{M}^{p.q}(X)$

For $X \in Sm^{ft}_k$, we have the commutative diagram with horizontal group exact sequences, characterizing the embedding: $\mathbb{M}^{p.q}_{sb}(k(X))_i \subseteq \mathbb{M}^{p.q}_{nr}(X)$, together with the canonical homomorphism: $\mathbb{M}^{p.q}(X) \to \mathbb{M}^{p.q}_{nr}(X)$

$$\begin{array}{ccccccc}
0 & \longrightarrow & \mathbb{M}^{p.q}_{sb}(k(X)) & \longrightarrow & \mathbb{M}^{p.q}(k(X)) & \longrightarrow & \oplus_{v\,:\,\text{divisorial valuation on } k(X)/k}\mathbb{M}^{p-1.q-1}(\kappa(v)) \\
 & & \downarrow & & \| & & \downarrow \\
0 & \longrightarrow & \mathbb{M}^{p.q}_{nr}(X) & \longrightarrow & \mathbb{M}^{p.q}(k(X)) & \longrightarrow & \oplus_{y\in X^{(1)}}\mathbb{M}^{p-1.q-1}(\kappa(v_y)) \\
 & & \nwarrow{\scriptstyle \exists 1} & & \uparrow & & \\
 & & & & \mathbb{M}^{p.q}(X) & &
\end{array}$$

Likewise, for any motivic spectra $\mathbb{X}, \mathbb{Y}$, the associated Nisnevich sheaf to the presheaf

$$Sm^{ft}_k \to AbelianGroups$$
$$X \mapsto Hom_{\mathcal{SH}(k)}\left(\left(\Sigma^{\infty}_{\mathbb{P}^1}X_+\right) \wedge \mathbb{Y}, \mathbb{X}\right)$$

also defines an unramified sheaf, and we can consider its stably birationalized Nisnevich subsheaf, though this is not a new example for we may write

$$Hom_{\mathcal{SH}(k)}\left(\left(\Sigma^{\infty}_{\mathbb{P}^1}X_+\right) \wedge \mathbb{Y}, \mathbb{X}\right) \cong Hom_{\mathcal{SH}(k)}\left(\Sigma^{\infty}_{\mathbb{P}^1}X_+, Map_{\mathcal{SH}(k)}(\mathbb{Y}, \mathbb{X})\right)$$

by use of the mapping spectrum $Map_{\mathcal{SH}(k)}(\mathbb{Y}, \mathbb{X})$.

Of course, our Theorem and Definition 1.12 implies $\mathbb{M}^{p.q}_{sb}$ is a stable birational invariant Nisnevich sheaf on $Sm^{ft}_k$ and $\mathbb{M}^{p.q}_{ur}$, which is equal to $\mathbb{M}^{p.q}_{sb}$ for smooth proper $k$-schemes, is a stable birational invariant of smooth proper $k$-schemes.

Also, there are many exapmles of the generalized motivic cohomology theories (see e.g. [M04][VRO07][L16]). Here, we focus upon just a couple;

For the motivic cohomology theory $H^{*,*}$ [VSF00][MVW06], represented by the motivic Eilenberg-MacLane spectrum $\mathbf{H}\mathbb{Z}$ [VRO07, p.162, 3.1] [L16, p.392–393], we shall also call the resulting unramified $H^{*,*}$-theory Nisnevich sheaf $H^{*,*}_{ur}$ the *unramified motivic cohomology Nisnevich sheaf*. Here, $H^{n,n}_{ur}$ for $n \in \mathbb{Z}_{\geq 1}$ is, by [NS89] [T92][K09], the *unramified Milnor K-theory Nisnevich sheaf* $\underline{\mathbf{K}}^M_n$ in weight $n$, which goes back to Kato [K80][K86]. This is the universal example of the unramified sheaves obtained from Rost's cycle modules as in (75). Furthermore, any Rost's cycle module $M_*$ gives rise to a $\mathbb{P}^1$-spectrum [M04, Example 5.3].

For the universal motivic cohomotopy theory $\pi^{*,*}$, represented by the motivic sphere spectrum $\mathbb{S}$ [VRO07, p.155, Remark 2.14] [L16, p.391], $\pi^{n,n}_{ur}$ for $n \in \mathbb{Z}$ is, by [M04][M12], the *unramified Milnor-Witt K-theory Nisnevich sheaf* $\underline{\mathbf{K}}^{MW}_n$ in weight $n$. This is the universal example of the unramified sheaves obtained from Milnor-Witt cycle modules as in (77).

In this way, motivic (stable) homotopy theory offers an unlimited supply of stable birational invariants of smooth proper $k$-schemes.

(v) ([S20, Theorem 0.2] [K21b, Corollary 1.17]) So far, all of our unramified sheaves were $\mathbb{A}^1$-invariant. However, Saito [S20, Theorem 0.2] and Koizumi [K21b, Corollary 1.17] showed reciprocity sheaves in the sense of [KSY16], which are not necessarily $\mathbb{A}^1$-invariant, are also unramified. Examples of reciprocity sheaves include $\mathbb{A}^1$-invariant Nisnevich sheaves with transfer treated in (ii)(i), smooth commutative algebraic groups over $k$, e.g. the additive group $\mathbb{G}_a$, and the modules of absolute Kähler differentials $\Omega^j$ and de Rham-Witt differentials $W_n\Omega^j$, as was observed by Rülling [KSY16, Appendix].

For more on the reciptocity sheaves, consult [KSY16][S20][BRS20] [KSY22].

In spite of such rich examples of unramified sheaves, unramified sheaves are not necessarily birational sheaves, as a matter of course. Now, as is already quoted over and over again, our main theorem of this paper shows how to convert such rich unramified sheaves into birational sheaves.

**Theorem and Definition 4.5.** (Theorem and Definition 1.12)
(i) Given an unramified sheaf $S$ on $Sm_k^{ft}$ (resp. on $\widetilde{Sm_k^{ft}}$), let us set

$$S_{sb}\left(K/k\right) := \cap_{\substack{v, \text{ divisorial} \\ \text{valuation of } K/k}} S(V_v) \tag{80}$$

for any finitely generated field extension $K/k$. Then the correspondence

$$U \mapsto S_{sb}(U) := S_{sb}(k(U)/k) = \cap_{\substack{v, \text{ divisorial} \\ \text{valuation of } k(U)/k}} S(V_v)$$

$$\left( \subseteq S(U) := \cap_{\substack{v, \text{ divisorial valuation of } k(U)/k \\ \text{s.t. } V_v = \mathcal{O}_{U,x} \text{ for some } x \in U^{(1)}}} S(V_v) = \cap_{x \in U^{(1)}} S(\mathcal{O}_{U,x}) \subset S(k(U)) \right)$$

defines a birational subsheaf $S_{sb}$ of an unramified sheaf $S$ on $Sm_k^{ft}$ (resp. on $\widetilde{Sm_k^{ft}}$), which we call the <u>*stably birationalized subsheaf*</u> $S_{sb}$ of an unramfied sheaf $S$, or simply <u>**SBNR**</u>, because $S_{ab}$ is, once its presheaf $\left(Sm_k^{ft}\right)^{op} \to \mathcal{C}$ property is established, immediately seen to be a birational presheaf, which is actually a Nisnevich sheaf by [AM11, Lemm.6.1.2] [MV99, §3, Proposition 1.4] as was recalled in Remark 1.10 (i), and is staby birational invariant by [KS15, p.330, Appendix A by Jean-Louis Colliot-Thélène] as was recalled in Remark 1.10 (ii).
(ii) For any smooth proper $k$-scheme $X$,

$$S_{sb}(X) = S(X).$$

Consequently, any unramified sheaf $S$ is stably birational invariant among smooth proper $k$-schemes.

From Theorem 1.11 and the above Theorem and Definition 1.12, we immediately obtain the following;

**Corollary 4.6.** (Corollary 1.15) *Let $S$ be an unramified Zariski sheaf, and $S_{sb}$ be its associated stably birational Nisnevich sheaf. Then, for any $X \in Sm_k^{ft}$, or more generally, those* ind $-k$ *varietis described in Remark 1.5(ii), including $BG$ and $X_G = EG \times_G X$, if $X$ is $(-i)$-retract rational, then $S_{sb}(X)$ is a direct summand of $S_{sb}(Z^i)$ for some smooth proper $Z^i$ of dimension $i$.* □

## 5. Proofs

We first take care of Proposition 1.4:

*Proof of Proposition 1.4.* (i): Obvious from the definitions, where $Z^0$ is always taken to be Spec $k$.

(ii): Set $\nu := n$ for the first, and $\nu := N$ for the last two. Then the claim follows immediately from the explicit birational equivalence: $\mathbb{A}^{\nu-i} = \mathbb{A}^{\nu-j} \times \mathbb{A}^{j-i}$.

(iii): The first implication is obvious. For the second implication, suppose a stable $(-i)$-rational equivalence is given by a regular equivalence at dense open subsets as follows:

$$\mathbb{A}^{N-n} \times X \ \supseteq \ V \ \overset{f}{\underset{\cong}{\to}} \ W \ \subseteq \ \mathbb{A}^{N-n} \times \mathbb{A}^{n-i} \times Z^i$$

Then, apply the construction in Remark 1.7(ii) to form the following diagram:

$$X \underset{\text{open dense}}{\hookleftarrow} U \underset{r}{\overset{i}{\rightleftarrows}} V' \underset{f^{-1}}{\overset{f}{\rightleftarrows}} f(V') \underset{\text{open dense}}{\hookrightarrow} \mathbb{A}^{N-n} \times \mathbb{A}^{n-i} \times Z^i$$

Of course, this shows $X$ is retract $(-i)$-rational,

(iv): The stably birational invariance of stable $(-i)$-rationality is obvious. Thus, we concentrate in proving the stably birational invariance of retract $(-i)$-rationality. So, let us suppose irreducible $X, Y$ are stable rational equivalent, given by a rational equivalence induced by an honest morphism $f$:

$$\mathbb{A}^r \times X \supseteq V \overset{f}{\underset{\cong}{\to}} W \subseteq \mathbb{A}^s \times Y \tag{81}$$

Now, suppose $Y$ is retract $(-i)$-rational, then clearly $\mathbb{A}^s \times Y$ is also retract $(-i)$-rationa, given by some honest morphisms defined on some dense open subsets by Remark 1.7 (i):

$$\mathbb{A}^s \times Y \underset{\text{open dense}}{\hookleftarrow} W' \underset{r'}{\overset{i'}{\rightleftarrows}} T \underset{\text{open dense}}{\hookrightarrow} \mathbb{A}^{N+s-i} \times Z^i \ . \tag{82}$$

Then, combining (81) and (82), we obtatin the following commutative diagram which shows the retract $(-i)$-rationality of $\mathbb{A}^r \times X$ :
(83)

$$\mathbb{A}^r \times X \underset{\text{open dense}}{\hookleftarrow} f^{-1}\big(W \cap W'\big) \underset{f^{-1}}{\overset{f}{\rightleftarrows}} W \cap W' \underset{r'}{\overset{i'}{\rightleftarrows}} (r')^{-1}\big(W \cap W'\big) \underset{\text{open dense}}{\hookrightarrow} \mathbb{A}^{N+s-i} \times Z^i$$

Finally, apply the construction in Remark 1.7(ii) with $V = f^{-1}(W \cap W')$. Then, combining with (83), we get the following commutative diagram: Then, apply the construction in Remark 1.7(ii) to form the following diagram:

$$X \underset{\text{open dense}}{\hookleftarrow} U \underset{r}{\overset{i}{\rightleftarrows}} V' \underset{f^{-1}}{\overset{f}{\rightleftarrows}} f(V') \underset{r'}{\overset{i'}{\rightleftarrows}} (r')^{-1}\big(f(V')\big) \underset{\text{open dense}}{\hookrightarrow} \mathbb{A}^{N+s-i} \times Z^i$$

This commutative diagram implies the retract $(-i)$-rationality of $X$, as was desired. $\square$

*Proof of Theorem and Definition 1.12 (i) for the case $\widetilde{Sm_k^{ft}}$, i.e. Proof of the claim:*

*For an unramified presheaf $S$ on $\widetilde{Sm_k^{ft}}$,*

$$\widetilde{Sm_k^{ft}} \ni X \ \mapsto \ S_{sb}(X) := \bigcap_{\substack{v,\ \boxed{geometric} \\ rank\ 1\ discrete \\ valuation\ on\ k(X)/k}} S(V_v)$$

*defines a subpresheaf $S_{sb}$ of $S$ on $\widetilde{Sm_k^{ft}}$.*

Any smooth morphism $f : Y \to X$ in $\widetilde{Sm_k^{ft}}$ induces the corresponding extension of function fields

$$k(f) : \ E := k(X) \ \subset \ k(Y) =: F.$$

Now, $k(f)$ induces $S(k(f))$, which Morel [M12, Proof of Proposition 2.8]. (see also [C19, Proof of Proposition 1.2] ) showed to restrict to $S(f)$ as in (73), and which we must show to restrict further to $S_{sb}\,(k(f))$ as the following commutative diagram:

$$\begin{array}{ccccc}
S\,(k(X)) & \hookleftarrow & S(X) & \hookleftarrow & S_{sb}\,(k(X)/k) \\
\Big\downarrow & & \Vert & & \Vert \\
 & & \bigcap_{x \in X^{(1)}} S\,(\mathcal{O}_{X,x}) & & \bigcap_{\substack{w,\ \boxed{\text{geometric}}\ \text{rank 1} \\ \text{discrete valuation on } k(X)/k}} S(V_w) \\
S(k(f)) & & \Big\downarrow \exists S(f) \text{ by Morel (73)} & & \Big\downarrow \exists S_{sb}(k(f))\ ? \\
 & & \bigcap_{y \in Y^{(1)}} S\,(\mathcal{O}_{Y,y}) & & \bigcap_{\substack{v,\ \boxed{\text{geometric}}\ \text{rank 1} \\ \text{discrete valuation on } k(Y)/k}} S(V_v) \\
\Big\downarrow & & \Vert & & \Vert \\
S\,(k(Y)) & \hookleftarrow & S(Y) & \hookleftarrow & S_{sb}\,(k(Y)/k)
\end{array}$$

For this purpose, it suffices to show, for each geometric rank 1 discrete valuation $v$ on $k(Y)/k$,

$$S(k(f)) \left( \bigcap_{\substack{w,\ \boxed{\text{geometric}}\ \text{rank 1} \\ \text{discrete valuation on } k(X)/k}} S(V_w) \right) \subseteq S(V_v). \tag{84}$$

Then, since the geometric rank 1 disrete valuation $v$ on $k(Y)/k$ is either trivial on $k(X)/k$ or restricts to a geometric discrete valation $w$ on $k(X)/k$ of rank 1 by [M08, Prop.1.4] [ZS60b, Chap.VI,§6,Lem.2,Cor.1], we have either one of the following commutative diagrams of essentially smooth $k$-algebras:

(85)
$$\begin{array}{ccc} k(X) & = & k(X) \\ {\scriptstyle S(f)}\downarrow & & \downarrow \\ k(Y) & \hookleftarrow & V_v \end{array} \qquad \begin{array}{ccc} k(X) & \hookleftarrow & V_w \\ {\scriptstyle S(f)}\downarrow & & \downarrow \\ k(Y) & \hookleftarrow & V_v \end{array}$$

Now, with this Morel's extension [M12, Proposition 2.8] of an unramified $\widetilde{\mathcal{F}_k}$-datum $S$ to an umramified sheaf on $\widetilde{Sm_k^{ft}}$ and its continuous extension to essential smooth $k$-schemes, guaranteed by **D1**, at hand, we may apply this extended $S$ to (85) to justify either one of the following commutative diagrams:

$$\begin{array}{ccc} S(k(X)) & = & S(k(X)) \\ {\scriptstyle S(k(f))}\downarrow & & \downarrow \\ S(k(Y)) & \hookleftarrow & S(V_v) \end{array} \qquad \begin{array}{ccc} S(k(X)) & \hookleftarrow & S(V_w) \\ {\scriptstyle S(k(f))}\downarrow & & \downarrow \\ S(k(Y)) & \hookleftarrow & S(V_v) \end{array}$$

Of course, this justifies (84), and hence, completes the proof.

*Proof of Theorem and Definition 1.12 (i) for the case $Sm_k^{ft}$, i.e. Proof of the claim:*

*For an unramified sheaf $S$ on $Sm_k^{ft}$,*

$$Sm_k^{ft} \ni X \mapsto S_{sb}(X) := \bigcap_{\substack{v,\ \boxed{geometric} \\ rank\ 1\ discrete \\ valuation\ on\ k(X)/k}} S(V_v)$$

*defines a birational subsheaf $S_{sb}$ of $S$ on $Sm_k^{ft}$.*

It suffices to show $S(f)$ restricts to $S_{sb}(f)$ for any morphism $f : X \to Y$ in $Sm_k$. Now, by Morel's proof of Theorem 4.3, it suffices to prove it for the cases of $f$ being smooth and $f$ being a codimension 1 regular immersion. However, the case of $f$ being smooth was just taken care of. Thus, it suffices to show $S(i) : S(X) \to S(Y)$ restricts to a map $S_{sb}(i) : S_{sb}(X) \to S_{sb}(Y)$, for any codimension 1 regular immersion $i : \overline{\{y\}} = Y \hookrightarrow X$ $(y \in X^{(1)})$ in $Sm_k$, inducing the canonical quotient

$$\pi_y : \mathcal{O}_{X,y} \twoheadrightarrow \kappa(y) := \mathcal{O}_{X,y}/\mathfrak{m}_{X,y}, \quad \text{a morphism in } Pro^{Aff}\left(Sm_k^{ft}\right).$$

More explicitly, $S(\pi_y)$, which Morel [M12, Lemma 2.12] (see also [C19, Theorem 1.3] ) showed to be identified with the <u>*specialization map*</u> $s_{v_y}$, where $v_y$ is the rank 1 geometric valuation corresponding to the discrete valuation ring $\mathcal{O}_{X,y}$, and to restrict to $S(i)$, as we recalled in (74) in our summary of Morel's proof of Theorem 4.3 (ii).

We then have to show $S(i)$ furhter restricts to $S_{sb}(i)$ :

$$\begin{array}{ccccc}
S(\mathcal{O}_{X,y}) & \hookleftarrow & S(X) & \hookleftarrow & S_{sb}(k(X)/k) \\
 & & \| & & \| \\
 & & \bigcap_{y' \in X^{(1)}} S(\mathcal{O}_{X,y'}) & & \bigcap_{w,\ \boxed{\text{geometric}}\ \text{rank 1 discrete valuation on } k(X)/k} S(V_w) \\
\downarrow s_{v_y} & & \vdots\ \exists S(i)\ \text{by Morel (74)} & & \vdots\ \exists S_{sb}(i)\ ? \\
 & & \bigcap_{z \in Y^{(1)}} S(\mathcal{O}_{Y,z}) & & \bigcap_{v,\ \boxed{\text{geometric}}\ \text{rank 1 discrete valuation on } k(Y)/k} S(V_v) \\
 & & \| & & \| \\
S(\kappa(y)) = S(k(Y)) & \hookleftarrow & S(Y) & \hookleftarrow & S_{sb}(k(Y)/k)
\end{array}$$

For this purpose, it suffices to show, for each geometric rank 1 discrete valuation $v$ on $k(Y)/k$,

(86)
$$s_{v_y}\left( \bigcap_{w,\ \boxed{\text{geometric}}\ \text{rank 1 discrete valuation on } k(X)/k} S(V_w) \right) \subseteq S(V_v).$$

Now, since $\mathcal{O}_{X,y}$ is the valuation ring of a rank 1 discrete valuation $v_y$ on $k(X)/k$ and $v$ is a rank 1 discrete valuation ring of $k(Y)/k = \kappa(v_y)/k$, we have the following commutative diagram to define the valuation ring of the composite valuation $v_y \circ v$ as in (27):

$$\begin{array}{ccccc}
k(X) & \hookleftarrow & \mathcal{O}_{X,y} = V_{v_y} & \hookleftarrow & V_{v_y \circ v} \\
 & & \downarrow \pi(v_y) & & \downarrow \\
 & & \kappa(y) = k(Y) = \kappa(v_y) & \hookleftarrow & V_v
\end{array}$$

Whereas the valuation ring $V_{v_y \circ v}$ of the rank 2 geometric valuation $v_y \circ v$ might not belong to $Pro^{Aff}(Sm_k^{ft})$, we may still apply Corollary 60 to obtain the following commutative diagram in $Pro^{Aff}\left(Sm_k^{ft}\right)$ :

(87)
$$\begin{array}{ccccccc}
k(X) & \hookleftarrow & \mathcal{O}_{X,y} = V_{v_y} = R_{\mathfrak{m}_{v_y} \cap R} & \hookleftarrow & R_{\mathfrak{m}_{v_y \circ v} \cap R} & \xrightarrow{\cong} & \mathcal{O}_{L, \mathfrak{m}_{v_y \circ v} \cap R} \\
 & & \downarrow \pi(v_y) & & \downarrow & & \\
 & & \kappa(y) = k(Y) = \kappa(v_y) = Frac\left(\frac{R}{\mathfrak{m}_{v_y} \cap R}\right) & \hookleftarrow & V_v = \left(\frac{R}{\mathfrak{m}_{v_y} \cap R}\right)_{\overline{\left(\mathfrak{m}_{v_y \circ v} \cap R\right)}} & &
\end{array}$$

Here, we have used the following notations and remarks:

- $R \subseteq V_{v_y \circ v}$ is a normal integral subdomain which is finitely generated over $k$ with $\mathfrak{m}_{v_y \circ v} \cap R \subset R$ its height 2 prime ideal, making $R_{\mathfrak{m}_{v_y \circ v} \cap R}$ a Noetherian local ring of dimension 2;
- $L$ is the smooth locus of Spec $R$, which is of finite type over $k$;
- Since $L$ contains $\mathfrak{m}_{v_y \circ v} \cap R$, we have an isomorphism of regular local rings:

$$R_{\mathfrak{m}_{v_y \circ v} \cap R} \xrightarrow{\cong} \mathcal{O}_{L, \mathfrak{m}_{v_y \circ v} \cap R} :$$

Then we may apply $S$ to the commutative diagram (87) in $Pro^{Aff}\left(Sm_k^{ft}\right)$ to obtain the following commutative diagram:

(88)

$$\begin{array}{ccccc}
S(k(X)) \hookleftarrow S(V_{v_y}) = S\left(R_{\mathfrak{m}_{v_y} \cap R}\right) & \hookleftarrow & S\left(R_{\mathfrak{m}_{v_y \circ v} \cap R}\right) & \xrightarrow{\cong} & S\left(\mathcal{O}_{L, \mathfrak{m}_{v_y \circ v} \cap R}\right) \\
\downarrow S(\pi(v_y)) & & \downarrow & & \\
S(\kappa(v_y)) = S\left(Frac\left(\frac{R}{\mathfrak{m}_{v_y} \cap R}\right)\right) & \hookleftarrow & S(V_v) = S\left(\left(\frac{R}{\mathfrak{m}_{v_y} \cap R}\right)_{\left(\mathfrak{m}_{v_y \circ v} \cap R\right)}\right) & &
\end{array}$$

Thus, since $S\left(\pi(v_y)\right) = s_{v_y}$, we see from (88) the following implications:

(89)

$$s_{v_y}\left(\bigcap_{\substack{w,\ \text{divisorial valuation} \\ \text{on } k(X)/k}} S(V_w)\right) \subseteq s_{v_y}\left(\bigcap_{\substack{w,\ \text{divisorial valuation on } k(X)/k = k(L)/k \\ \text{s.t. } V_w = \mathcal{O}_{L,x} \text{ for some } x \in L^{(1)}}} S(V_w)\right) = s_{v_y}\left(S(L)\right)$$
$$\subseteq s_{v_y}\left(S\left(\mathcal{O}_{L, \mathfrak{m}_{v_y \circ v} \cap R}\right)\right) = s_{v_y}\left(S\left(R_{\mathfrak{m}_{v_y \circ v} \cap R}\right)\right) \subseteq S(V_v).$$

Consequently we have shown (86), and the proof is complete.

*Proof of Theorem and Definition 1.12 (ii), i.e. Proof of the claim:*

*For an unramified sheaf $S$ on $Sm_k^{Sm}$, the following hold:*
(1) *for any proper $X \in Sm_k$,*

$$S_{sb}(X) = S(X)$$

(2) *$S$ is a stable birational invariant of proper smooth $k$ schemes of finite type.*

Of course, (2) would follow from Theorem and Definition 1.12 (i), which we have just proved above, and (1). Thus, it suffices to prove (1). However, (1) follows from (the proof of) [CT95, Proposition 2.1.8 (e)]. Thus the proof is complete. □

## References

[A56a] Shreeram Abhyankar, *On the valuations centered in a local domain,* Amer. J. Math. 78 (1956), 321–348.

[A56b] Shreeram Abhyankar, *Local uniformization on algebraic surfaces over ground fields of char-acter istic $p \neq 0$*, Ann. of Math. 63 (1956), 491–526.

[A74] Michel André, *Homologie des algèbres commutatives*, Die Grundlehren der mathematischen Wissenschaften, Band 206. Springer-Verlag, Berlin-New York, 1974. xv+341 pp.

[AC12] C. Araujo and A-M. Castravet, *Polarized minimal families of rational curves and higher Fano manifolds*, American J. Math., 134(1) (2012), 87-107.

[AGV73] M. Artin, A. Grothendieck, J.-L. Verdier avec la participation de N. Bourbaki, P. Deligne, B. Saint-Donat, *Théorie des Topos et Cohomologie Étale des Schémas. Seminaire de Geometrie Algebrique du Bois-Marie 1963-1964 (SGA 4): Tome 3*, Lecture Notes in Mathematics, Vol. 305. Springer-Verlag, Berlin-New York, 1973. vi+640 pp.

[A86a] M. Artin, *Néron models,* Arithmetic geometry (Storrs, Conn., 1984), 213–230, Springer, New York, 1986.

[A86b] M. Artin, *Lipman's proof of resolution of singularities for surfaces,* Arithmetic geometry (Storrs, Conn., 1984), 267–287, Springer, New York, 1986.

[AM72] M. Artin, D. Mumford, *Some elementary examples of unirational varieties which are not rational,* Proc. London Math. Soc. (3) 25 (1972), 75–95.

[AM11] Aravind Asok, Fabien Morel, *Smooth varieties up to* $\mathbb{A}^1$*-homotopy and algebraic h-cobordisms*, Adv. Math. 227 (2011), no. 5, 1990–2058.

[A13] Aravind Asok, *Rationality problems and conjectures of Milnor and Bloch-Kato,* Compos. Math. 149 (2013), no. 8, 1312–1326.

[AH62] Michael Francis Atiyah, Friedrich Ernst Peter Hirzebruch, *Analytic cycles on complex manifolds,* Topology 1 (1962), 25–45.

[AM69] Michael Francis Atiyah, Ian G. Macdonald, *Introduction to commutative algebra,* Addison-Wesley Publishing Co., Reading, Mass.-London-Don Mills, Ont. 1969, ix+128 pp.

[AB17] Asher Auel, Marcello Bernardara, *Cycles, derived categories, and rationality,* Surveys on recent developments in algebraic geometry, 199–266, Proc. Sympos. Pure Math., 95, Amer. Math. Soc., Providence, RI, 2017.

[B21] Tom Bachmann, *The zeroth* $\mathbb{P}^1$*-stable homotopy sheaf of a motivic space*, J. Inst. Math. Jussieu (2021), 1–25.

[B16] Arnaud Beauville, *The Lüroth problem,* Rationality problems in algebraic geometry, 1–27, Lecture Notes in Math., 2172, Fond. CIME/CIME Found. Subser., Springer, Cham, 2016.

[BRS20] Federico Binda, Kay Rülling, Shuji Saito, On the cohomology of reciprocity sheaves , **arXiv:2010.03301**.

[BO74] Spencer Bloch, Arthur Ogus, *Gersten's conjecture and the homology of schemes,* Ann. Sci. École Norm. Sup. (4) 7 (1974), 181–201.

[BS83] Spencer J. Bloch, Vasudevan Srinivas, *Remarks on correspondences and algebraic cycles,* Amer. J. Math. 105 (1983), no. 5, 1235–1253.

[BK85] F. A. Bogomolov, P. I. Katsylo, *Rationality of some quotient varieties,* Math. USSR Sbornik TOM 126(168) (1985), Bun. 4 Vol. 54(1986), No. 2, 571–576.

[B87] F. A. Bogomolov, *The Brauer group of quotient spaces of linear representations,* Izv. Akad. Nauk SSSR Ser. Mat. 51 (1987), 485–516, 688.

[B88] F. A. Bogomolov, *The Brauer group of quotient spaces of linear representations*, (Russian) Izv. Akad. Nauk SSSR Ser. Mat. 51 (1987), no. 3, 485–516, 688; translation in Math. USSR-Izv. 30 (1988), no. 3, 455–485.

[B89] F. A. Bogomolov, *Brauer groups of the fields of invariants of algebraic groups,* Mat. Sb. 180 (1989), 279–293.

[BLR90] Siegfried Bosch, Werner Lütkebohmert, Michel Raynaud, *Néron models*, Ergebnisse der Mathematik und ihrer Grenzgebiete (3) 21, Springer-Verlag, Berlin, 1990, x+325 pp.

[BDPP13] Sébastien Boucksom, Jean-Pierre Demailly, Mihai Păun, Thomas Peternell, *The pseudo-effective cone of a compact Kähler manifold and varieties of negative Kodaira dimension,* J. Algebraic Geom. 22 (2013), no. 2, 201–248.

[B72] Nicolas Bourbaki, *Elements of mathematics. Commutative algebra,* Translated from the French. Hermann, Paris; Addison-Wesley Publishing Co., Reading, Mass., 1972. xxiv+625 pp.

[C92] F. Campana, *Connexité rationnelle des variétés de Fano,* Ann. Sci. École Norm. Sup. (4) 25 (1992), no. 5, 539–545.

[CDP15] F. Campana, J.-P. Demailly, Th. Peternell, *Rationally connected manifolds and semipositivity of the Ricci curvature.* Recent advances in algebraic geometry, 71–91, London Math. Soc. Lecture Note Ser., 417, Cambridge Univ. Press, Cambridge, 2015.

[C19] Robin Carlier, *Notes: Unramified Sheaves of Set,* 23 octobre 2019, 13pp. **https://deglise.perso.math.cnrs.fr/docs/2019/Morel/Carlier1.pdf**

[CJ19] Gunnar Carlsson, Roy Joshua, *Atiyah-Segal Derived Completions for Equivariant Algebraic G-Theory and K-Theory*, 39pp, **arXiv:1906.06827v2**

[CL17] Andre Chatzistamatiou, Marc Levine, *Torsion orders of complete intersections,* Algebra Number Theory 11 (2017), no. 8, 1779–1835.

[CGR06] V. Chernousov, P. Gille, Z. Reichstein, *Resolving G-torsors by abelian base extensions,* J. Algebra 296 (2006), no. 2, 561–581.

[CMS10] Suyoung Choi, Mikiya Masuda, Dong Youp Suh, Topological classification of generalized Bott towers. Trans. Amer. Math. Soc. 362 (2010), no. 2, 1097–1112.

[C08] Utsav Choudhury, *Homotopy theory of schemes and $\mathbb{A}^1$-fundamental groups,* ALGANT Erasmus Mundus Master Thesis, June 2008, 70pp.
**https://algant.eu/documents/theses/choudhury.pdf**

[CTHK97] Jean-Louis Colliot-Thélène, Raymond T. Hoobler, Bruno Kahn, *The Bloch-Ogus-Gabber theorem,* Algebraic K-theory (Toronto, ON, 1996), 31–94, Fields Inst. Commun., 16, Amer. Math. Soc., Providence, RI, 1997.

[CTO89] Jean-Louis Colliot-Thélène, Manuel Ojanguren, *Variétés unirationnelles non rationnelles: au-delà de l'exemple d'Artin et Mumford,* Invent. Math. 97 (1989), no. 1, 141–158.

[CT95] Jean-Louis Colliot-Thélène, *Birational invariants, purity and the Gersten conjecture, K*-theory and algebraic geometry: connections with quadratic forms and division algebras (Santa Barbara, CA, 1992), 1–64, Proc. Sympos. Pure Math., 58, Part 1, Amer. Math. Soc., Providence, RI, 1995.

[CTS07] Jean-Louis Colliot-Thélène, Jean-Jacques Sansuc, *The rationality problem for fields of invariants under linear algebraic groups (with special regards to the Brauer group). Algebraic groups and homogeneous spaces,* 113–186, Tata Inst. Fund. Res. Stud. Math., 19, Tata Inst. Fund. Res., Mumbai, 2007.

[CTV12] Jean-Louis Colliot-Thélène, Claire Voisin, *Cohomologie non ramifiée et conjecture de Hodge entière*, Duke Math. J. 161 (2012), no. 5, 735–801.

[CLO12] Brian Conrad, Max Lieblich, Martin Olsson, *Nagata compactification for algebraic spaces*, J. Inst. Math. Jussieu 11 (2012), no. 4, 747–814.

[C21] Steven Dale Cutkosky, *Local Uniformization of Abhyankar Valuations,* Michigan Math. J. Advance Publication 1–33, 2021.

[D06] Frédéric Déglise, *Transferts sur les groupes de Chow à coefficients,* Math. Z. 252 (2006), no. 2, 315–343.

[D11] F. Déglise, *Modules homotopiques,* Doc. Math. 16 (2011), 411–455.

[DFJK21] Frédéric Déglise, Jean Fasel, Fangzhou Jin, Adeel A. Khan, *On the rational motivic homotopy category*, J. Éc. polytech. Math., 8 (2021), 533–583.

[dJS07] A. J. de Jong and J. Starr, *Higher Fano manifolds and rational surfaces,* Duke Math. J. 139 (2007), no. 1, 173–183.

[D87] I. V. Dolgachev, *Rationality of fields of invariants,* in: Algebraic Geometry, Bowdoin, 1985, Brunswick, Maine, 1985, Amer. Math. Soc., Providence, RI, 1987, pp. 3–16.

[EG98] Dan Edidin, William Graham, *Equivariant intersection theory,* Invent. Math. 131 (1998), no. 3, 595–634.

[EKW21] Elden Elmanto, Girish Kulkarni, Matthias Wendt, *$\mathbb{A}^1$-connected components of classifying spaces and purity for torsors*, **arXiv:2104.06273**.

[EM73] S. Endo, T. Miyata, *Invariants of finite abelian groups,* J. Math. Soc. Japan 25 (1973) 7–26.

[F20] Niels Feld, *Feld, Niels (F-GREN-IF) Milnor-Witt cycle modules*, J. Pure Appl. Algebra 224 (2020), no. 7, 106298, 44 pp.

[F21a] Niels Feld, *Morel homotopy modules and Milnor-Witt cycle modules*, Doc. Math. 26 (2021), 617–659.

[F21b] Niels Feld, *Milnor-Witt homotopy sheaves and Morel generalized transfers,* Adv. Math. 393 (2021), Paper No. 108094.

[F915] E. Fischer, *Die Isomorphie der Invariantenkörper der endlichen Abelschen Gruppen linearen Transformationen,* Nachr. König. Ges. Wiss. Göttingen (1915), 77–80.

[F98] William Fulton, Intersection Theory, Second edition, Springer, 1998.

[G94] Ofer Gabber, *Gersten's conjecture for some complexes of vanishing cycles,* Manuscripta Math. 85 (1994), no. 3–4, 323–343.

[GMS03] Skip Garibaldi, Alexander Merkurjev, Jean-Pierre Serre, *Cohomological invariants in Galois cohomology,* University Lecture Series, 28. American Mathematical Society, Providence, RI, 2003. viii+168 pp.

[GK94] M. Grossberg, Y. Karshon, Bott towers, complete integrability, and the extended character of representations, Duke Math. J., 76 (1994), 23–58.

[G60] Alexander Grothendieck, *Éléments de géométrie algébrique. I. Le langage des schémas*, Inst. Hautes Études Sci. Publ. Math., No. 4 (1960), 5–228.

[G64] Alexander Grothendieck, *Éléments de géométrie algébrique. IV. Étude locale des schémas et des morphismes de schémas, I*, Inst. Hautes Études Sci. Publ. Math., No. 20 (1964), 5–259.

[G65] Alexander Grothendieck, *Éléments de géométrie algébrique. IV. Étude locale des schémas et des morphismes de schémas, II*, Inst. Hautes Études Sci. Publ. Math., No. 24 (1965), 5–231.

[G66] Alexander Grothendieck, *Éléments de géométrie algébrique. IV. Étude locale des schémas et des morphismes de schémas, III*, Inst. Hautes Études Sci. Publ. Math., No. 28 (1966), 5–255.

[G67] Alexander Grothendieck, *Éléments de géométrie algébrique. IV. Étude locale des schémas et des morphismes de schémas, IV*, Inst. Hautes Études Sci. Publ. Math., No. 32 (1967), 5–361.

[G68] Alexander Grothendieck, *Le groupe de Brauer. I, Algèbres d'Azumaya et interprétations diverses; Le groupe de Brauer. II. Théorie cohomologique, Dix exposés sur la cohomologie des schémas; Le groupe de Brauer. III. Exemples et compléments, Dix exposés sur la cohomologie des schémas*, Adv. Stud. Pure Math. 46–188, North-Holland, Amsterdam, 1968.

[GD71] A. Grothendieck, J. A. Dieudonné, *Éléments de géométrie algébrique. I.* Grundlehren der mathematischen Wissenschaften , 166, Springer-Verlag, Berlin, 1971. ix+466 pp.

[H77] R. Hartshorne,, Algebraic Geometry, Graduate Texts in Mathematics 52, Springer Verlag, New York, 1977.

[HW19] Christian Haesemeyer,Charles A.Weibel, *The norm residue theorem in motivic cohomology,* Annals of Mathematics Studies, 200. Princeton University Press, Princeton, NJ, 2019. xiii+299 pp.

[H64] Heisuke Hironaka, *Resolution of singularities of an algebraic variety over a field of characteristic zero. I, II,* Ann. of Math. (2) 79 (1964), 109–203; ibid. (2) 79 1964 205–326.

[H14] Akinari Hoshi, *Rationality problem for quasi-monomial actions,* Algebraic number theory and related topics 2012, 203-227, RIMS Kôkyûroku Bessatsu B51, Res. Inst. Math. Sci. (RIMS), Kyoto, (2014).

[H20] Akinari Hoshi, *Noether's problem and rationality problem for multiplicative invariant fields: a survey,* Algebraic Number Theory and Related Topics 2016, 29–53, RIMS Kôkyûroku Bessatsu, B77, Res. Inst. Math. Sci. (RIMS), Kyoto, 2020.

[HKY20] Akinari Hoshi, Ming-chang Kang, Aiichi Yamasaki, *Degree three unramified cohomology groups and Noether's problem for groups of order* 243, Journal of Algebra 544 (2020) 262–301.

[HKK17] Annette Huber, Stefan Kebekus, Shane Kelly, *Differential forms in positive characteristic avoiding resolution of singularities*, Bull. Soc. Math. France 145 (2017), no. 2, 305–343.

[HM04] J.-M. Hwang and N. Mok, *Birationality of the tangent map for minimal rational curves*, Asian J. Math. 8 (2004), no.1, 51–63.

[IWX20] Daniel C. Isaksen, Guozhen Wang, Zhouli Xu, *Stable homotopy groups of spheres,* Proc. Natl. Acad. Sci. USA 117 (2020), no. 40, 24757–24763

[JLY02] Christian U. Jensen, Arne Ledet, Noriko Yui, *Generic polynomials. Constructive aspects of the inverse Galois problem,* Mathematical Sciences Research Institute Publications, 45. Cambridge University Press, Cambridge, 2002. x+258 pp.

[KN14] Bruno Kahn, Nguyen Thi Kim Ngan, *Sur l'espace classifiant d'un groupe algébrique linéaire, I,* J. Math. Pures Appl. (9) 102 (2014), no. 5, 972–1013.

[KN16] Bruno Kahn, Nguyen Thi Kim Ngan, *Modules de cycles et classes non ramifiées sur un espace classifiant,* Algebr. Geom. 3 (2016), no. 3, 264–295.

[KSY16] Bruno Kahn, Shuji Saito, Takao Yamazaki, *Reciprocity sheaves.* With two appendices by Kay Rülling, Compos. Math. 152 (2016), no. 9, 1851–1898.

[KSY22] Bruno Kahn, Shuji Saito, Takao Yamazaki, *Reciprocity sheaves, II.* Homology Homotopy Appl. 24 (2022), no. 1, 71–91.

[KS15] Bruno Kahn, Ramdorai Sujatha, *Birational geometry and localisation of categories. With appendices by Jean-Louis Colliot-Thélène and Ofer Gabber,* Doc. Math. 2015, Extra vol.: Alexander S. Merkurjev's sixtieth birthday, 277–334

[KS16] Bruno Kahn, Ramdorai Sujatha, *Birational motives, I: pure birational motives,* Ann. K-Theory 1, 379–440 (2016).

[KS17] Bruno Kahn, Ramdorai Sujatha, *Birational motives, II: Triangulated birational motives.* Int. Math. Res. Not. IMRN 2017, no. 22, 6778–6831.

[KOY21] Wataru Kai, Shusuke Otabe, Takao Yamazaki, *Unramified logarithmic Hodge-Witt cohomology and* $\mathbb{P}^1$*-invariance*, **arXiv:2105.07433**.

[K80] Kazuya Kato, *A generalization of local class field theory by using* $K$*-groups. II,* J. Fac. Sci.Univ. Tokyo Sect. IA Math. 27 (3) (1980), 603–683.

[K86] Kazuya Kato, *Milnor* $K$*-theory and the Chow group of zero cycles,* In Applications of algebraic $K$-theory to algebraic geometry and number theory (Boulder, Colo., 1983), Part I, Contemp. Math. 55, Amer. Math. Soc., Providence, RI, 1986, 241–253.

[K09] Moritz Kerz, *The Gersten conjecture for Milnor* $K$*-theory,* Invent. Math. 175 (2009), no. 1, 1–33.

[KK05] Hagen Knaf, Franz-Viktor Kuhlmann, *Abhyankar places admit local uniformization in any characteristic,* Ann. Sci. École Norm. Sup. (4) 38 (2005), no. 6, 833–846.

[K21a] Junnosuke Koizumi, *Zeroth* $\mathbb{A}^1$*-homology of smooth proper varieties*, **arXiv:2101.04951**.

[K21b] Junnosuke Koizumi, *Steinberg symbols and reciprocity sheaves*, **arXiv:2108.04163**.

[KMM92] János Kollár, Yoichi Miyaoka, Shigefumi Mori, *Rationally connected varieties,* J. Algebraic Geom. 1 (1992), no. 3, 429–448.

[K92] János Kollár, *"Trento examples"*, in Ballico, E.; Catanese, F.; Ciliberto, C. (eds.), Classification of irregular varieties, Lecture Notes in Math., vol. 1515, Springer, 1992, p. 134.

[K95] János Kollár, *Nonrational hypersurfaces,* J. Amer. Math. Soc. 8 (1995), no. 1, 241–249.

[K96] János Kollár, *Rational curves on algebraic varieties*, Ergebnisse der Mathematik und ihrer Grenzgebiete. 32, Springer-Verlag, 1996.

[KM98] János Kollár, Shigefumi Mori, *Birational geometry of algebraic varieties,* With the collaboration of C. H. Clemens and A. Corti. Translated from the 1998 Japanese original. Cambridge Tracts in Mathematics, 134. Cambridge University Press, Cambridge, 1998. viii+254 pp.

[K18] Amalendu Krishna, *The completion problem for equivariant* $K$*-theory,* J. Reine Angew. Math. 740 (2018), 275–317.

[L74] H.W. Lenstra Jr., *Rational functions invariant under a finite abelian group,* Invent. Math. 25 (1974) 299–325.

[L16] Marc Levine, *An overview of motivic homotopy theory,* Acta Math. Vietnam. 41 (2016), no. 3, 379–407.

[L05] Martin Loren, *Multiplicative invariant theory*, Encyclopaedia of Mathematical Sciences, 135. Invariant Theory and Algebraic Transformation Groups, VI. Springer-Verlag, Berlin, 2005. xii+177 pp.

[L875] J. Lüroth, *Beweis eines Satzes über rationale Curven*, Math. Ann. (1875), no. 2, 163–165.

[MR10] Javier Majadas, Antonio G. Rodicio, *Smoothness, regularity and complete intersection*, London Mathematical Society Lecture Note Series, 373. Cambridge University Press, Cambridge, 2010. vi+134 pp.

[M80] Hideyuki Matsumura, *Commutative algebra, Second edition,* Mathematics Lecture Note Series, 56. Benjamin/Cummings Publishing Co., Inc., Reading, Mass., 1980. xv+313 pp.

[M89] Hideyuki Matsumura, *Commutative ring theory,* Translated from the Japanese by M. Reid. Second edition. Cambridge Studies in Advanced Mathematics, 8. Cambridge University Press, Cambridge, 1989. xiv+320 pp.

[MVW06] Carlo Mazza, Vladimir Voevodsky, Charles Weibel, *Lecture notes on motivic cohomology,* Clay Mathematics Monographs, 2. American Mathematical Society, Providence, RI; Clay Mathematics Institute, Cambridge, MA, 2006. xiv+216 pp.

[M08] Alexander Merkurjev, *Unramified elements in cycle modules,* J. Lond. Math. Soc. (2) 78 (2008), no. 1, 51–64.

[M17] Alexander S. Merkurjev, *Invariants of algebraic groups and retract rationality of classifying spaces.* Algebraic groups: structure and actions, 277–294, Proc. Sympos. Pure Math., 94, Amer. Math. Soc., Providence, RI, 2017.

[M70] John Milnor, *Algebraic* $K$*-theory and quadratic forms*, Invent. Math. 9 (1970), 318–344.

[M19] Norihiko Minami, *Higher uniruledness, Bott towers and "Higher Fano Manifolds"*, RIMS Kôkyûroku No.2135, The theory of transformation groups and its applications, 20, 15pp, 2019. **https://www.kurims.kyoto-u.ac.jp/~kyodo/kokyuroku/contents/pdf/2135-20.pdf**

[M20] Norihiko Minami, *From Ohkawa to strong generation via approximable triangulated categories—a variation on the theme of Amnon Neeman's Nagoya lecture series. (English summary) Bousfield classes and Ohkawa's theorem,* 17–88, Springer Proc. Math. Stat., 309, Springer, Singapore, 2020,

[MM86] Yoichi Miyaoka, Shigefumi Mori, *A numerical criterion for uniruledness*, Ann. of Math. (2), 124, 1986, 1, 65–69.

[MV99] F. Morel and V. Voevodsky, *$\mathbb{A}^1$-homotopy theory of schemes,* I. H. E. S Publ. Math., 90, (1999), 45–143 (2001).

[M04] Fabien Morel, *An introduction to $\mathbb{A}^1$-homotopy theory,* Contemporary developments in algebraic $K$-theory, 357–441, ICTP Lect. Notes, XV, Abdus Salam Int. Cent. Theoret. Phys., Trieste, 2004.

[M05a] Fabien Morel, *Milnor's conjecture on quadratic forms and mod 2 motivic complexes,* Rend. Sem. Mat. Univ. Padova 114 (2005), 63–101.

[M05b] Fabien Morel, *The stable $\mathbb{A}^1$-connectivity theorems,* K-Theory 35 (2005), no. 1–2, 1–68.

[M12] Fabien Morel, *$\mathbb{A}^1$-algebraic topology over a field.* Lecture Notes in Mathematics, 2052. Springer, Heidelberg, 2012. x+259 pp.

[M79] Shigefumi Mori, *Projective manifolds with ample tangent bundles*, Ann. of Math. 110, 593–606, 1979.

[NS89] Yu. P. Nesterenko, A. A. Suslin, *Homology of the general linear group over a local ring, and Milnor's $K$-theory,* (Russian) Izv. Akad. Nauk SSSR Ser. Mat. 53 (1989), no. 1, 121–146; translation in Math. USSR-Izv. 34 (1990), no. 1, 121–145.

[NS14] Josnei Novacoski, Mark Spivakovsky, *Reduction of local uniformization to the rank one case*, Valuation theory in interaction, 404–431, EMS Ser. Congr. Rep., Eur. Math. Soc., Zürich, 2014.

[NS16] Josnei Novacoski, Mark Spivakovsky, *On the local uniformization problem* , Algebra, logic and number theory, 231–238, Banach Center Publ., 108, Polish Acad. Sci. Inst. Math., Warsaw, 2016.

[O94] Tetsushi Ogoma, *General Néron desingularization based on the idea of Popescu*, J. Algebra 167 (1994), no. 1, 57–84.

[P93] Emmanuel Peyre, *Unramified cohomology and rationality problems,* Math. Ann. 296 (1993), no. 2, 247–268.

[P18] Alena Pirutka, *Varieties that are not stably rational, zero-cycles and unramified cohomology,* Algebraic geometry: Salt Lake City 2015, 459–483, Proc. Sympos. Pure Math., 97.2, Amer. Math. Soc., Providence, RI, 2018.

[P17] Bjorn Poonen, *Rational points on varieties,* Graduate Studies in Mathematics, 186. American Mathematical Society, Providence, RI, 2017. xv+337 pp.

[P86] Dorin Popescu, *General Néron desingularization and approximation,* Nagoya Math. J. 104 (1986), 85–115.

[P89] Dorin Popescu, *Polynomial rings and their projective modules,* Nagoya Math. J. 113 (1989), 121–128.

[P90] Dorin Popescu, *Letter to the editor: "General NÃ©ron desingularization and approximation"*, Nagoya Math. J. 118 (1990), 45–53.

[P19] Dorin Popescu, *On a question of Swan,* With an appendix by Kęstutis Česnavičius, Algebr. Geom. 6 (2019), no. 6, 716–729.

[RS16] Andreas Rosenschon, V. Srinivas, (2016), *Étale motivic cohomology and algebraic cycles,* Journal of the Institute of Mathematics of Jussieu, 15 (2016), (3): 511–537.

[R96] Markus Rost, *Chow groups with coefficients,* Doc. Math. 1 (1996), No. 16, 319–393.

[R98] Markus Rost, *Chain lemma for splitting fields of symbols,* (1998). **https://www.math.uni-bielefeld.de/~rost/chain-lemma.html**

[R02] Markus Rost, *Norm varieties and algebraic cobordism,* Proceedings of the International Congress of Mathematicians, Vol. II (Beijing, 2002), 77–85, Higher Ed. Press, Beijing, 2002.

[R15] David Rydh, *Noetherian approximation of algebraic spaces and stacks*, J. Algebra 422 (2015), 105–147.

[S20] Shuji Saito, *Purity of reciprocity sheaves,* Adv. Math. 366 (2020), 107067, 70 pp.

[S84a] David J. Saltman, *Retract rational fields and cyclic Galois extensions,* Israel J. Math. 47 (1984), no. 2–3, 165–215.

[S84b] David J. Saltman, *Noether's problem over an algebraically closed field,* Invent. Math. 77 (1984), no. 1, 71–84.

[S85] David J. Saltman, *The Brauer group and the center of generic matrices* J. Algebra 97 (1985), no. 1, 53–67.

[S19] Stefan Schreieder, *Stably irrational hypersurfaces of small slopes,* J. Amer. Math. Soc. 32 (2019), no. 4, 1171–1199.

[S21b] Stefan Schreieder, *Unramified cohomology, algebraic cycles and rationality*, 47pp, **https://arxiv.org/pdf/2106.01057.pdf**

[S97] Jean-Pierre Serre, *Galos cohomology,* Translated from the French by Patrick Ion and revised by the author., Springer 1997.

[S00] Jean-Pierre Serre, *Local algebra.* Translated from the French by CheeWhye Chin and revised by the author. Springer Monographs in Mathematics. Springer-Verlag, Berlin, 2000. xiv+128 pp.

[S72] Stephen S. Shatz, *Profinite groups, arithmetic, and geometry*, Annals of Mathematics Studies, No. 67. Princeton University Press, Princeton, N.J.; University of Tokyo Press, Tokyo, 1972. x+252 pp

[S919] A. Speiser, *Zahlentheoretische SÃ¤tze aus der Gruppentheorie,* Math. Z. 5 (1919) 1–6.

[S99] Mark Spivakovsky, *A new proof of D. Popescu's theorem on smoothing of ring homomorphisms*, J. Amer. Math. Soc. 12 (1999), no. 2, 381–444.

[sp18] The Stacks Project Authors, *Stacks Project*, **https://stacks.math.columbia.edu**, 2018.

[S17] Matt Stevenson, *Abhyankar and quasi-monomial valuations*, October 2, 2017, 9pp. **http://www-personal.umich.edu/~stevmatt/abhyankar.pdf**

[SJ06] Andrei Suslin, Seva Joukhovitski, *Norm varieties,* J. Pure Appl. Algebra 206 (2006), no. 1–2, 245–276.

[S83] Richard G. Swan, *Noether's problem in Galois theory,* Emmy Noether in Bryn Mawr (Bryn Mawr, Pa., 1982), 21–40, Springer, New York-Berlin, 1983.

[S98] Richard G. Swan, *Néron-Popescu desingularization*, Algebra and geometry (Taipei, 1995), 135–192, Lect. Algebra Geom., 2, Int. Press, Cambridge, MA, 1998.

[S21] Taku Suzuki, *Higher order minimal families of rational curves and Fano manifolds with nef Chern characters,* J. Math. Soc. Japan 73 (2021), no. 3, 949–964.

[T95] Bernard Teissier, *Résultats récents sur l'approximation des morphismes en algèbre commutative (d'après André, Artin, Popescu et Spivakovsky)*, Séminaire Bourbaki, Vol. 1993/94, Astérisque, 227, (1995), Exp. No. 784, 4, 259–282,

[T14] Bernard Teissier, *Overweight deformations of affine toric varieties and local uniformization, Valuation theory in interaction,* 474–565, EMS Ser. Congr. Rep., Eur. Math. Soc., Zürich, 2014.

[T11] Michael Temkin, *Relative Riemann-Zariski spaces* , Israel J. Math. 185 (2011), 1–42.

[T13] Michael Temkin, *Inseparable local uniformization*, J. Algebra 373 (2013), 65–119.

[TT90] R. W. Thomason, Thomas Trobaugh, *Higher algebraic $K$-theory of schemes and of derived categories,* The Grothendieck Festschrift, Vol. III, 247–435, Progr. Math., 88, Birkhäuser Boston, Boston, MA, 1990.

[T92] Burt Totaro, *Milnor $K$-theory is the simplest part of algebraic $K$-theory,* $K$-Theory 6 (1992), no. 2, 177–189.

[T99] Burt Totaro, *The Chow ring of a classifying space,* Algebraic K-theory (Seattle, WA, 1997), 249–281, Proc. Symposia in Pure Math, 67, AMS, Providence, (1999).

[T16] Burt Totaro, *Hypersurfaces that are not stably rational,* J. Amer. Math. Soc. 29 (2016), no. 3, 883–891.

[v14] *Valuation theory in interaction*, Proceedings of the Second International Conference and Workshop on Valuation Theory held in Segovia and El Escorial, July 18–29, 2011. Edited by Antonio Campillo, Franz-Viktor Kuhlmann and Bernard Teissier. EMS Series of Congress Reports. European Mathematical Society (EMS), Zürich, 2014. xiv+656 pp.

[V00] Michel Vaquié, *Valuation*, Resolution of singularities (Obergurgl, 1997), 539–590, Progr. Math., 181, Birkhäuser, Basel, 2000.

[V06] Michel Vaquié, *Valuations and local uniformization,* Singularity theory and its applications, 477–527, Adv. Stud. Pure Math., 43, Math. Soc. Japan, Tokyo, 2006.

[V98] Vladimir Voevodsky, *$\mathbb{A}^1$-homotopy theory,* Proceedings of the International Congress of Mathematicians, Vol. I (Berlin, 1998). Doc. Math. 1998, Extra Vol. I, 579–604.

[VSF00] Vladimir Voevodsky, Andrei Suslin, Eric M. Friedlander, *Cycles, transfers, and motivic homology theories,* Annals of Mathematics Studies, 143. Princeton University Press, Princeton, NJ, 2000. vi+254 pp.

[V03] Vladimir Voevodsky, *Motivic cohomology with $\mathbb{Z}/2$-coefficients,* Publ. Math. Inst. Hautes Études Sci. 98 (2003), 59–104.

[VRO07] Vladimir Voevodsky, Oliver Röndigs, Paul Arne Østvær, *Voevodsky's Nordfjordeid lectures: motivic homotopy theory,* Motivic homotopy theory, 147–221, Universitext, Springer, Berlin, 2007.

[V11] Vladimir Voevodsky, *Motivic cohomology with $\mathbb{Z}/l$-coefficients,* Annals of Math. (2) 174 (2011), 401–438.

[V14] Claire Voisin, *Chow rings, decomposition of the diagonal, and the topology of families,* Annals of Mathematics Studies, 187. Princeton University Press, Princeton, NJ, 2014. viii+163 pp.

[V19] Claire Voisin, *Birational Invariants and Decomposition of the Diagonal,* Birational Geometry of Hypersurfaces, Andreas Hochenegger, Manfred Lehn, Paolo Stellari Editors, Lecture Notes of the Unione Matematica Italiana 26, Springer 2019.

[W09] Charles A. Weibel, *The norm residue isomorphism theorem,* J. Topol. 2 (2009), 346–372.

[Z39] Oscar Zariski, *The reduction of singularities of an algebraic surface,* Ann. Math., 40 (1939),639–689.

[Z40] Oscar Zariski, *Local uniformization on algebraic varieties,* Ann. of Math. 41 (1940), 852–896.

[ZS60a] Oscar Zariski, Pierre Samuel, *Commutative algebra. Vol. I.* With the cooperation of I. S. Cohen. Corrected reprinting of the 1958 edition. Graduate Texts in Mathematics, No. 28. Springer-Verlag, New York-Heidelberg-Berlin, 1975. xi+329 pp.

[ZS60b] Oscar Zariski, Pierre Samuel, *Commutative algebra. Vol. II.* Reprint of the 1960 edition. Graduate Texts in Mathematics, Vol. 29. Springer-Verlag, New York-Heidelberg, 1975. x+414 pp.

Nagoya Institute of Technology, Gokiso, Showa-ku, Nagoya 466-8555, JAPAN
*Email address*: nori@nitech.ac.jp